\documentclass[11pt,letterpaper]{article}
\usepackage[T1]{fontenc}
\usepackage{lmodern,amsmath,amssymb,amsthm,booktabs,array,tabularx}
\usepackage[margin=1.05in]{geometry}
\usepackage{graphicx,tikz,microtype}
\usepackage[round,authoryear]{natbib}
\usepackage[hidelinks,bookmarksnumbered]{hyperref}
\hypersetup{pdftitle={Uniform FDR consistency of the Benjamini-Hochberg procedure under Gaussian dependence},pdfauthor={Zijian Wei},pdfsubject={Multiple testing under Gaussian dependence}}

\theoremstyle{plain}
\newtheorem{theorem}{Theorem}[section]
\newtheorem{proposition}[theorem]{Proposition}
\newtheorem{lemma}[theorem]{Lemma}
\newtheorem{corollary}[theorem]{Corollary}
\newcommand{\Prob}{\mathbb P}
\newcommand{\E}{\mathbb E}
\newcommand{\R}{\mathbb R}
\newcommand{\FDR}{\operatorname{FDR}}
\newcommand{\FDP}{\operatorname{FDP}}
\newcommand{\one}{\mathbf 1}
\newcommand{\Simes}{\mathsf S}
\newcommand{\G}{\gamma_d}
\newcommand{\msdoi}[2]{\href{https://doi.org/#1}{#2}}
\newcommand{\msarxiv}[2]{\href{https://arxiv.org/abs/#1}{arXiv:#2}}
\numberwithin{equation}{section}
\title{Uniform FDR consistency of the Benjamini--Hochberg procedure under Gaussian dependence}
\author{Zijian Wei\\[3pt]{\normalsize School of Mathematical Sciences, Shanghai Jiao Tong University}}
\date{November 29, 2026}
\begin{document}
\maketitle
\begin{abstract}
For Gaussian mean tests with known marginal variances, Su asked whether lowering the nominal level of the ordinary Benjamini--Hochberg procedure forces its false discovery rate (FDR) to zero uniformly over the dimension, covariance and signals. We prove a uniform upper bound of order $q\log(1/q)$ as the nominal level $q$ tends to zero. Only the true-null statistics need to be jointly Gaussian. The proof uses a maximal inequality for integrated Gaussian exceedance counts, which bounds the global-null Simes rejection probability by $e q$ for two-sided tests. When the entire vector follows a model with a fixed number of common Gaussian factors and independent residuals, we determine the sharp first-order worst-case FDR. It is asymptotic to $q\sqrt{\log(1/q)}/\sqrt{\pi}$ for one-sided tests and to half that quantity for two-sided tests. Both constants also hold for a slowly growing number of factors.
\end{abstract}


\section{Introduction}\label{sec:intro}

Can the Benjamini--Hochberg procedure be made reliable for correlated Gaussian tests simply by lowering its nominal level? For a fixed number of hypotheses the answer follows from the harmonic bound of \citet{BY}. With $m$ hypotheses and a proportion $\pi_0$ of true nulls, that bound is $\pi_0 q H_m$, where $q$ is the nominal level and $H_m=\sum_{j=1}^m j^{-1}$. If $m$ is allowed to increase as $q$ decreases, the bound need never fall below one. \citet[Definition~1.3 and p.~24]{Su} asked whether the Gaussian mean problem with known covariance has a stronger property. Writing $\mathcal C_m$ for all $m\times m$ correlation matrices, the question is whether
\begin{equation}\label{eq:question}
 \lim_{q\downarrow0}\ \sup_{m\ge1}\ \sup_{\theta\in\mathbb R^m,\ \Sigma\in\mathcal C_m}
           \FDR_{\theta,\Sigma}(q)=0 .
\end{equation}
Here the procedure uses the original critical values $qk/m$ of \citet{BH}. Its input consists of the marginal $p$-values alone. The supremum includes singular matrices and allows both the signals and their dependence to change with $q$.

We prove \eqref{eq:question} for one- and two-sided Gaussian tests and give an upper bound of order $q\log(1/q)$. A constant multiple of $q$ cannot hold uniformly. \citet{Dobriban} gave a correlated two-sided Gaussian example with FDR above the nominal level, and \citet{Lei} constructed Gaussian models for which the ratio $\FDR/q$ grows at least as fast as $\sqrt{\log(1/q)}$. Together these results describe a distinction with a practical consequence. A sufficiently small nominal level guarantees a small absolute FDR over the entire Gaussian model, although the required reduction in level can grow relative to the target error rate.

The difficulty in \eqref{eq:question} already appears when every null hypothesis is true. In that case a false discovery is a crossing of the Simes boundary at any one of $m$ ranks \citep{Simes}. Write $\widehat F$ for the empirical distribution function of the null $p$-values. At a fixed threshold $t$, $\E\widehat F(t)=t$ for boundary nulls. Summing bounds for the individual crossings would reintroduce a dependence on $m$. We instead integrate $\widehat F(t)/t$ up to a variable scale $s$. The resulting statistic has mean one and, for two-sided tests, symmetric convex sublevel sets in the underlying Gaussian space. We prove that the probability that it ever reaches $h$ is at most $1/h$, simultaneously over scales, with a bound independent of the number of coordinates and their covariance. The factor $e$ in our two-sided Simes inequality comes from a deterministic comparison of a single boundary crossing with this integral. Su's FDR-linking argument then gives \eqref{eq:question}, even when the alternative $p$-values may depend arbitrarily on the Gaussian nulls.

How much does joint Gaussianity of the signals matter? When both the nulls and the signals depend on $d$ common factors and independent Gaussian residuals, we determine the first-order worst-case FDR. For every fixed $d\ge1$, its ratio to $q\sqrt{\log(1/q)}$ tends to $1/\sqrt{\pi}$ for one-sided tests and to $1/(2\sqrt{\pi})$ for two-sided tests. \citet{Lei} established these constants for one common factor and supplied the examples attaining them. Our contribution is the matching upper bound for several factors, uniformly over their loadings and over the number of tests. The two-sided bound presents a specific obstacle. A comparison statistic may have almost zero loading along the direction of a given null and still have a substantial loading in an orthogonal direction. Section~\ref{sec:multiple-factors} shows how shifting that orthogonal component constrains the ranks that BH can select. This also explains why the same constants persist for a slowly increasing number of factors.

Two results address different restrictions on the same Gaussian experiment. \citet{ShangGuo} proved dimension-free FDR bounds for Gaussian arrays whose correlations decay at a specified inverse-logarithmic rate with index distance. The maximal inequality below does not require an ordering of the coordinates. \citet{Zhang} proved the exact two-sided Simes inequality for all covariance matrices at sufficiently small levels in each fixed dimension, and for growing dimensions under sparsity or block conditions. Zhang also established validity near $0.05$ for all three-dimensional correlations. Our Simes bound applies simultaneously to all dimensions and levels, with a constant $e$. Whether the constant can be reduced to one in this generality remains open.

\section{A maximal inequality for Gaussian exceedance counts}\label{sec:max}

Consider a centered Gaussian vector $Z$ with standard normal marginals and any correlation matrix. Write $\bar\Phi$ for the standard normal survival function. With $P_i=2\bar\Phi(|Z_i|)$, let $\widehat F(t)=m^{-1}\sum_i\mathbf1\{P_i\le t\}$. For a vector $p$ of length $l$, let
\[
 \Simes(p)=\min_{1\le k\le l}\frac{l\,p_{(k)}}{k}.
\]
The event $\{\Simes(P)\le t\}$ says that $\widehat F(kt/m)\ge k/m$ for at least one $k$. To study all these inequalities at once, integrate the empirical count on a logarithmic threshold scale:
\begin{equation}
 M(s)=\frac1s\int_0^s \frac{\widehat F(t)}{t}\,dt
     =\frac1{ms}\sum_{i=1}^m\left(\log\frac{s}{P_i}\right)_+ .
 \label{eq:integrated}
\end{equation}

Here $y_+=\max(y,0)$. The quantity $2msM(s)$ is the soft-threshold Fisher combination of \citet{TFisher}. Each $P_i$ is uniform, so $\E M(s)=1$ at every fixed $s$. The next result controls the continuum of scales at the same cost.

\begin{theorem}[Gaussian maximal inequality]\label{lem:max}
For any $m\ge1$, any Gaussian correlation matrix and any $h>0$,
\[
 \Pr\!\left\{\sup_{0<s\le1}M(s)\ge h\right\}\le\min\{1,h^{-1}\}.
\]
\end{theorem}

The proof orders Gaussian measures tilted by $M(s)$ on symmetric convex sets. The event of no previous crossing is such a set, so the order applies at the first crossing.

\begin{proof}
Represent $Z_i=a_i^\top G$, where $G$ is standard Gaussian in $\mathbb R^d$ and the $a_i$ are unit vectors. This representation also holds for singular covariance, with $d$ equal to its rank.
Put $\psi(v)=-\log(2\bar\Phi(v))$ for $v\ge0$ and $Y_i(g)=\psi(|a_i^\top g|)$. Indeed, $2\bar\Phi(|a_i^\top G|)$ is uniform on $(0,1)$, so $\Pr\{Y_i(G)>y\}=e^{-y}$ for $y\ge0$. In particular, $\E(Y_i(G)-b)_+=\int_b^\infty e^{-y}\,dy=e^{-b}$. The function $\psi$ is increasing and convex. Indeed, with $\phi$ the standard normal density and $\lambda=\phi/\bar\Phi$, we have $\psi'=\lambda$ and $\lambda'(v)=\lambda(v)\{\lambda(v)-v\}>0$, since $\phi(v)=\int_v^\infty x\phi(x)\,dx>v\bar\Phi(v)$. For $b\ge0$ define the function
\[
 H_b(g)=\frac{e^b}{m}\sum_{i=1}^m(Y_i(g)-b)_+ .
\]
It is even and convex as a function on $\R^d$, and $H_b(G)=M(e^{-b})$ by \eqref{eq:integrated}.

Throughout, symmetry of a set means $K=-K$. We first show that for every symmetric convex Borel set $K$, and for standard Gaussian measure $\G$,
\begin{equation}
 \int_K H_{b'}\,d\G\ \le\ \int_K H_b\,d\G,\qquad 0\le b\le b'.
 \label{eq:tilt}
\end{equation}
Fix one direction $a=a_i$ and write $G=a(a^\top G)+G_\perp$. The vector $G_\perp$ is an independent standard Gaussian on $a^\perp$. Define
\[
 k(v)=\Prob\{va+G_\perp\in K\},\qquad v\in\R .
\]
The symmetry of $K$ and $G_\perp$ makes $k$ even. The integrand $\one_K(va+y)\exp(-\|y\|^2/2)$ is log-concave in $(v,y)$, so Pr\'ekopa's marginal theorem makes $k$ log-concave \citep{Prekopa}. An even log-concave function is nonincreasing on $[0,\infty)$. In dimension one the same property follows directly because $K$ is a symmetric interval. Thus $\ell(y)=k(\psi^{-1}(y))$ is nonincreasing on $[0,\infty)$.

Since $k$ is even, conditioning on $a^\top G$ gives $\ell(Y_i(G))$ as the conditional probability of $K$. The exponential law of $Y_i(G)$ and the substitution $w=y-b$ yield
\[
 \begin{aligned}
 \int_K e^b(Y_i(g)-b)_+\,d\G(g)
 &= \int_b^\infty \ell(y)e^b(y-b)e^{-y}\,dy\\
 &= \int_0^\infty \ell(b+w)\,w e^{-w}\,dw .
 \end{aligned}
\]
The density $we^{-w}$ does not depend on $b$. Increasing $b$ decreases $\ell(b+w)$ at every $w$ because $\ell$ is nonincreasing. Averaging this inequality over the $m$ directions proves \eqref{eq:tilt}.

For a finite set of scales, order $b_j=-\log s_j$ increasingly, $j=1,\ldots,J$. Let $T$ be the first index $j$ with $H_{b_j}(G)\ge h$, taking $T=J$ if there is no crossing. The event that the first $j$ scores are below $h$ is $\{G\in K_j\}$, where
\[
 K_j=\bigcap_{l\le j}\{g:H_{b_l}(g)<h\}.
\]
These sets are centrally symmetric and convex by the evenness and convexity of the scores. The increment from $j-1$ to $j$ is included exactly when none of the first $j-1$ scores has crossed. Summing these increments up to the first crossing (or to $J$ if there is none) gives
\[
 H_{b_T}(G)=H_{b_1}(G)+
 \sum_{j=2}^{J}\one_{\{G\in K_{j-1}\}}\big[H_{b_j}(G)-H_{b_{j-1}}(G)\big].
\]
On the event of a crossing it is at least $h$, and it is always nonnegative. Taking expectations gives
\[
 \begin{aligned}
 h\Pr\!\left\{\max_{1\le j\le J}H_{b_j}(G)\ge h\right\}
 &\le \E H_{b_T}(G)\\
 &= 1+\sum_{j=2}^{J}\left[
      \int_{K_{j-1}}H_{b_j}\,d\G-\int_{K_{j-1}}H_{b_{j-1}}\,d\G
      \right]\le1,
 \end{aligned}
\]
where the last inequality uses \eqref{eq:tilt} on each $K_{j-1}$. The finite-scale bound follows.

By continuity of $M(s)$ for $s>0$, the supremum over $(0,1]$ equals the supremum over any fixed countable dense subset. For $0<h'<h$, a supremum at least $h$ implies a crossing of $h'$ at some point of that subset. Apply the finite-scale bound to increasing finite subsets, then let $h'\uparrow h$. Together with the bound by one, this proves the stated inequality.
\end{proof}

\subsection{A global-null Simes bound}
\begin{corollary}[Global-null Gaussian Simes bound]\label{thm:simes}
Let $Z$ be any centered jointly Gaussian vector with standard normal marginal distributions. If $P_i=2\bar\Phi(|Z_i|)$, then
\begin{equation}
 \Prob\{\Simes(P)\le t\}\le \min(1,et),\qquad 0<t<1.
 \label{eq:simes}
\end{equation}
For prespecified one-sided Gaussian $p$-values the right side can be replaced by $\min(1,2et)$.
\end{corollary}
\begin{proof}[Proof of Corollary~\ref{thm:simes}]
For two-sided tests, take $0<t<1/e$ and put $t_k=kt/m$. If $P_{(k)}\le t_k$, then $\widehat F(t_k)\ge k/m=t_k/t$. Since $\widehat F$ is nondecreasing and $et_k\le1$, on this event
\begin{equation}
 M(et_k)\ge \frac1{et_k}\int_{t_k}^{et_k}\frac{\widehat F(t_k)}{s}\,ds
 =\frac{\widehat F(t_k)}{et_k}\ge\frac1{et}.
 \label{eq:rank}
\end{equation}
Simes rejects exactly when $P_{(k)}\le t_k$ for some $k$, so
\[
 \{\Simes(P)\le t\}\subseteq
 \left\{\max_{1\le k\le m}M(et_k)\ge(et)^{-1}\right\}.
\]
Theorem~\ref{lem:max} bounds the probability on the right by $et$. For $t\ge1/e$ the probability bound by one suffices.

For right-sided values $U_i=\bar\Phi(Z_i)$, the two-sided values satisfy $P_i=2\min(U_i,1-U_i)\le2U_i$. The order statistics then give $\Simes(P)\le2\Simes(U)$, so the two-sided bound at $2t$ gives the claimed one-sided bound for $t<1/(2e)$. For $t\ge1/(2e)$ use the bound by one. A coordinate sign change accommodates a prespecified left-sided test while preserving the centered Gaussian law with unit marginals.
\end{proof}

For a dilation by $c>1$, the lower bound on $M(ct_k)$ in the proof has multiplier $(\log c)/c$, maximized at $c=e$. The factor $e$ is sharp for this deterministic smoothing step: when all the $p$-values equal some $p<1/e$, the maximum of $\widehat F(t)/t$ is $1/p$, while $\sup_s M(s)=1/(ep)$. This calculation identifies where the constant enters. Sharpness of the Gaussian rejection probability is a separate question.

\section{Uniform consistency of ordinary BH}\label{sec:main}\label{sec:linking}\label{sec:model}

Write $X=\theta+Z$, where $Z\sim N_m(0,\Sigma)$ and $\Sigma$ is a correlation matrix. Centering and scaling give this model for any known positive marginal variances, including singular joint laws. For two-sided tests of $\theta_i=0$ use $P_i=2\bar\Phi(|X_i|)$. For right-sided tests of $\theta_i\le0$ use $P_i=\bar\Phi(X_i)$. Prespecified left-sided tests are handled by changing signs. A coordinate observed with zero variance can be tested exactly and assigned $p$-value one when its null holds.

Let $I_0$ denote the set of true nulls, $m_0=|I_0|$, and $\pi_0=m_0/m$. Write $P_{(1)}\le\cdots\le P_{(m)}$ for the ordered $p$-values. At level $q\in(0,1)$ the BH rejection count is
\[
 R=\max\big(\{k:P_{(k)}\le qk/m\}\cup\{0\}\big).
\]
The procedure rejects all $P_i\le qR/m$ when $R>0$. If $V$ rejected hypotheses are true nulls, $\FDP=V/(R\vee1)$ and $\FDR=\E[\FDP]$.

Under the global null the FDP equals the indicator of at least one rejection, so the FDR is $\Pr\{\Simes(P)\le q\}$. With signals present we use the Simes statistic of the true nulls alone.
\begin{theorem}[Gaussian FDR consistency]\label{thm:main}
Suppose the true-null test statistics follow the standardized Gaussian model above. The nonnull $p$-values may have any joint distribution with these statistics. Take $C=e$ if all true-null tests are two-sided, and $C=2e$ for one-sided tests or a prespecified mixture of sidedness. Put $u=C\pi_0 q$. For $0<x\le1$,
\begin{equation}
 \Prob(\FDP\ge x)\le \min(1,u/x).
 \label{eq:fdx}
\end{equation}
Consequently,
\begin{equation}
 \FDR\le B(u),\qquad
 B(u)=\begin{cases}
 0,&u=0,\\
 u\{1+\log(1/u)\},&0<u<1,\\
 1,&u\ge1 .
 \end{cases}
 \label{eq:fdr}
\end{equation}
In particular, if $\mathcal C_m$ is the set of all $m$-dimensional correlation matrices, then for Gaussian mean testing
\begin{equation}
 \lim_{q\downarrow0}\ 
 \sup_{m\ge1}\ \sup_{\Sigma\in\mathcal C_m}\ \sup_{\theta\in\R^m}
 \FDR_{\theta,\Sigma}(q)=0 .
 \label{eq:consistent}
\end{equation}
\end{theorem}
\begin{proof}
The two-sided Simes bound \eqref{eq:simes} is the Gaussian input. Corollary~\ref{thm:simes} obtained it from the maximal inequality in Theorem~\ref{lem:max}. If $m_0=0$, the FDP is zero. Assume $m_0>0$.

For each true null with positive variance, let $Q_i=2\bar\Phi(|Z_i|)$ be its centered two-sided reference $p$-value. For a deterministic true null with $P_i=1$, adjoin an independent standard Gaussian coordinate and use its two-sided $p$-value for $Q_i$. The resulting reference vector comes from a centered Gaussian vector with unit marginal variances, so \eqref{eq:simes} applies even if its correlation matrix is singular.

For a two-sided true null with positive variance, $\theta_i=0$ gives $Q_i=P_i$. For a right-sided true null, $\theta_i\le0$ gives
\[
 \frac{Q_i}{2}=\bar\Phi(|Z_i|)\le\bar\Phi(Z_i)
 \le \bar\Phi(\theta_i+Z_i)=P_i .
\]
Changing the signs of $Z_i$ and $\theta_i$ gives the same inequality for a left-sided true null. At a deterministic true null, $Q_i\le1=P_i$. Thus $Q_i\le(C/e)P_i$ for every true null, with $C=e$ when all are two-sided and $C=2e$ otherwise.

Let $S_0=\Simes((P_i)_{i\in I_0})$. Order statistics preserve coordinatewise inequalities, and multiplying every $p$-value by the same positive constant multiplies its Simes statistic by that constant. Hence
\[
 \Simes((Q_i)_{i\in I_0})\le(C/e)S_0 .
\]
In particular, $S_0\le t$ implies that the reference Simes statistic is at most $(C/e)t$. Applying \eqref{eq:simes} at this threshold, and using the bound by one when the threshold is at least one, gives
\[
 \Prob(S_0\le t)\le\min\{1,Ct\},\qquad t>0.
\]

On any outcome with $V>0$, BH has rejected $V$ true nulls at the threshold $qR/m$. The $V$th smallest true-null $p$-value is therefore at most $qR/m$. For $x>0$, on the event $\{\FDP\ge x\}$ we have $V\ge xR>0$, so
\[
 S_0\le\frac{m_0}{V}\frac{qR}{m}
 \le\frac{\pi_0q}{x}.
\]
This is the pathwise comparison in \citet[Theorem~1]{Su}. It holds for every realization of the nonnull $p$-values. Taking probabilities and applying the true-null Simes bound gives
\[
 \Prob(\FDP\ge x)\le\Prob(S_0\le\pi_0q/x)
 \le\min\{1,C\pi_0q/x\},
\]
which is \eqref{eq:fdx}.

Since $\FDP\in[0,1]$, its mean satisfies
\[
 \FDR=\int_0^1\Prob(\FDP\ge x)\,dx
 \le \int_0^1\min\{1,u/x\}\,dx=B(u).
\]
For $0<u<1$ the integral splits at $x=u$ and equals $u+u\log(1/u)$. The cases $u=0$ and $u\ge1$ give the other two branches in \eqref{eq:fdr}. Finally, $B$ is nondecreasing and $u=C\pi_0q\le2eq$ for every model. Thus
\[
 \sup_{m\ge1}\ \sup_{\Sigma\in\mathcal C_m}\ \sup_{\theta\in\R^m}
 \FDR_{\theta,\Sigma}(q)\le B(2eq)\longrightarrow0
 \qquad(q\downarrow0).
\]
\end{proof}
\subsection{What the logarithmic factor measures}

Combining Theorem~\ref{thm:main} with the harmonic bound gives $\min\{\pi_0qH_m,B(C\pi_0q)\}$ at every finite $m$. It applies without a decay envelope or a constraint on the covariance rank. Under the inverse-logarithmic decay assumption of \citet{ShangGuo}, a comparison with their explicit finite-sample bound and a sharp asymptotic analysis are given in Appendix~\ref{S:sec:decay-sharp}.

It is instructive to locate the logarithm in the proof. On an outcome with $V>0$ and at most $s$ alternatives, $\FDP\ge1/(s+1)$. Hence the portion of the tail integral that yields $\log(1/q)$ is truncated once $s$ is fixed. Conversely, even with independent nulls, signals allowed to depend on the nulls can produce a logarithm. For $n$ boundary nulls set $a=nq/(n+s)$ and $K=\lfloor\min(s+1,1/a)\rfloor$. If the smallest null $p$-value falls in $((j-1)q/m,jq/m]$ with $1\le j\le K$, assign $j-1$ signal $p$-values zero and the rest one. There are at least $j$ rejections, including at least one true null. Additional rejections are true nulls and can only increase the resulting FDP above $1/j$. The interval has probability at least $a/e$, because its upper endpoint is at most $1/n$ and the minimum has density $n(1-t)^{n-1}$. Summing the contributions gives at least $(a/e)H_K$. This is the finite-count form of the construction of \citet{Su}, also studied by \citet{LyuWang}. Appendix~\ref{S:app:finite} gives the matching order for all $n,s$ and the leading coefficient in several regimes.

This construction chooses the signals after seeing the smallest null $p$-value. In a Gaussian mean model each signal instead consists of a deterministic mean and a Gaussian noise coordinate. The difference is visible even for two tests. With one two-sided Gaussian null and one unrestricted alternative, the worst FDR is exactly $3q/4$. When the alternative must be a joint Gaussian mean shift, the corresponding supremum is asymptotic to $5q/8$ (Proposition~\ref{S:prop:one-null}). We next quantify what this restriction buys in models with common Gaussian factors.

\section{Sharp FDR under Gaussian factor dependence}\label{sec:multiple-factors}\label{sec:sharp}

The bound in Section~\ref{sec:main} allows the alternative $p$-values to be arbitrary functions of the nulls. We now return to Gaussian mean shifts. For $d\ge1$, consider
\begin{equation}\label{eq:factor-d-model}
 T_i=\theta_i+b_i^\top G+\sqrt{1-\|b_i\|^2}\,\varepsilon_i,\qquad
 b_i\in\mathbb R^d,\quad \|b_i\|\le1 ,
\end{equation}
where $G\sim N_d(0,I_d)$, and the $\varepsilon_i$ are independent standard normals independent of $G$. Loadings may have either sign and residual variances may be zero. The covariance of $(T_1,\ldots,T_m)$ can have full rank for any $d$. Let $B_{d,b}(q)$ be the supremum of the BH FDR over $m$, $\theta$ and all loadings in \eqref{eq:factor-d-model}. We use $b=1$ for right-sided tests and $b=2$ for two-sided tests. Set $Q=\bar\Phi$ and $u_b=Q^{-1}(q/b)$.

\begin{theorem}[Sharp factor-model bounds]\label{thm:multifactor-sharp}
For fixed $d\ge1$, as $q\downarrow0$,
\begin{align}
 \FDR_{\rm one}&\le\pi_0q\,\frac{u_1}{\sqrt{2\pi}}+O_d(q),\label{eq:one-fdr}\\
 \FDR_{\rm two}&\le\pi_0q\,\frac{u_2}{2\sqrt{2\pi}}+O_d(q\sqrt{u_2}).\label{eq:two-fdr}
\end{align}
The remainders are uniform in $m$, $\theta$ and the loadings. In particular, for $b=1,2$,
\begin{equation}\label{eq:factor-sharp}
 \lim_{q\downarrow0}\frac{B_{d,b}(q)}{q\sqrt{\log(1/q)}}=\frac{1}{b\sqrt{\pi}} .
\end{equation}
\end{theorem}

The constants in \eqref{eq:factor-sharp} are attained asymptotically by the one-factor sequences of \citet[Theorems~2.1 and~3.1 and Propositions~2.2 and~3.2]{Lei}. Embedding their loading vectors in the first coordinate gives the lower bounds for every $d$. The work in the theorem lies in the upper bounds: with arbitrarily many loading directions, independent residuals and unrestricted means, an additional factor cannot increase either leading coefficient.

The following finite model illustrates the difference between consistency and exact nominal control. With $m=200$, take $T_i=0.3G+\sqrt{0.91}\,\varepsilon_i$ for $i\le199$ and $T_{200}=5.5-G$, with all inputs independent standard normals. The nulls have pairwise correlation $0.09$ and correlation $-0.3$ with the signal. At $q=0.05$, numerical integration gives $\FDR/q=0.98785$ for one-sided tests and $0.99403$ for two-sided tests. At $q=0.005$ the corresponding ratios are $1.01096$ and $1.00311$ for the same mean and covariance. Appendix~\ref{S:app:quadrature} gives the numerical integration for this example.

\subsection{Conditioning on the common factors}\label{app:factors}

For a true null $i$ put $r_i=\|b_i\|$. For $r_i>0$, write $v_i=b_i/r_i$ and $Y_i=v_i^\top G$. At the boundary of the null, the conditional distribution of $T_i$ given $G$ is that of $r_iY_i+\sqrt{1-r_i^2}\,\varepsilon_i$. Denote the corresponding conditional distribution function of $P_i$, given $Y_i=y$, by $F^{(b)}_{r_i,y}$, and set
\[
 D_r(y)=\sup_{0<t\le q}\frac{F^{(b)}_{r,y}(t)}{t}.
\]
For $r_i=0$ take any unit $v_i$, so $D_0=1$. A negative null mean in a right-sided test only reduces its conditional distribution function.

Let $R_i$ be the rejection count after replacing $P_i$ by zero. Once $G$ is fixed, $R_i$ depends only on the residuals of the other tests, which are independent of the residual for $i$. Thus $R_i$ and $P_i$ are conditionally independent. The procedure rejects $i$ precisely when $P_i\le qR_i/m$, in which case $R=R_i$. Since $R_i\ge1$, the conditional expected contribution of null $i$ to the FDP, given $G$ and $R_i$, is bounded by $F^{(b)}_{r_i,Y_i}(qR_i/m)/R_i\le(q/m)D_{r_i}(Y_i)$. Summing on an event $A$ determined by $G$ and using $\FDP\le1$ on $A^c$ gives
\begin{equation}\label{eq:factor-cutoff}
 \FDR\le \Pr(A^c)+\frac qm\sum_{i\in I_0}\E[D_{r_i}(Y_i)\mathbf1_A].
\end{equation}
The same event $A$ truncates every null projection at a total FDR cost of $\Pr(A^c)$.

\begin{proposition}[A scalar bound for either choice of tails]\label{thm:fixed-factors}
For fixed $d$, with $u=u_b$, both choices of tails satisfy
\begin{equation}\label{eq:factors-main}
  \FDR\le\pi_0 q\,\frac{u}{\sqrt{2\pi}}+O_d(q).
\end{equation}
If $\|b_i\|=1$ for all true nulls, the FDR is $O_d(q[1+\log\log(1/q)])$.
\end{proposition}

\begin{proof}
Write $\phi$ for the standard normal density. The two scalar estimates needed below, for $u\ge1$, are
\begin{align}
 \int_{-u}^u \phi(y)D_r(y)\,dy
   &\le \frac12+\sqrt{\frac{1-r}{1+r}}\,\frac{u}{\sqrt{2\pi}}\quad (r<1),\label{eq:factor-strip}\\
 D_r(y)&\le 1/Q(|y|)\quad (0\le r\le1).\label{eq:factor-outer}
\end{align}
The central integral is zero at $r=1$. To see the first inequality, put $s=\sqrt{1-r^2}$ and exchange $x$ and $y$ in the bivariate normal density. For a right-sided threshold $x\ge u$, the product of $\phi(y)$ and the conditional tail ratio equals the average of $s^{-1}\phi((y-rt)/s)$ over a standard normal $t$ restricted to $t\ge x$. Its supremum is at most $s^{-1}\phi((ru-y)_+/s)$. Integrating over $y\le u$ gives \eqref{eq:factor-strip}. For two tails average the bounds at $y$ and $-y$ over $[-u,u]$ \citep[Lemma~4.3]{Lei}. This argument also bounds $\phi(y)D_r(y)$ by $1/(s\sqrt{2\pi})$.

For \eqref{eq:factor-outer} write $z=|y|$. A threshold $x\le z$ satisfies the bound directly. If $x>z$, completing a square gives $(x-rz)^2/(1-r^2)\ge x^2-z^2$. For $z>0$, $Q(\sqrt{x^2-z^2})/Q(x)$ decreases in $x>z$ and its limit at $z$ is $1/(2Q(z))$. Indeed, the Gaussian hazard $\lambda=\phi/Q$ satisfies $(\lambda(v)/v)'<0$ by the Mills bound $\lambda(v)<v+1/v$. The case $z=0$ is immediate, and at $r=1$ the conditional law is a point mass.

Take $v$ to be the larger of $u$ and the upper $q$-quantile of $\|G\|$. Gaussian and chi tail estimates at fixed $d$ give
\begin{equation}\label{eq:factor-quantiles}
 v-u=O_d\big((1+\log u)/u\big),\qquad
 \log\{Q(u)/Q(v)\}=O_d(1+\log u).
\end{equation}
In \eqref{eq:factor-cutoff} take $A=\{\|G\|\le v\}$. Each projection then lies in $[-v,v]$. For $r\le1/2$ the uniform bound on $\phi D_r$ makes the integral over $u<|y|\le v$ of order $O_d(1)$. For $r>1/2$, \eqref{eq:factor-outer} bounds this integral by $2\log\{Q(u)/Q(v)\}$, while the coefficient of $u/\sqrt{2\pi}$ in \eqref{eq:factor-strip} is at most $1/\sqrt3$. Its deficit absorbs the shell contribution as $u$ grows. When all $r_i=1$, only the shell and $\Pr(A^c)\le q$ contribute.
\end{proof}

\subsection{Why the two-sided coefficient is half as large}\label{sec:multifactor-sharp}

In \eqref{eq:factor-cutoff} the supremum over $t$ allows the rejection threshold to depend arbitrarily on the common factor. For small $r$ the scalar bound retains a coefficient twice the one in \eqref{eq:two-fdr}. The improvement comes from a constraint on the threshold itself: if BH rejects $k$ hypotheses, every larger rank must have failed. We describe how that constraint limits the part of the Gaussian sample space on which a near-maximal conditional ratio can be realized.

Write $p(x)=2Q(x)$, $u=Q^{-1}(q/2)$, $t_k=qk/m$ and $x_k=p^{-1}(t_k)$. Fix a null $i$ with $0<r=\|b_i\|<3/5$ and decompose $G=zv_i+w$ into independent standard Gaussians on the line and its orthogonal complement. The other coordinates have the representation
\[
 T_j(z,w)=\theta_j+a_jz+c_j^\top w+s_j\varepsilon_j,\qquad
 a_j=v_i^\top b_j,\quad c_j=b_j-a_jv_i,\quad s_j^2=1-\|b_j\|^2.
\]
The mean exceedance count and the null tail ratio are
\[
\begin{aligned}
 M_{z,w}(x)&=m^{-1}\sum_{j\ne i}\Pr(|T_j(z,w)|\ge x\mid z,w),\\
 \rho_r(x,z)&=\Pr(|rz+\sqrt{1-r^2}\,\varepsilon_i|\ge x)/p(x).
\end{aligned}
\]
For a given $z,w$, call $x_k$ admissible if $k>J=\lceil\sqrt u\rceil$ and
\begin{equation}\label{eq:mf-feasible}
 M_{z,w}(x_k)\ge\frac{p(x_k)}{2q},\qquad
 M_{z,w}(y)\le\frac{3p(y)}q\quad(u\le y<x_k).
\end{equation}
Given $G$, the exceedance counts are sums of independent indicators. On $R_i=k$ the count of other coordinates at $x_k$ is $k-1$ and the count at $x_l$ is at most $l-2$ for $l>k$. Summing Chernoff bounds over the ranks shows that $R_i>J$ selects an admissible threshold except on an event of conditional probability at most $Ce^{-cJ}$. The upper inequality extends from the grid to every intermediate $y$ by monotonicity. Each rank $k\le J$ costs at most $q/m$ after integrating in $z$, since $\int\phi(z)\rho_r(x_k,z)\,dz=1$.

Put $K=\sqrt{8\log u}$, $L=K+\sqrt u$. For fixed $w$, let $\mathcal W_r(w)$ be the set of $z$ with $L\le|z|\le u$ for which some admissible $x$ satisfies $\big||z|-rx\big|\le K$. Write $|\mathcal W_r(w)|$ for its Lebesgue measure. The scalar Gaussian tail estimates of \citet[Lemma~4.4]{Lei} concentrate the integral of $\phi(z)\sup_x\rho_r(x,z)$ on these sets when the supremum is over admissible $x$. The integral omitted on $[-u,u]$ is $O(L+ue^{-K^2/2})$. At a contributing pair $(z,x)$ the integrand $\phi(z)\rho_r(x,z)$ is at most
\begin{equation}\label{eq:mf-height}
 \frac{1+C(K+1)/u}{2\sqrt{2\pi}\sqrt{1-r^2}}+Ce^{-u}.
\end{equation}
The following estimate controls how often this height can be attained.

\begin{lemma}[Length of the admissible set]\label{lem:length}
For each fixed $d$, uniformly over $0<r<3/5$, the Gaussian means, the loadings and $m$,
\[
 \mathbb E_w|\mathcal W_r(w)|\le (1-r)u+O_d(\sqrt u).
\]
\end{lemma}

Here is the geometric step in the lemma. Split the range of $|z|$ into dyadic intervals $[a,\min(2a,u)]$ starting at $L$. On such an interval let $\eta=a^{-1/2}$. At an admissible pair, at least a third of the required mass in \eqref{eq:mf-feasible} comes from one of three groups: coordinates with $|a_j|\ge\eta r/a$, those with smaller $|a_j|$ and total loading at most $r/20$, and those with smaller $|a_j|$ and larger total loading. In the first group, moving $z$ in one of its two directions increases the absolute values of a fixed mass of coordinates. Condition \eqref{eq:mf-feasible} then precludes simultaneous witnesses of both signs at larger magnitudes, apart from an interval of length $O(K+\sqrt a)$. The small-total-loading group can contribute only a further $O(K+1)$ length. A convexity bound on its folded Gaussian density ratio gives this estimate. Appendix~\ref{S:app:admissible-proof} proves the bounds for all three groups and their sum in Lemma~\ref{lem:length}.

The third group is the one that appears when there is more than one factor. Its coordinates have $\|c_j\|\ge r/30$ for large $u$. Among the $2(d-1)$ signed basis directions in $v_i^\perp$ one can therefore choose $e$ so that a fixed fraction of the relevant mass increases in absolute value when $w$ moves along $e$. For this direction write $E_e(w)$ for the part of the set supported by that fraction on this interval. For $t\ge h=C_d\eta$, shifting $w$ to $w+te$ increases the corresponding absolute values by at least $\eta r$ even after $z$ is changed to any other point with $|z|$ in the interval. Since $p(x+\eta r)\le e^{-x\eta r}p(x)$ and $x\eta r\ge\eta(a-K)$, the shifted mass exceeds the upper threshold in \eqref{eq:mf-feasible} at $x+\eta r$. Thus, whenever both sets are nonempty,
\begin{equation}\label{eq:mf-shift}
 \sup_{z'\in E_e(w+te)}|z'|\le \inf_{z\in E_e(w)}|z|+2K+\eta r^2 .
\end{equation}

To see how \eqref{eq:mf-shift} yields a bound after averaging, write $w=w_0+se$ and, for fixed $w_0$, let $H(s)$ be the supremum of $|z'|$ over $z'\in E_e(w_0+s'e)$ with $s'\ge s$, using $a$ if none are occupied. This nonincreasing function takes values in $[a,2a]$. Up to endpoints the set at $s$ has length at most $2[H(s)-H(s+h)+2K+\eta r^2]$. Its Gaussian average is $O_d(ha+K)$, since
\[
 \int_\mathbb R\phi(s)[H(s)-H(s+h)]\,ds
 \le h\|\phi\|_\infty a
\]
by integration against the variation measure of $H$. There are only $2(d-1)$ directions at fixed $d$. Their total contribution on this interval is $O_d(\sqrt a+K)$. The first group contributes the length of the interval intersected with $[ru-K,u]$ once. Summing over $a$ gives the leading length $(1-r)u$ and an error $O_d(\sqrt u)$, since $\sum_a\sqrt a=O(\sqrt u)$ and $K\log u=O(\sqrt u)$.

\begin{proof}[Proof of Theorem~\ref{thm:multifactor-sharp}]
Proposition~\ref{thm:fixed-factors} gives \eqref{eq:one-fdr}. For two-sided tests with $0<r<3/5$, apply the conditional rank bound \eqref{S:eq:mf-conditional} from Appendix~\ref{S:app:admissible-proof} to one true null on $|z|\le u$. Its first two terms account for the first $J$ ranks and the exceptional event. The integral outside the witnesses costs $O(L+ue^{-K^2/2})$. On the witnesses, multiply the height bound \eqref{eq:mf-height} by the expected length bound in Lemma~\ref{lem:length}. In units of $q/m$, the central contribution is then at most
\[
 \begin{aligned}
 &J+O\!\left((1+u)e^{-cJ}+L+ue^{-K^2/2}\right)\\
 &\quad+\left\{\frac{1+C(K+1)/u}{2\sqrt{2\pi}\sqrt{1-r^2}}+Ce^{-u}\right\}
       \bigl\{(1-r)u+O_d(\sqrt u)\bigr\}.
 \end{aligned}
\]
The first line is $O(\sqrt u)$. Since $(1-r)/\sqrt{1-r^2}\le1$ and $(1-r^2)^{-1/2}\le5/4$, the second line is at most $u/(2\sqrt{2\pi})+O_d(\sqrt u)$.

For $3/5\le r<1$, \eqref{eq:factor-strip} already suffices: $\sqrt{(1-r)/(1+r)}\le1/2$. For $r=1$ there is no central contribution, and a null with $r=0$ costs at most $q/m$ altogether. The radial cutoff in \eqref{eq:factor-quantiles} costs $q$ once for the entire FDP. Its shell adds at most $O_d(1+\log u)$ per null in units of $q/m$. Summation proves \eqref{eq:two-fdr}. Since $u_b=\sqrt{2\log(1/q)}+o(1)$, Lei's one-factor examples give the matching coefficients. On the one-sided class their $O(q)$ lower remainder also gives $B_{d,1}(q)=q\sqrt{\log(1/q)}/\sqrt\pi+O_d(q)$.
\end{proof}

\begin{corollary}[A growing number of factors]\label{cor:mf-growing-sharp}
If $d=d(q)\ge1$ and $(d+2)\log u_2=o(u_2)$, both limits in \eqref{eq:factor-sharp} hold along this sequence.
\end{corollary}

Appendix~\ref{S:app:growing-two-sided} proves the two-sided assertion by replacing the finite collection of transverse directions with a random unit vector. Averaging over that direction bounds the error in \eqref{eq:two-fdr} by a constant times $q\{d^{1/4}\sqrt{u_2}+(d+2)\log u_2\}$. The second summand comes from the Gaussian radial cutoff. Appendix~\ref{S:sec:growing-one-sided} applies the cutoff in the scalar proof to obtain the one-sided assertion. In one dimension the transverse group is empty. Corollary~\ref{S:cor:sharp-cone-remainder} in Appendix~\ref{S:app:benefit-factors} uses this fact to bound the two-sided upper remainder by $O(q[1+\log\log(1/q)]^{3/2})$ after the leading term is written as $q\sqrt{\log(1/q)}/(2\sqrt\pi)$.

\section{Discussion}\label{sec:discussion}

Two bounds now delimit the worst-case FDR for a fully Gaussian vector without structural assumptions. With $L=\log(1/q)$, Lei's examples give a lower bound of order $q\sqrt L$ and Theorem~\ref{thm:main} gives an upper bound of order $qL$. Theorem~\ref{thm:multifactor-sharp} identifies the lower scale, including its constant, in every fixed-dimensional factor model with independent residuals. Removing the restriction on factor dimension would require control of how many different null directions can influence the BH rank. For example, take $K$ independent bivariate blocks with the same nonzero within-block correlation. The cross-covariance between the first and second coordinates of the blocks is a nonsingular $K\times K$ diagonal matrix. In any factor representation with independent residuals this submatrix is the product of one loading matrix and the transpose of another, and has rank at most $d$. Thus even these sparse Gaussian models can require $d\ge K$.

A restriction on the rank of the null covariance itself has a different effect. With rank at most a fixed $d\ge2$ and arbitrary alternatives, the worst FDR over all dimensions is asymptotic to $\tfrac{d-1}{2}q\log\log(1/q)$. If instead the ordered null correlations decay at least as fast as $[\log(e+|i-j|)]^{-1-\epsilon}$ for fixed $\epsilon>0$ and a fixed envelope constant, the worst FDR with $n$ true nulls and $s$ unrestricted alternatives is asymptotic to $a H_{\lfloor\min(s+1,1/a)\rfloor}$, uniformly over $n,s$, with $a=nq/(n+s)$. Independent nulls and signals chosen as functions of them attain the latter coefficient. Appendices~\ref{S:sec:purerank} and~\ref{S:sec:decay-sharp} give the proofs.

The unit-constant two-sided Simes problem is a second open question. Even an affirmative answer in full generality would, by Su's linking theorem alone, still give an FDR upper bound of order $q\log(1/q)$ under arbitrary alternatives. The rate question asks how closely the signals can adapt to the null order statistics while the whole vector remains Gaussian. The Simes question asks how often Gaussian nulls alone can cross several rank thresholds. At each fixed dimension $m\ge2$ there is a neighborhood of independence on which the two-sided Simes size is strictly below $q$ at all levels for every nonidentity correlation matrix (Theorem~\ref{S:thm:local-all-level}). The size at independence is exactly $q$.

\section*{Computational materials}
The ancillary R script reproduces the finite Gaussian examples in Appendix~\ref{S:sec:numerics} using the rank recursions in Appendix~\ref{S:app:quadrature}. During preparation, the finite-scale maximal argument, the Gaussian Simes bounds in Corollary~\ref{thm:simes}, and the BH bounds in Theorem~\ref{thm:main} were also formalized in Lean.

\section*{Disclosure of generative AI use}\hypertarget{ai-disclosure}{}
The author used OpenAI's GPT-6.2 Astra through ChatGPT and Codex in the research and preparation of this preprint. Its contributions included developing and revising the mathematical arguments for the Gaussian maximal inequality, the factor-model bounds, and the results in the appendices. It also assisted in searching the literature, writing and checking the Lean formalization, preparing and running the R code, and drafting and editing the text. The scope of the formal verification is described above.

\appendix
\renewcommand{\thesection}{S\arabic{section}}
\renewcommand{\theHsection}{supp.\arabic{section}}
\setcounter{figure}{0}\setcounter{table}{0}
\renewcommand{\thefigure}{S\arabic{figure}}
\renewcommand{\thetable}{S\arabic{table}}
\renewcommand{\theHfigure}{supp.\arabic{figure}}
\renewcommand{\theHtable}{supp.\arabic{table}}
\section*{Technical appendices}
\pdfbookmark[1]{Technical appendices}{supplement-start}
\noindent\hyperref[S:app:benefit-factors]{\textbf{\ref*{S:app:benefit-factors}}\quad Proofs and extensions for Gaussian factors}\par
\noindent\hyperref[S:app:finite]{\textbf{\ref*{S:app:finite}}\quad The number of alternatives and the dimension}\par
\noindent\hyperref[S:sec:purerank]{\textbf{\ref*{S:sec:purerank}}\quad A bound on the rank of the null covariance}\par
\noindent\hyperref[S:sec:decay-sharp]{\textbf{\ref*{S:sec:decay-sharp}}\quad Sharp FDR with correlations decaying in the index distance}\par
\noindent\hyperref[S:sec:planar-main]{\textbf{\ref*{S:sec:planar-main}}\quad Exact control for two-sided Gaussian Simes tests}\par
\noindent\hyperref[S:sec:numerics]{\textbf{\ref*{S:sec:numerics}}\quad Finite Gaussian examples}\par
\noindent\hyperref[S:app:more]{\textbf{\ref*{S:app:more}}\quad Further properties of the maximal inequality}\par
\noindent\hyperref[S:app:ext]{\textbf{\ref*{S:app:ext}}\quad Weights and Gaussian subspace tests}\par
\bigskip
\section{Proofs and extensions for Gaussian factors}\label{S:app:benefit-factors}
\subsection{The admissible-set length and the sharp two-sided bound}\label{S:app:admissible-proof}

Starting from the conditional BH cutoff in \eqref{eq:factor-cutoff}, we restrict attention to ranks that the other tests can support. This reduction gives \eqref{S:eq:mf-conditional}. The geometric estimates then give the expected length bound used to control the remaining integral.
\begin{proof}[Proof of Lemma~\ref{lem:length} and the two-sided bound in Theorem~\ref{thm:multifactor-sharp}]
Write $p(x)=2\bar\Phi(x)$ and $u=\bar\Phi^{-1}(q/2)$. Put $t_k=qk/m$, $x_k=p^{-1}(t_k)$, and $J=\lceil\sqrt u\rceil$. We work for sufficiently large $u$ depending at most on $d$. Fix a true null $i$ with $0<r=\|b_i\|<3/5$. Decompose $G=zv_i+w$, $v_i=b_i/r$, so that $z$ and $w$ are independent standard Gaussians on the line and its orthogonal complement. Conditional on $(z,w)$ a comparison statistic has the form
\[
 Y_j(z,w)=\theta_j+a_j z+c_j^\top w+s_j\varepsilon_j,\quad a_j=v_i^\top b_j,\quad c_j=b_j-a_jv_i,\quad s_j^2=1-\|b_j\|^2 .
\]
Let $M_{z,w}(x)=m^{-1}\sum_{j\ne i}\Pr\{|Y_j(z,w)|\ge x\}$, where the probability integrates only the residual. We regard this as a measure of total mass at most one on affine paths, putting mass $1/m$ on each index and sampling its residual. Let $\rho_r(x,z)=\Pr\{|rz+\sqrt{1-r^2}\varepsilon_i|\ge x\}/p(x)$.

\paragraph{Admissible ranks.} Call $x_k$, $k>J$, admissible at $(z,w)$ if
\begin{equation}\label{S:eq:mf-feasible}
 M_{z,w}(x_k)\ge p(x_k)/(2q),\qquad M_{z,w}(y)\le3p(y)/q\quad (u\le y<x_k).
\end{equation}
Pin $P_i$ at zero and let $R_i$ be the resulting BH count. Outside an event of conditional probability $Ce^{-cJ}$, $R_i=k>J$ implies admissibility. Indeed, the number $N_l$ of comparison statistics exceeding $x_l$ in absolute value is Poisson-binomial with mean $mM_{z,w}(x_l)$. Chernoff bounds and a sum over $l\ge J$ exclude simultaneously $N_l\ge l-1$ with mean $<l/2$, and $N_l\le l-2$ with mean $>2l$, except on the stated event. On $R_i=k$, $N_k=k-1$ and $N_l\le l-2$ for $l>k$. Interpolation gives the second inequality: for $x_l\le y<x_{l-1}$, $l>k$, the bound $M_{z,w}(y)\le2l/m$ is at most $3p(y)/q$ since $l\ge3$.

Relative to the common factor cutoff in \eqref{eq:factor-cutoff}, the central part $|z|\le u$ of the $i$th contribution is at most $q/m$ times
\begin{equation}\label{S:eq:mf-conditional}
 J+O((1+u)e^{-cJ})+
 \mathbb E_w\int_{-u}^u\phi(z)\sup_{x\text{ admissible at }(z,w)}\rho_r(x,z)\,dz .
\end{equation}
In fact each rank $k\le J$ costs at most one since $\int\phi(z)\rho_r(x_k,z)\,dz=1$. The exceptional-event bound follows by multiplying its conditional probability by the unrestricted scalar envelope and using \eqref{eq:factor-strip}. An empty supremum equals zero. Set
\[
 K=\sqrt{8\log u},\qquad L=K+\sqrt u.
\]
Call $(z,w,x)$ a witness if $x$ is admissible, $L\le|z|\le u$ and $\big||z|-rx\big|\le K$. By \citet[Lemma 4.4, (205)--(206)]{Lei}, the parts of the integral with $|z|<L$ or $\big||z|-rx\big|>K$ cost $O(L+u e^{-K^2/2})$. On witnesses its integrand $\phi(z)\rho_r(x,z)$ is at most
\begin{equation}\label{S:eq:mf-height}
 \frac{1+C(K+1)/u}{2\sqrt{2\pi}\sqrt{1-r^2}}+Ce^{-u}.
\end{equation}
Here $L=K+\sqrt u$ is the lower endpoint in the scalar estimate. It suffices to bound the expected length in $z$ of the witness set.

Split $[L,u]$ into dyadic intervals $D_a=[a,\min(2a,u)]$, $a=2^j L$, retaining only $a<u$. Endpoints have zero length and may be counted twice. On each such interval take $\eta=a^{-1/2}$ and $b=\eta r/a$. Partition the comparison indices into three classes:
\[
 A=\{j:|a_j|\ge b\},\quad S=\{j:|a_j|<b,\ \|b_j\|\le r/20\},\quad T=\{j:|a_j|<b,\ \|b_j\|>r/20\}.
\]
At every witness some class contributes at least $p(x)/(6q)$ to $M_{z,w}(x)$. Denote the three possible witness sets restricted to $|z|\in D_a$ by $E_A(w),E_S(w),E_T(w)$ respectively, allowing ties. For every $x>0$ and $t\ge0$ the Gaussian hazard lower bound gives $p(x+t)\le e^{-xt}p(x)$.

\paragraph{Large projected loadings.} Fix $w$. Any witness in $E_A(w)$, written $(z_0,x_0)$ with $y_0=|z_0|$, has a mass at least $p(x_0)/(12q)$ of affine paths exceeding $x_0$ whose absolute value grows when $z$ moves in the same direction $\sigma\in\{-1,1\}$. Along that direction for distance $t$ their absolute values increase by at least $bt$. Put
\[
 D=2K+4+\frac{10a}{\eta(a-K)}.
\]
For any $y\in D_a$ with $y>y_0+D$, the two points $\pm y$ cannot both have witnesses. Indeed, one of them lies in the favored direction from $z_0$, at distance $t\in[y-y_0,4a]$. If its witness has threshold $x'$, then $r(x'-x_0)\ge y-y_0-2K>4$ whereas $rbt\le4\eta r^2<2$, so $x_0+bt<x'$. Also $x_0bt\ge\eta(a-K)(y-y_0)/a>10$. The favored mass at the new point therefore exceeds $3p(x_0+bt)/q$, contradicting \eqref{S:eq:mf-feasible}. Taking $y_0$ down to the infimum of magnitudes of $E_A(w)$, both signs can contribute only over an initial interval of length at most $D$. All magnitudes belong to $D_a\cap[ru-K,u]$. Thus
\begin{equation}\label{S:eq:mf-axial}
 |E_A(w)|\le |D_a\cap[ru-K,u]|+D=|D_a\cap[ru-K,u]|+O(K+\sqrt a).
\end{equation}

\paragraph{Small total loadings.} Fix $w$ and a witness $(z_0,x_0)$ in $E_S(w)$. Use only the submixture indexed by $S$; denote its tail at $z_0$ by $M^S(x)$ and the ratio $M^S(x)/p(x)$ by $Q(x)$. We claim there is $v\in[x_0,x_0(1+C/u)]$ with $Q(v)\ge e^{10u}/q$. The following verification supplies the needed two-parameter version of the small-coefficient argument in \citet[Lemma 4.6]{Lei}.

Let $\kappa(x)$ be the density of the folded submixture at $x\ge0$ divided by $2\phi(x)$. Put $B=r/20$, $\beta=B^2/(1-B^2)$. Each component, before multiplying by its mass, equals
\[
 s_j^{-1}\exp\{-\mu_j^2/(2s_j^2)-(s_j^{-2}-1)x^2/2\}\cosh(\mu_j x/s_j^2),\qquad \mu_j=\theta_j+a_jz_0+c_j^\top w.
\]
Since $s_j^2\ge1-B^2$, each component multiplied by $\exp(\beta x^2/2)$ is log-convex. H\"older's inequality preserves log-convexity under summation, and $\kappa(0)\le(1-B^2)^{-1/2}$. Thus $\log\kappa(x)+\beta x^2/2$ is convex. The ratio $Q$ is the average of $\kappa$ against the Gaussian tail above $x$, and $Q'= (\phi/\bar\Phi)(Q-\kappa)$. Put $R=12e^{10u}$. Until $Q$ reaches $RQ(x_0)$ or $\kappa$ reaches $Q/2$, $(\log Q)'\ge x_0/2$. Consequently one of these events occurs by some $y\le x_0+2\log R/x_0\le x_0(1+C/u)$, with $Q(y)\ge Q(x_0)\ge 1/(6q)$. The first event proves the claim. In the second case $H=\log[\kappa(y)/\kappa(0)]\ge u^2/3$ for large $u$. On $y\le x\le2y$ convexity gives
\[
 (\log\kappa)'(x)\ge H/y-\beta x\ge H/(2y),
\]
because $rx_0\le u+K$ and $2\beta y^2\le 2(u+K)^2(1+C/u)^2/399\le H/2$. Take $v=y(1+100/u)$. Then $v+1/v\le2y$ for large $u$, and on $[v,v+1/v]$ the density ratio is at least $e^{100H/(2u)}\kappa(y)$. This interval captures at least $1-e^{-1}$ of the Gaussian tail above $v$ by the hazard bound. The claim follows from $\kappa(y)\ge1/(12q)$.

For arbitrary $|z'|\in D_a$, couple each comparison path with its value at $z_0$. Its change on $S$ has magnitude at most $\Delta=4ba\le2r$. Thus $M_{z',w}(v-\Delta)\ge p(v)Q(v)$. Since $rv\le(u+K)(1+C/u)$ and $v-\Delta\ge u-2$, the hazard upper bound $\phi(x)/\bar\Phi(x)\le x+1/x$ gives $p(v)/p(v-\Delta)\ge e^{-4u}$. Hence $M_{z',w}(v-\Delta)>3p(v-\Delta)/q$. This also forces $v-\Delta\ge u$ since $M_{z',w}\le1$. By \eqref{S:eq:mf-feasible} every admissible threshold there is at most $v-\Delta$. Its witnesses have $|z'|\le rv+K\le |z_0|+2K+C$. Taking the infimum of $|z_0|$ in $E_S(w)$ bounds $|E_S(w)|\le 4K+C$.

\paragraph{Strong transverse loadings.} This class is empty when $d=1$. For $d\ge2$, take an orthonormal basis for $v_i^\perp$ and let $\mathcal E$ be the set of its $2(d-1)$ signed coordinate directions. Every $j\in T$ satisfies $\|c_j\|^2=\|b_j\|^2-a_j^2\ge r^2/900$ for sufficiently large $u$. For each such path exceeding $x$ there is an $e\in\mathcal E$ for which $\operatorname{sign}(Y_j)c_j^\top e\ge\gamma r$, with $\gamma=1/(30\sqrt d)$. Thus at a transverse witness $(z,w,x)$ there exists an $e$ for which the mass of such paths is at least $p(x)/(12(d-1)q)$. Denote by $E_e(w)$ the set of $z$ admitting a witness with at least that mass for the specified direction, using only $j\in T$. These sets cover $E_T(w)$.

Put $h=8\eta/\gamma$ and $\delta=\eta r$. For $t\ge h$, the same paths at $(z',w+te)$, for every $|z'|\in D_a$, have absolute value at least $x+\gamma rt-4ba\ge x+\delta$. Since $x\delta\ge\eta(a-K)\to\infty$, their mass exceeds $3p(x+\delta)/q$ for large $u$ depending on $d$. It follows from \eqref{S:eq:mf-feasible} that every admissible threshold at $(z',w+te)$ is at most $x+\delta$. In particular, if both $E_e(w)$ and $E_e(w+te)$ are nonempty, then
\begin{equation}\label{S:eq:mf-shift}
 \sup_{z'\in E_e(w+te)}|z'|\le \inf_{z\in E_e(w)}|z|+2K+r\delta,\qquad t\ge h .
\end{equation}
The infimum is valid since the threshold associated to each $z$ on the right satisfies $rx\le |z|+K$.

Write $w=w_0+ye$, $w_0\perp e$. Fix $w_0$ and take $M=\sqrt{8\log u}$. Where the fiber $E_e(w_0+ye)$ is nonempty let $I_y$ be the interval between its smallest and largest absolute values, using infimum and supremum. On empty fibers assign zero length. Each fiber has Lebesgue measure at most $2|I_y|$. For any translated lattice of step $h$ in $[-M,M]$, enumerate its occupied sites by $y_1<\cdots<y_N$. By \eqref{S:eq:mf-shift} the upper end of $I_{y_{j+1}}$ is bounded by the lower end of $I_{y_j}$ plus $2K+r\delta$. Consequently
\[
 \sum_{j=1}^N |I_{y_j}|\le a+(2K+1)N,\qquad N\le1+2M/h .
\]
Integrating the lattice offset over $[0,h)$ gives $\int_{-M}^M |I_y|\,dy\le ha+(2K+1)(h+2M)$. Multiplying by the upper bound on the Gaussian density, using $|E_e|\le 2a$ on $|y|>M$, then integrating over $w_0$ yields
\[
 \mathbb E_w |E_e(w)|=O_d(ha+(1+K)(1+M)+a\Pr\{|N(0,1)|>M\})=O_d(\sqrt a+\log u).
\]
The admissibility tests can use rational intermediate thresholds together with $u$, by left continuity of the inclusive tails. For the integrations, there are finitely many candidate ranks, their conditional tails are Borel, and the projected fiber endpoints are universally measurable. In particular the bound applies to mixtures with deterministic comparison coordinates. Summing over the finitely many signed directions proves $\mathbb E_w |E_T(w)|=O_d(\sqrt a+\log u)$.

\paragraph{Completion.} Summing the three estimates over the dyadic intervals gives expected total witness length at most $(1-r)u+O_d(\sqrt u)$. Indeed their disjoint interiors give total main term $|[L,u]\cap[ru-K,u]|\le(1-r)u+K$. Their number is $O(\log u)$, the geometric sum $\sum_a\sqrt a=O(\sqrt u)$, and $(\log u)^2=O(\sqrt u)$. Substituting in \eqref{S:eq:mf-height} and \eqref{S:eq:mf-conditional}, and using $(1-r)/\sqrt{1-r^2}\le1$ and $(1-r^2)^{-1/2}\le5/4$, bounds the integrated central contribution by $u/(2\sqrt{2\pi})+O_d(\sqrt u)$ before multiplying by $q/m$. For $3/5\le r<1$ use \eqref{eq:factor-strip} and $\sqrt{(1-r)/(1+r)}\le1/2$; for $r=1$ the central integral vanishes, and for $r=0$ the entire contribution is at most $q/m$. The common radial cutoff of \eqref{eq:factor-quantiles} costs at most $q$ in FDR outside its ball, while the scalar shell between $u$ and its radius costs $O_d(1+\log u)$ per null before multiplying by $q/m$, by \eqref{eq:factor-outer}. Summing proves the upper bound. The matching asymptotic lower bound follows from the embedded one-factor constructions of \citet[Theorem 2.1 and Proposition 2.2]{Lei}.
\end{proof}

\subsection{Growing dimension and the number of loading directions}\label{S:sec:growing-one-sided}
The scalar proof of Proposition~\ref{thm:fixed-factors} gives the one-sided assertion of Corollary~\ref{cor:mf-growing-sharp} under the directional condition below.
The proof of \eqref{eq:factor-cutoff} only uses $G$-conditional independence of each null $p$-value from all other $p$-values. Thus the scalar upper bound \eqref{eq:factors-main} also allows alternative mechanisms depending arbitrarily on the common factors and on auxiliary randomness independent of the null residuals. The lower constructions cited above already belong to the fully Gaussian class in \eqref{eq:factor-d-model}.

\paragraph{A growing set of factor directions.}
The same leading upper bound holds for a sequence of models if there exists $v\ge u$ for which
\[
 \Pr\{\max_{i\in I_0:\,r_i>0}|Y_i|>v\}\le q,
 \qquad \log\{Q(u)/Q(v)\}=o(u).
\]
The maximum over an empty index set is taken to be zero. Indices with $r_i=0$ individually contribute at most $q/m$ and can be dropped from the maximum. The hazard lower bound $h(x)\ge x$ implies $u(v-u)\le\log\{Q(u)/Q(v)\}$, and hence $v-u=o(1)$. For $r_i\le1/2$ the shell integral is then $o(1)$. For $r_i>1/2$ it is $o(u)$ and the same fixed deficit in \eqref{eq:factor-strip} absorbs it. In particular the $O(q)$ remainder in \eqref{eq:factors-main} has an absolute constant eventually along such a sequence.

If the nonzero null loading vectors span at most $d=d(q)$ dimensions, a sufficient condition is $(d+1)\log u=o(u)$. Project $G$ onto their span and set $v^2=u^2+4(d+2)\log u$. The exponential moment of $\chi_d^2$ at $1/2-1/u^2$ gives, for $u>\sqrt2$,
\[
 \Pr(\chi_d>v)\le(u^2/2)^{d/2}\exp(-v^2/2+v^2/u^2)=o(Q(u)) .
\]
Here $v^2/u^2=1+o(1)$ and $Q(u)\asymp e^{-u^2/2}/u$. Mills' bounds give $\log\{Q(u)/Q(v)\}=O((d+2)\log u)=o(u)$. Alternatively, if there are $B$ distinct unoriented nonzero null loading directions, choose $Q(v)=q/(2B)$. A union bound gives the required tail probability, and $\log B=o(u)$ suffices in every ambient dimension. Neither criterion restricts the number of nonnull loading directions.

\paragraph{A criterion in terms of Gaussian regressions.}
For a fully Gaussian vector of any covariance, the directional maximum criterion also applies to normalized leave-one-out predictors: let $r_iY_i=\mathbb E[T_i-\theta_i\mid T_{-i}]$ with $r_i^2$ the variance of the predictor, omitting $r_i=0$. Gaussian regression makes the residual independent of $T_{-i}$, including for a singular joint law. If $A=\{\max_i|Y_i|\le v\}$, use
\[
 \mathbb E[\mathrm{FDP}\,\mathbf1_A]
 \le\sum_{i\in I_0}\mathbb E\left[
 \frac{\mathbf1_{\{P_i\le qR_i/m\}}}{R_i}\,\mathbf1_{\{|Y_i|\le v\}}\right],
\]
where the final indicator is taken as one when $r_i=0$. Inside each expectation $Y_i$ and $R_i$ are measurable with respect to $T_{-i}$. Thus conditioning on $T_{-i}$ gives the same scalar bound as above. This argument requires no mutual independence of the different regression residuals. A given model may admit more economical auxiliary factor directions than the leave-one-out predictors, so either representation can be used.

\subsection{The growing-dimensional two-sided bound}\label{S:app:growing-two-sided}
\begin{proof}[Proof of the two-sided assertion of Corollary~\ref{cor:mf-growing-sharp}]
Write $K=\sqrt{8\log u}$. In the proof of Theorem~\ref{thm:multifactor-sharp} use the lower cutoff $L=K+d^{1/4}\sqrt u$. On each dyadic shell $D_a=[a,\min(2a,u)]$, $a=2^jL<u$, take $\eta_a=d^{-1/4}a^{-1/2}$ and $b=\eta_a r/a$. Under the weaker assumption $d=o(u^2)$, $L=o(u)$ and
\[
 \min_a\eta_a(a-K)\ \gtrsim\ (u^2/d)^{1/8}\longrightarrow\infty .
\]
The scalar bounds carry absolute constants. The discarded central interval costs $O(L)$. The additive error in \eqref{S:eq:mf-axial} on this shell is $O(K+1/\eta_a)=O(K+d^{1/4}\sqrt a)$, and the small-total-loading part still costs $O(K+1)$ since $4ba\le 2r$ for large $u$. We estimate the strong transverse contribution without a sum over $d$ coordinate directions.

For $d\ge2$ take $e$ uniformly distributed on the unit sphere of $v_i^\perp$, independently, and put $k=d-1$. For any fixed unit vector $v_0$ there, the coordinate $U=v_0^\top e$ satisfies $\E U^2=1/k$, $\E U^4=3/[k(k+2)]$. Paley--Zygmund applied to $U^2$ and symmetry give $\Pr\{U\ge(2k)^{-1/2}\}\ge1/24$, also for $k=1$. Every path in class $T$ has $\|c_j\|\ge r/30$ for sufficiently large $u$. At a fixed transverse witness $(z,w,x)$ in this shell, write $W\ge p(x)/(6q)$ for the mass of its exceeding $T$ paths. Let $W_e$ be the submass for which $\operatorname{sign}(Y_j)c_j^\top e\ge r/(50\sqrt d)$. Then $0\le W_e\le W$ and $\E_e W_e\ge W/24$. Consequently $\Pr_e\{W_e\ge W/48\}\ge1/48$. Let $E_e(w)$ be the set of $z$ admitting a transverse witness for which $W_e\ge p(x)/(288q)$. This gives the pointwise inequality
\[
 \mathbf1_{\{z\in E_T(w)\}}\le48\E_e\mathbf1_{\{z\in E_e(w)\}}.
\]
The same inequality holds when several ranks witness a given point, by choosing the first qualifying rank in the finite list.

For fixed $e$, set $h=400\eta_a\sqrt d$ and $\delta=\eta_a r$. On any favored path, a shift $w\mapsto w+te$ with $t\ge h$ increases its signed value by at least $8\eta_a r$, while arbitrary $|z'|\in D_a$ in place of $z$ costs at most $4ba$. At the new point those paths exceed $x+\delta$ and carry mass at least $p(x)/(288q)>3p(x+\delta)/q$, since $x\delta\ge\eta_a(a-K)\to\infty$ uniformly in the shell. Admissibility therefore forces every future witness in this shell to have magnitude $|z'|\le |z|+C_0$, $C_0=2K+r\delta$. On each line $w=w_0+se$, $w_0\perp e$, define $H(s)$ to be the supremum of $|z'|$ over witnesses in this shell at coordinates $s'\ge s$, taking $a$ for an empty set. This function is nonincreasing and takes values in $[a,2a]$. Each nonempty oriented fiber has all its magnitudes between $H(s+h)-C_0$ and $H(s)$. Thus
\[
 |E_e(w_0+se)|\le2[H(s)-H(s+h)+C_0],\qquad
 \int\phi(s)[H(s)-H(s+h)]\,ds\le h\|\phi\|_\infty a.
\]
The last bound follows by integrating the variation measure of $H$ over intervals of length $h$. Conditional on $w_0$, the line coordinate is standard normal. Integrate first over the line, then over $w_0,e$, to obtain $\E_w|E_T(w)|=O(ha+K)=O(d^{1/4}\sqrt a+K)$ with absolute constants. For $d=1$ this class is empty. The disjoint shell interiors count their axial main lengths only once, $\sum_a\sqrt a=O(\sqrt u)$, and $K\log u=O(\sqrt u)$. The central contribution is therefore at most $u/(2\sqrt{2\pi})+O(d^{1/4}\sqrt u)$ per null before its $q/m$ multiplier.

Finally take $v^2=u^2+4(d+2)\log u$. The exponential moment of $\chi_d^2$ at $1/2-u^{-2}$ gives
\[
 \Pr\{\|G\|>v\}\le(u^2/2)^{d/2}\exp(-v^2/2+v^2/u^2)=o(q),
\]
since $v^2/u^2=1+o(1)$. Mills' bounds imply $\log\{\bar\Phi(u)/\bar\Phi(v)\}=O((d+2)\log u)$. By \eqref{eq:factor-outer} this bounds the additional shell integral up to an absolute factor. We have shown
\[
 \FDR\le\pi_0 q\frac{u}{2\sqrt{2\pi}}+O\!\left(q\{d^{1/4}\sqrt u+(d+2)\log u\}\right).
\]
The remainder is $o(qu)$ under the stated assumption. Embedding the one-factor lower construction of \citet[Theorem~2.1 and Proposition~2.2]{Lei} gives the matching lower bound.
\end{proof}

The central estimate in this proof only requires $d=o(u^2)$. Thus the same sharp upper coefficient applies along any sequence of fully Gaussian models of factor dimension $d=o(u^2)$ having a deterministic $v\ge u$ such that $\Pr\{\max_{i\in I_0:r_i>0}|v_i^\top G|>v\}\le q$ and $\log\{\bar\Phi(u)/\bar\Phi(v)\}=o(u)$. Use this simultaneous directional cutoff in place of the full radius in the preceding proof. For example, the conditions hold when the null loading vectors span at most $k$ dimensions with $(k+2)\log u=o(u)$, by the same radial calculation restricted to that span, or when the number $N$ of distinct unoriented nonzero null loading directions obeys $\log(N+1)=o(u)$, by taking $\bar\Phi(v)=q/(2N)$ when $N\ge1$. With $N=0$, all nulls have zero loadings and the FDR is at most $\pi_0q$.

\subsection{An improved remainder for aligned loadings}
\begin{corollary}\label{S:cor:sharp-cone-remainder}
Fix $d$ and $c>0$. Suppose that for every true null with $0<r_i=\|b_i\|<3/5$ and each $j\ne i$, $|b_i^\top b_j|\ge c r_i\|b_j\|$. Then the two-sided model satisfies
\[
 \FDR\le\pi_0q\,\frac{u}{2\sqrt{2\pi}}+O_{d,c}\!\big(q(1+\log u)^{3/2}\big),\qquad u=\bar\Phi^{-1}(q/2).
\]
The condition holds for all one-factor models with $c=1$. Consequently, as $q\downarrow0$ and $D=\log(1/q)$,
\[
 B_{1,2}(q)=\frac{q\sqrt D}{2\sqrt\pi}+O\!\big(q[1+\log D]^{3/2}\big).
\]
\end{corollary}

\begin{proof}
In the dyadic proof take $J=\lceil C_0\log u\rceil$, with $C_0$ large enough that $u e^{-c_0J}\le1$ for the exponent $c_0$ in the conditional count bound. Keep $K=\sqrt{8\log u}$ and put $L=2K$. The estimate \eqref{S:eq:mf-height} on $|z|\ge L$ now has error $O(u^{-36})$ in place of $Ce^{-u}$. Indeed, the Mills-ratio identity in \citet[(207)--(210)]{Lei} bounds the term with the matching sign by the same leading factor as \eqref{S:eq:mf-height}. For the term with the opposite sign its prefactor is bounded and $|z|+rx\ge 2L-K=3K$, so its exponential factor is at most $e^{-9K^2/2}=u^{-36}$. The central region $|z|<L$ costs $O(K)$ by the unrestricted density bound for $r<3/5$.

On each dyadic interval $D_a$ take $\eta=1/4$ and $b=r/(4a)$. For large enough $u$ depending on $c$, $a\ge L\ge 5/c$. Thus $|a_j|<b$ implies $\|b_j\|\le |a_j|/c<r/(4ac)\le r/20$, and $T$ is empty. The axial exclusion width is $2K+4+40a/(a-K)=O(K+1)$ because $a\ge 2K$, and the small-group displacement is $4ab=r$. The dyadic small-total-loading argument therefore gives the same $O(K+1)$ length in each interval. The total error in the witness length is $O((K+1)\log u)$; the axial main lengths again sum to at most $(1-r)u+K$. Ranks up to $J$ cost $O(\log u)$ and the common radial shell costs $O_d(1+\log u)$, proving the assertion. Lei's Theorem~2.1 and Proposition~2.2 give the matching lower bound with $O(q)$ error in the one-factor class. Finally $u=\sqrt{2D}+O((1+\log D)/\sqrt D)$.
\end{proof}

\subsection{Dependent Gaussian residuals}\label{S:sec:dependent-residuals}

In \eqref{eq:factor-d-model} allow the residual vector to have an arbitrary Gaussian correlation matrix, still independent of $G$. We say the conditional law is weakly PRDS on the true nulls if, for almost every $g$, for every true null $i$ and coordinatewise increasing Borel set $D$, the function
\[
 t\longmapsto \Prob(P\in D\mid P_i\le t,G=g)
\]
is nondecreasing wherever its conditioning event has positive conditional probability. Only the sets $D=\{R\le j\}$ will be needed.

\begin{theorem}[One-sided tests with dependent residuals]\label{S:thm:dependent-residuals}
Fix $d\ge1$. If the right-sided $p$-values are conditionally weakly PRDS on the true nulls, the upper bound \eqref{eq:factors-main} holds uniformly. This condition holds whenever every true null has nonnegative residual covariance with every other coordinate. The supremum over either of these two model classes is $q\sqrt{\log(1/q)}/\sqrt\pi+O_d(q)$. Under conditional weak PRDS the same upper bound \eqref{eq:factors-main} holds for two-sided tests.
\end{theorem}

\begin{theorem}[Sharp two-sided constant on sparse residual graphs]\label{S:thm:residual-blocks}
Fix $d\ge1$ and use two-sided $p$-values in the factor model with dependent Gaussian residuals. Form the residual covariance graph with an edge for each nonzero off-diagonal residual covariance. Suppose its vertices can be properly colored with $B\ge1$ colors, each true null has at most $B-1$ neighbors, and the conditional $p$-values given $G$ are weakly PRDS on the true nulls. The graph conditions hold, in particular, if the maximum degree is at most $B-1$, or if the residuals split into independent blocks of size at most $B$. For $u=\bar\Phi^{-1}(q/2)$ and all sufficiently small $q$ depending only on $d$,
\begin{equation}\label{S:eq:two-block-bound}
 \FDR\le q+\pi_0q\left\{\frac{u}{2\sqrt{2\pi}}
       +C_d\big[\sqrt u+\log(B+1)\big]\right\}.
\end{equation}
If $B=B(q)\ge1$ and $\log(B+1)=o(\sqrt{\log(1/q)})$, the supremum over this class is $(1/(2\sqrt\pi)+o(1))q\sqrt{\log(1/q)}$.
\end{theorem}

\paragraph{A positive definite example.}
This class includes positive definite covariance matrices of unbounded independent-residual factor dimension. Take $d=1$ and write $Z=G$. Within each of $K$ independent residual pairs take correlation $0<\eta<1$, assign both coordinates the same loading $0<a<1$, and give them a common mean $\theta_h$; any subset of the pairs may be true nulls. All coordinates then share the common factor. Every covariance matrix in this construction is positive definite. The cross-covariance of the $K$ first coordinates with the $K$ second coordinates is
\[
 a^2\mathbf1\mathbf1^\top+(1-a^2)\eta I_K.
\]
The two index sets are disjoint, so in an independent-coordinate-residual representation this whole matrix is a product of two $K\times d$ loading matrices, one transposed. Its rank $K$ forces $d\ge K$.

For completeness, conditional weak PRDS holds in each folded pair even with the nonzero shift induced by $Z$. Write $X,Y$ for its two conditional statistics, with common mean $\mu$, marginal variance $\sigma^2=1-a^2$, and correlation $\eta$. A joint sign change permits $\mu\ge0$. Given $|X|=x$, the probability of $X=x$ is $w(x)=(1+e^{-2\mu x/\sigma^2})^{-1}$. For $v\ge0$, set $g(m)=\Prob(|N(m,s^2)|\ge v)$, where $s^2=\sigma^2(1-\eta^2)$, and write $c=(1-\eta)\mu$, $t=\eta x$. The conditional survival probability of $|Y|$ is $w(x)g(c+t)+(1-w(x))g(c-t)$. Its derivative is
\[
 w'(x)[g(c+t)-g(c-t)]
   +\eta\{w(x)g'(c+t)-(1-w(x))g'(c-t)\}.
\]
The first term is nonnegative since $g$ is even and nondecreasing on $[0,\infty)$. The second is nonnegative when $c\le t$. For $c>t$ and $v>0$, use $g'(m)\propto e^{-m^2/(2s^2)}\sinh(vm/s^2)$ to get
\[
 \frac{w(x)g'(c+t)}{(1-w(x))g'(c-t)}
 \ge \exp\{2\mu x/\sigma^2-2ct/s^2\}\ge1,
\]
because $ct/s^2=\eta\mu x/[(1+\eta)\sigma^2]$. For $v=0$ the survival function is constant. Therefore $|Y|$ given $|X|=x$ is stochastically nondecreasing in $x$. By the coordinatewise monotonicity of $2\bar\Phi(|X|)$ and by independence from other pairs given $Z$, the conditional two-sided $p$-values are weakly PRDS on every true null.

\subsection{Proof of Theorem~\ref{S:thm:dependent-residuals}}
\begin{proof}
Write $t_k=qk/m$, $F_{i,g}(t)=\Prob(P_i\le t\mid G=g)$, and, whenever $F_{i,g}(t_k)>0$, put
\[
 w_{ik}(g)=\Prob(R=k\mid P_i\le t_k,G=g).
\]
Put $w_{ik}(g)=0$ otherwise. The step-up definition gives the identity
\begin{equation}\label{S:eq:conditional-rank-mass}
 \E[\FDP\mid G=g]=\frac qm\sum_{i\in I_0}\sum_{k=1}^m
 w_{ik}(g)\frac{F_{i,g}(t_k)}{t_k}.
\end{equation}
For each $i$ let $A_k=\Prob(R\le k\mid P_i\le t_k,G=g)$ when defined and zero otherwise, and put $A_0=0$. Since $\{R\le k-1\}$ is increasing in the $p$-values, weak PRDS gives $w_{ik}(g)\le A_k-A_{k-1}$, including the first index with $F_{i,g}(t_k)>0$. Thus $\sum_k w_{ik}(g)\le1$. Replacing the CDF ratio in \eqref{S:eq:conditional-rank-mass} by the boundary-null envelope $D_{r_i}(Y_i)$ gives \eqref{eq:factor-cutoff} for every factor-measurable event. The scalar estimates and radius truncation in Section~\ref{app:factors} now apply verbatim.

For right-sided tests with the stated covariance sign condition, conditional on $G$ and on a nondegenerate $T_i=x$, the other statistics have an $x$-independent Gaussian noise term and a mean nondecreasing coordinatewise in $x$, by Gaussian regression. The right-sided $p$-values decrease in their statistics, so the conditional probability of an increasing $p$-value event is nonincreasing in $x$. Truncating $T_i$ at its upper-tail cutoff yields weak PRDS. A coordinate degenerate given $G$ makes the truncation event either null or certain. The one-sided lower construction of \citet[Theorem~3.1 and Proposition~3.2]{Lei} has independent residuals and belongs to both classes.
\end{proof}

\subsection{Proof of Theorem~\ref{S:thm:residual-blocks}}
\begin{proof}
Let $u_k=\bar\Phi^{-1}(t_k/2)$ and write $g$ for a realization of $G$. Use the weights from \eqref{S:eq:conditional-rank-mass}. Define
\[
 M_g(v)=\frac1m\sum_{\ell=1}^m
 \Prob\bigl(\lvert\theta_\ell+b_\ell^\top g+\sqrt{1-\|b_\ell\|^2}\,\varepsilon\rvert\ge v\bigr),
 \qquad \varepsilon\sim N(0,1).
\]
Call rank $k$ admissible at $g$ if $mM_g(u_k)\ge k/2$ and $mM_g(u_j)\le2j$ for every $j>k$. These discrete conditions imply \eqref{S:eq:mf-feasible} with comparison mixture $M_g$ by the same interpolation as in Theorem~\ref{thm:multifactor-sharp}.

For $k\ge4B$ that is not admissible we claim, uniformly in $i,g$,
\begin{equation}\label{S:eq:residual-block-deviation}
 w_{ik}(g)\le e^{-ck/B}
\end{equation}
for a universal $c>0$. It suffices to consider conditioning events of positive probability. Remove $i$ and its neighbors, a set of size at most $B$. At threshold $u_j$ let $S_j$ count the remaining exceedances and $\mu_j=\E(S_j\mid G=g)$. Conditional on $G=g$, the whole remaining residual vector is independent of $P_i$ by Gaussianity and zero cross-covariance, including when the covariance is singular. Its joint law therefore survives further conditioning on $P_i\le t_k$. Write $S_j=\sum_{h=1}^B S_{jh}$ for the color sums. Within each color, centered Gaussian residuals are mutually independent. H\"older's inequality across colors, followed by the Bernoulli moment bound within colors, gives for any real $\lambda$
\[
 \E(e^{\lambda S_j/B}\mid G=g,P_i\le t_k)
 \le \prod_{h=1}^B\E(e^{\lambda S_{jh}}\mid G=g)^{1/B}
 \le \exp\{(e^\lambda-1)\mu_j/B\}.
\]
If $mM_g(u_k)<k/2$, then $\mu_k<k/2$ and $R=k$ forces $S_k\ge k-B\ge3k/4$. Taking $\lambda=\log(3/2)$ yields an exponent at most $-[\tfrac34\log(3/2)-\tfrac14]k/B$. Otherwise, fix any $j>k$ with $mM_g(u_j)>2j$. Then $\mu_j>2j-B\ge7j/4$ and $R=k$ forces $S_j<j$. Taking $\lambda=-\log(3/2)$ gives an exponent at most $-[\tfrac7{12}-\log(3/2)]j/B$. Both bracketed constants are positive. This proves \eqref{S:eq:residual-block-deviation} uniformly in the probability of the conditioning event.

Take $J=\max\{\lceil\sqrt u\rceil,\lceil C_0 B[1+\log(B+1)+\log(u+1)]\rceil\}$ with $C_0$ a sufficiently large universal constant. Summation of \eqref{S:eq:residual-block-deviation} gives
\[
 \sum_{k>J:\,k\ {\rm not\ admissible}}w_{ik}(g)
 \le C_1 B e^{-cJ/B}\le (u+1)^{-1}.
\]
Fix a true null $i$ and decompose $G=zv_i+w$ along its loading as in Theorem~\ref{thm:multifactor-sharp}. When $r_i=0$ choose any unit $v_i$. For $|z|\le u$, integrating the contribution of inadmissible ranks above $J$ against the joint Gaussian density costs $O(1)$ before multiplication by $q/m$: their weights sum to at most $(u+1)^{-1}$ and \eqref{eq:factor-strip} bounds the integral of the scalar envelope $D_{r_i}(z)$ by $1/2+u/\sqrt{2\pi}$. On admissible ranks the weights sum to at most one, so their weighted sum of ratios is bounded by the supremum of $\rho_{r_i}(u_k,z)$ over thresholds satisfying \eqref{S:eq:mf-feasible}. The deterministic central estimate in the proof of Theorem~\ref{thm:multifactor-sharp} bounds this integrated supremum by $u/(2\sqrt{2\pi})+O_d(\sqrt u)$. The proof uses only the univariate Gaussian marginal on each affine path and total comparison mass at most one. Both properties hold for the pooled measure used here, including its atom for $i$. For $r_i=0$ use the envelope one, and for $r_i\ge3/5$ apply \eqref{eq:factor-strip} directly.

For $d=1$ the same scalar bound follows directly as $u/(2\sqrt{2\pi})+O(\sqrt u)$ by \citet[Lemma~4.7, equations~(211) and (261)--(264)]{Lei}: use the probability measure $m^{-1}\sum_\ell\delta_{(\theta_\ell,b_\ell)}$, the admissible weights, and $r=|b_i|$. The requisite support and mass conditions hold by construction.

Ranks $1\le R\le J$ contribute at most $\pi_0q H_{J\wedge m}$ to the unconditional FDR: when $t_{j-1}<P_i\le t_j$ and $i$ is rejected, $R\ge j$, and each null interval has probability $q/m$. Here $t_0=0$. The harmonic sum is at most $C[1+\log(B+1)+\sqrt u]$. For the remaining ranks discard the event $\|G\|>v$, with $v$ as in \eqref{eq:factor-quantiles}, only once for the whole FDP at cost $q$. For each $i$, its part in the ball with $u<|z|\le v$ costs $O_d(1+\log u)$ before multiplication by $q/m$, using the unit mass bound and \eqref{eq:factor-outer}. These terms give \eqref{S:eq:two-block-bound}. The independent-residual lower construction of \citet[Theorem~2.1 and Proposition~2.2]{Lei} embeds in every $d\ge1$.
\end{proof}

\section{The number of alternatives and the dimension}\label{S:app:finite}
\subsection{The number of alternatives}\label{S:sec:finite}

Theorem~\ref{thm:main} allows the nonnull $p$-values to be arbitrary functions of the true-null statistics. Within this larger model one can quantify the effect of the number of signals. Let $n\ge1$ be the number of true nulls, let $s\ge0$ be the number of false nulls, and set $m=n+s$. For a fixed testing convention, write $W_{n,s}(q)$ for the supremum of the ordinary BH FDR over the allowed Gaussian true-null laws and all joint completions by $s$ nonnull $p$-values. Put
\[
 a=\frac{nq}{m},\qquad L=\min\{s+1,1/a\},\qquad K=\lfloor L\rfloor,\qquad
 H_K=\sum_{r=1}^K\frac1r .
\]
Here $s$ counts the available nonnull hypotheses, whether or not BH rejects them.

\begin{theorem}[Finite-sample order with arbitrary alternatives]\label{S:thm:finite}
For two-sided Gaussian tests take $C=e$. For one-sided tests take $C=2e$. For every $q\in(0,1)$, $n\ge1$ and $s\ge0$,
\begin{equation}\label{S:eq:finite}
 \frac{a}{e}H_K
 \ \le\ W_{n,s}(q)\
 \le\ \min\{B(Ca),\;Ca[1+\log(s+1)]\}.
\end{equation}
In particular, with absolute constants independent of all three parameters,
\[
 W_{n,s}(q)\asymp a\big[1+\log\min\{s+1,1/a\}\big].
\]
\end{theorem}

The upper bound uses a counting constraint. If at least one true null is rejected and at most $s$ signals can be rejected, a positive FDP must be at least $1/(s+1)$. This truncates the integral of the tail bound in \eqref{eq:fdx}. For the lower bound, choose independent Gaussian nulls and use the smallest null $p$-value to determine how many alternatives to reject. This construction specializes the adversary of \citet[Theorem~4]{Su}; a finite-count formulation is given by \citet[Lemma~5.1]{LyuWang}. Section~\ref{S:app:finite} gives the bounds and their refinements for independent nulls.

With $s$ fixed the bound is linear in $q$, uniformly in $n$. If $a\to0$ and $s+1\ge1/a$, the worst-case order is $a\log(1/a)$. The constructions that attain this larger rate let the alternatives depend on the observed nulls. Section~\ref{S:sec:one-null} gives a case in which restricting the signals to Gaussian mean shifts changes the leading coefficient.

\subsection{The leading coefficient as the dimension grows}\label{S:sec:window}

We retain the model of Theorem~\ref{S:thm:finite}, with $W_{n,s}$ taken for any fixed assignment of one- and two-sided null tests. Composite one-sided true nulls can be included. The bound on dimension below links growth in the number of tests to the scale of a Gaussian extreme.
\begin{theorem}[Sharp coefficient for growing dimension]\label{S:thm:window}
For sequences $n=n(q)\ge1$ and $s=s(q)\ge0$ satisfying
\[
 \log^2(m+1)=o\bigl(\log(1/q)\bigr),\qquad q\downarrow0,
\]
we have
\[
 \frac{W_{n,s}(q)}{(nq/m)H_{s+1}}\longrightarrow1 .
\]
\end{theorem}

In particular, for every fixed $n,s$, the supremum divided by $nq/m$
tends to $H_{s+1}$. When $s=0$ the result gives first-order validity
uniformly over correlation matrices throughout the stated range of
growing dimensions. For fixed dimension with all tests two-sided, \citet[Theorem~2.1]{Zhang} proves the stronger exact inequality at sufficiently
small levels. For a family of one-sided tests, \citet[Proposition~2.2 and its quantitative estimate]{LyuWang} prove first-order validity at fixed dimension and give a
quantitative bound implying it when $n^2\log(n+1)=o(\log(1/q))$.

For instance, with two unrestricted alternatives the worst FDR is asymptotic to $(11/6)nq/m$ throughout the stated regime. To obtain the upper bound, we condition on one small null $p$-value and distribute the FDP among null indices in proportion to $1/P_i$. Gaussian regression bounds how quickly the remaining ratios $P_i/P_j$ can change within the relevant tail range. The resulting extremum over an ordered simplex occurs when a single coordinate carries the mass. Independent nulls attain the bound to first order. Section~\ref{S:app:window} gives the proof.

\subsection{Finite-sample bounds and independent nulls}\label{S:app:finite-bounds}

\begin{proof}[Proof of Theorem~\ref{S:thm:finite}]
On $\{V>0\}$, BH rejects at most $s$ false nulls. Thus $\FDP\ge1/(s+1)$ whenever it is positive. For $0<x\le1/(s+1)$ the event $\{\FDP\ge x\}$ coincides with $\{\FDP\ge1/(s+1)\}$. Its probability is at most $Ca(s+1)$ by \eqref{eq:fdx}. On $1/(s+1)<x\le1$, integrate the bound $Ca/x$. This gives the second upper bound in \eqref{S:eq:finite}; the first is Theorem~\ref{thm:main}.

For the lower bound, take $n$ independent Gaussian nulls, so their $p$-values are independent uniforms for either testing convention. Let $U$ be their minimum. For $r=1,\ldots,K$ let $I_r=((r-1)q/m,rq/m]$. When $U\in I_r$, assign $r-1$ nonnull $p$-values the value zero and the remainder one. Outside these intervals assign all nonnulls one. This uses only $s$ alternatives. On $I_r$, rank $r$ is feasible, all $r-1$ zeros are rejected, and no $p$-value equal to one is rejected since $q<1$. Any further rejections are true nulls. Hence $\FDP\ge 1/r$. Integrating the density of the minimum gives the explicit lower bound
\begin{equation}\label{S:eq:lower-sum}
 W_{n,s}(q)\ \ge\
 \sum_{r=1}^K\frac{(1-(r-1)q/m)^n-(1-rq/m)^n}{r}.
\end{equation}
For $n\ge2$, every interval has upper endpoint at most $Kq/m\le1/n$. On these intervals the density $n(1-x)^{n-1}$ is at least $n(1-1/n)^{n-1}\ge n/e$. For $n=1$ the density equals one. Each interval therefore has probability at least $a/e$, which gives the lower bound in \eqref{S:eq:finite}.

Finally, $H_{\lfloor L\rfloor}\ge\max(1,\log L)\ge(1+\log L)/2$. Also $B(Ca)\le Ca[1+\log(1/a)]$: for $Ca<1$ this follows from $C\ge1$, and for $Ca\ge1$ it follows from $B(Ca)=1$. Taking the smaller upper bound proves the uniform order assertion.
\end{proof}

The lower construction uses the smallest null $p$-value to choose the number of rejections supplied by the alternatives. It specializes the adversarial construction in \citet[Theorem~4 and Section~2.2]{Su}. The finite-count oracle of \citet[Lemma~5.1 and the finite-count integral bound]{LyuWang} also accounts for how many alternatives are available. When $n=1$, one may use all $m$ intervals in \eqref{S:eq:lower-sum}; the exact worst-case value is $(q/m)H_m$, as obtained in their Theorem~5.2 for one-sided Gaussian dependence. In fact the same identity for our enlarged model follows because on $I_r$ any rejection of the unique null requires at least $r$ total rejections, and the construction attains exactly $1/r$.

Independent true nulls give a sharper comparison. If $W_{n,s}^{\mathrm{ind}}(q)$ denotes the same supremum with independent uniform true-null $p$-values, then
\begin{equation}\label{S:eq:iid}
 \big|W_{n,s}^{\mathrm{ind}}(q)-a\log L\big|\le a .
\end{equation}
Indeed, the exact independent Simes identity \citep{Simes} gives the upper bound with $C=1$. For the lower bound, Bernoulli's inequality on $I_r$ gives $(1-x)^{n-1}\ge1-(n-1)x\ge1-ra$. Hence the sum in \eqref{S:eq:lower-sum} is at least $aH_K-a^2K\ge a(\log L-1)$. In particular the leading coefficient of $a\log L$ equals one whenever $L\to\infty$ in this independent-null model.

A second regime has an explicit leading constant. For every fixed $s\ge0$,
\begin{equation}\label{S:eq:iidfixed}
 W_{n,s}^{\mathrm{ind}}(q)=aH_{s+1}+O_s(a^2)
 \quad\text{as } a\downarrow0,\quad\text{uniformly over }n .
\end{equation}
To see this, write $U_{(k)}$ for the $k$th null order statistic and consider the event that $U_{(k)}\le(k+s)q/m$ for some $k\ge2$. A binomial union bound gives this event probability at most
\[
 \sum_{k=2}^{n}\binom nk\big((k+s)q/m\big)^k
 \le \sum_{k=2}^{n}\big[e a(1+s/2)\big]^k=O_s(a^2).
\]
Outside this event at most one true null can be rejected, whatever values the alternatives take. On $I_r$, the FDP in that case is at most $1/r$, since rejecting the smallest null requires at least $r$ rejections. Above $I_{s+1}$ the sole null rejection would be impossible. Hence the sum in \eqref{S:eq:lower-sum}, with $K=s+1$, bounds $W^{\mathrm{ind}}_{n,s}(q)$ to within $O_s(a^2)$. Each of its $s+1$ interval probabilities lies between $a(1-ra)$ and $a$. This proves \eqref{S:eq:iidfixed}. For example, with one unrestricted alternative the first-order worst FDR is $3a/2$ under independent true nulls.

\subsection{Proof of Theorem~\ref{S:thm:window}}\label{S:app:window}

\paragraph{Reduction to the true-null values.}
Set $\tau=q/m$, write $P_1,\ldots,P_n$ for the true-null values and
$P_{(1)}\le\cdots\le P_{(n)}$ for their order statistics. The largest
attainable FDP given these values is
\[
 G(P)=\max_{\substack{1\le k\le n\\0\le a\le s}}
       \frac{k}{k+a}\,{\bf1}\{P_{(k)}\le(k+a)\tau\}.
\]
Indeed, any nonempty set of rejections with $k$ nulls and $a$
alternatives satisfies the displayed threshold. Conversely, for a pair
attaining a positive maximum, set $a$ alternative values to zero and
the others to one. BH rejects all these zeros and at least $k$ true
nulls. Any further true-null rejections increase or preserve the FDP,
and the ones are never rejected because $q<1$. A fixed tie rule makes
this choice measurable. This is the finite-alternative oracle of \citet[Lemma~5.1]{LyuWang}, here attained because the alternatives
are unrestricted. For $0<t\le1$, define
\[
 d_k(t)=\min\{k+s,\lfloor k/t\rfloor\}.
\]
Then $G(P)\ge t$ if and only if $P_{(k)}\le\tau d_k(t)$ for some $k$.
These integer thresholds also occur in their finite-count linking bound.

\paragraph{An inequality for a fixed vector.}
For $x_1\ge\cdots\ge x_n\ge0$ with $\sum x_j=1$, put
\[
 I(x)=\int_0^1\max_{1\le k\le n}d_k(t)x_k\,dt .
\]
The function $I$ is convex on the ordered simplex. Its vertices $v_l$
have $l$ initial entries equal to $1/l$ and the rest zero. Indeed, with
$x_{n+1}=0$ the identity $x=\sum_{l=1}^n l(x_l-x_{l+1})v_l$ is a
convex combination. Since $d_k(t)$ is nondecreasing in $k$,
\[
 I(v_l)=\frac1l\int_0^1 d_l(t)\,dt
       =1+\sum_{j=1}^s\frac1{l+j}\le H_{s+1}.
\]
The middle identity follows by writing
$d_l(t)=l+\sum_{j=1}^s{\bf1}\{t\le l/(l+j)\}$. Hence
$I(x)\le H_{s+1}$ throughout the simplex.

\paragraph{Conditioning on one coordinate.}
Since $G(P)$ is nonincreasing in each null $p$-value, coupling any interior one-sided mean to its boundary value gives a pointwise upper bound. We hence take each null on its boundary, so $P_i$ is uniform on $(0,1)$. Allocate the weight
\[
 A_i=\frac{1/P_i}{\sum_{j=1}^n1/P_j}
\]
to coordinate $i$. The weights sum to one almost surely. Conditional
on $P_i=\tau z$, define $r_j(z)=P_i/P_j$ and order the ratios
decreasingly as $r_{[1]}(z)\ge\cdots\ge r_{[n]}(z)$. Then
\[
 A_i=\frac1{\sum_j r_j(z)},\qquad
 \{G\ge t\}
 =\left\{z\le\max_k d_k(t)r_{[k]}(z)\right\}.
\]
The change of variable $P_i=\tau z$ gives
$\mathbb E G/(n\tau)=n^{-1}\sum_i\int_0^{1/\tau}
\mathbb E[A_iG\mid P_i=\tau z]\,dz$.

Fix $0<\delta<1/2$ and put
\[
 L=nm^2/\delta,\qquad T=\log(n/\delta),\qquad
 M=\sqrt{2\log(2nL/\delta)}.
\]
In the last average the range $0<z<\delta$ contributes at most
$\delta$. On $\{G>0\}$ one has $\min_jP_j\le q=m\tau$. Thus
$\sum_{i:P_i>L\tau}A_i\le nm/L$ there, and
$\Pr(G>0)\le nq$. The total contribution from $z>L$, after division
by $n\tau$, is at most $nm^2/L=\delta$.

To couple all the remaining conditional laws, orient each one-sided
statistic so that its $p$-value uses the upper normal tail $\bar\Phi$.
Let $c_j=1$ for a one-sided null and $c_j=2$ for a two-sided null.
Under the conditioning above, write $y=\bar\Phi^{-1}(\tau z/c_i)$.
The anchor statistic is $y$ for a one-sided null and is $\pm y$
with equal probabilities for a two-sided null. Condition also on this
sign. Gaussian regression represents each coordinate as
\[
 U_j(z)=b_j y+\sqrt{1-b_j^2}\,E_j .
\]
Here $b_j$ is the correlation with the anchor multiplied by the
conditioned sign. The residual vector $E$ can be held fixed as $z$
varies. Each nonzero residual coefficient multiplies a standard normal
variable, with arbitrary dependencies between residuals. A zero
coefficient can be assigned an unused standard normal variable.
This representation covers singular covariance matrices. With
$v_j(z)=U_j(z)$ for one-sided nulls and $v_j(z)=|U_j(z)|$ for
two-sided nulls, it gives
\[
 r_j(z)=\frac{c_i\bar\Phi(y)}{c_j\bar\Phi(v_j(z))}.
\]
In particular $r_i(z)=1$. The event $\max_j|E_j|>M$ has probability
at most $2n\exp(-M^2/2)$. Its contribution from $\delta\le z\le L$
is at most $L\,2n\exp(-M^2/2)=\delta$ after normalization, since
$A_iG\le1$.

\paragraph{Control of the conditional ratios over the entire interval.}
It suffices to treat the range where $\tau L<1/2$ and
$y_*=\bar\Phi^{-1}(\tau L)$ satisfies
$y_*^2\ge C_0(1+T+M^2)$ for a sufficiently large absolute constant
$C_0$. Then $y\ge y_*>0$ for all anchors and all $\delta\le z\le L$.
On $\{\max_j|E_j|\le M\}$ define
$\widetilde r_j(z)=\max\{r_j(z),e^{-T}\}$. There is an absolute
constant $C$ such that throughout this interval
\begin{equation}\label{S:eq:window-clip}
 e^{-\eta}\widetilde r_j(1)\le\widetilde r_j(z)
 \le e^\eta\widetilde r_j(1),\qquad
 \eta=\frac{C(1+T+M^2)\log(L/\delta)}{y_*^2}.
\end{equation}
Here are the details. Write $h(x)=\phi(x)/\bar\Phi(x)$. For $x>0$
the Mills inequalities give $x\le h(x)\le x+1/x$. Whenever
$r_j(z)\ge e^{-T}$ we have
$\log\{\bar\Phi(v_j)/\bar\Phi(y)\}\le T+\log2$. If $v_j\le y/2$,
the same logarithm is at least
$\int_{y/2}^y h(x)\,dx\ge3y^2/8$, which contradicts the choice of
$C_0$. If $y/2<v_j<y$, integration of $h(x)\ge y/2$ gives
$v_j\ge y-2(T+\log2)/y$. The same lower bound holds for $v_j\ge y$.
In either case Cauchy--Schwarz gives
$v_j\le |b_j|y+\sqrt{1-b_j^2}\,M\le\sqrt{y^2+M^2}$.
Let $c=b_j$ for a one-sided coordinate and
$c=b_j\operatorname{sign}(U_j)$ for a two-sided coordinate.
Because $v_j\le cy+\sqrt{2(1-c)}M$, the lower bound just proved gives
\[
 (1-c)y\le 2(T+\log2)/y+\sqrt{2(1-c)}M,
 \qquad 1-c\le C(1+T+M^2)/y^2 .
\]
The second inequality follows by absorbing half of the quadratic term
in $y\sqrt{1-c}$. The Mills inequalities and the bounds on $v_j$
give $|h(v_j)/h(y)-1|\le C(1+T+M^2)/y^2$. Differentiating with
respect to $\log z$ now gives
\[
 \left|\frac{d}{d\log z}\log r_j(z)\right|
  =\left|1-c\frac{h(v_j(z))}{h(y)}\right|
  \le \frac{C(1+T+M^2)}{y_*^2}
\]
where $r_j(z)\ge e^{-T}$. Below the clipping threshold the clipped
logarithm is constant. The clipped logarithm is absolutely continuous
on this compact interval; a possible sign change at $U_j=0$ occurs
inside the constant region. Integration proves \eqref{S:eq:window-clip}.

\paragraph{Completing the upper bound.}
We have $\sum_jr_j(z)\ge1$ and
$\sum_j(\widetilde r_j(z)-r_j(z))\le ne^{-T}=\delta$. Thus \eqref{S:eq:window-clip}
implies $A_i\le(1+\delta)e^\eta/\sum_j\widetilde r_j(1)$.
On $\{G\ge t\}$ it also gives
$z\le e^\eta\max_k d_k(t)\widetilde r_{[k]}(1)$.
Integrating these two bounds first over $z\in[\delta,L]$ and then
over $t\in(0,1]$, with the good residual vector fixed, bounds the
integral by $(1+\delta)e^{2\eta}I(x)$, where
$x_k=\widetilde r_{[k]}(1)/\sum_j\widetilde r_j(1)$.
The bound is independent of the residual vector and the anchor sign.
Restoring the three discarded contributions yields
\begin{equation}\label{S:eq:window-upper}
 \frac{W_{n,s}(q)}{n\tau}
 \le (1+\delta)e^{2\eta}H_{s+1}+3\delta .
\end{equation}

\paragraph{Attainment to first order.}
Take independent uniform true-null $p$-values, realized by independent
Gaussian statistics with the prescribed sides. Whenever exactly one
true-null value lies in $((j-1)\tau,j\tau]$ for some $1\le j\le s+1$
and all other true-null values exceed $q$, set $j-1$ alternatives
to zero and the rest to one. These disjoint events yield FDP $1/j$.
Consequently for every $q\in(0,1)$,
\begin{equation}\label{S:eq:window-lower}
 \frac{W_{n,s}(q)}{n\tau}
 \ge (1-q)^{n-1}H_{s+1}.
\end{equation}
For $n=1$ this attains the exact single-null formula of \citet[Theorem~5.2]{LyuWang}. For general $n$, the lower construction uses
the permitted dependence of the alternative values on all true nulls.

Finally put $D=\log(1/q)$ and take $\delta=\exp(-D^{1/4})$. Under the
assumption on $m$, the quantities $1+T+M^2$ and $\log(L/\delta)$
are $o(\sqrt D)$, $\log(1/(\tau L))=D+o(D)$ and the normal Mills
inequalities give $y_*^2\sim2D$. Thus \eqref{S:eq:window-upper} applies for small enough
$q$ and $\eta\to0$. Since $nq\to0$ and $H_{s+1}\ge1$, \eqref{S:eq:window-upper} and \eqref{S:eq:window-lower}
prove the asserted limit.

\section{A bound on the rank of the null covariance}\label{S:sec:purerank}

Let $W_d(q)$ be the supremum of the FDR of ordinary BH at nominal
level $q$ over all finite numbers $n\ge1$ of true nulls and $s\ge0$
of alternatives, where the true-null statistics are jointly centered
Gaussian with unit variances and covariance rank at most $d$. The alternative $p$-values may depend arbitrarily on the true nulls. Thus $W_d$ differs from the fully Gaussian factor-model supremum $B_{d,b}$ in Section~\ref{sec:multiple-factors} in both the covariance constraint and the range of alternative laws.
Each null can use a prespecified one-sided or two-sided test. The
theorem holds separately when all tests are one-sided, all are
two-sided, or the choices of sides are mixed.

\begin{theorem}[Sharp rate at fixed null rank]\label{S:thm:purerank}
For each fixed integer $d\ge2$,
\[
 \lim_{q\downarrow0}
 \frac{W_d(q)}{q\log\log(1/q)}=\frac{d-1}{2}.
\]
For $d=1$ and exclusively two-sided tests, $W_1(q)=q$ for every
$0<q<1$. \end{theorem}
The rank constraint applies to the complete covariance of
the true-null statistics, including any independent noise terms.

The proof gives a geometric interpretation of the coefficient. Represent each null statistic as a projection of $G\sim N_d(0,I_d)$. An upper cutoff on $\|G\|$ limits how small any null $p$-value can be. For a lower bound, choose many directions on the unit sphere with sufficiently small joint exceedance probabilities. At the relevant Gaussian tail scale $u=\bar\Phi^{-1}(q)$, one can use on the order of $(u/\log u)^{d-1}$ such directions. The logarithm of their number is asymptotic to $(d-1)\log u$, which accounts for the constant in the theorem. Section~\ref{S:app:rank} gives both bounds.

\subsection{Proof of Theorem~\ref{S:thm:purerank}}\label{S:app:rank}

\paragraph{Upper bound.}
Put $m=n+s$ and $\pi_0=n/m$. There are unit vectors
$a_1,\ldots,a_n\in\mathbb R^d$ and $G\sim N(0,I_d)$ such that the
null statistics are $Z_i=a_i^\top G$. Write $\bar\Phi$ for the
standard normal survival function, orienting each one-sided test to
use the upper tail. If $\|G\|\le v$, then both the one-sided
$p$-value $\bar\Phi(Z_i)$ and the two-sided $p$-value
$2\bar\Phi(|Z_i|)$ are at least $\bar\Phi(v)$.

If BH makes $R>0$ rejections, every rejected null satisfies
$P_i\le qR/m\le q$ and therefore $1/R\le q/(mP_i)$. Choose $v$
so that $\Pr(\|G\|>v)=q$. On $\{\|G\|>v\}$ use $\mathrm{FDP}\le1$;
on its complement, apply the preceding two inequalities to each
rejected null. Since each null $p$-value is uniform and
$\bar\Phi(v)\le q$, we obtain
\[
 \mathrm{FDR}\le q+\frac qm\sum_{i=1}^n
 \mathbb E\!\left[\frac{{\bf1}\{\bar\Phi(v)\le P_i\le q\}}{P_i}\right]
 =q+\pi_0 q\log\frac q{\bar\Phi(v)} .
\]
For fixed $d$, integrating the $\chi_d$ density in the upper tail gives
$q\sim [2^{1-d/2}/\Gamma(d/2)]v^{d-2}e^{-v^2/2}$. By the normal
Mills asymptotic $\bar\Phi(v)\sim (2\pi)^{-1/2}v^{-1}e^{-v^2/2}$.
Hence $q/\bar\Phi(v)$ is asymptotic to a positive constant depending
only on $d$ times $v^{d-1}$. Also
$v\sim u:=\bar\Phi^{-1}(q)\sim\sqrt{2\log(1/q)}$. In particular,
uniformly in the numbers of tests and the rank-constrained null laws,
\[
 \mathrm{FDR}\le \pi_0(d-1)q\log u+O_d(q).
\]

\paragraph{Lower bound.}
We give the construction for every sufficiently small $q$, with
$d\ge2$ fixed. In a spherical patch of angular radius less than
$\pi/8$ in $S^{d-1}$ choose directions $a_1,\ldots,a_n$ with
pairwise angular distance at least $h=(\log u)/u$ and with
$n\asymp_d h^{-(d-1)}$. For example, take a maximal $h$-separated
set in a closed patch of fixed positive spherical area. Caps of
radius $h$ cover the patch, whereas caps of radius $h/2$ centered
at the chosen points have disjoint interiors. The spherical area of
a cap of radius $h$ is comparable to $h^{d-1}$ for small $h$,
giving the stated cardinality. The patch diameter ensures
$|a_i^\top a_j|\le\cos h\le1-c h^2$ for $i\ne j$, with an absolute
$c>0$. Let $Z_i=a_i^\top G$.

For either choice of sides, $P_i\le q$ implies $|Z_i|\ge u$. For
$i\ne j$ write $\rho=a_i^\top a_j$. For each of the four choices
of signs, simultaneous signed exceedances of $u$ imply that the
signed sum exceeds $2u$ and its variance is at most $2(1+|\rho|)$.
Consequently
\[
 \Pr(P_i\le q,\ P_j\le q)
 \le4\bar\Phi\left(u\sqrt{2/(1+|\rho|)}\right)
 \le4q\exp(-c' u^2h^2)
\]
for an absolute $c'>0$. The second inequality uses
$\bar\Phi(au)/\bar\Phi(u)\le\exp[-(a^2-1)u^2/2]$ for $a\ge1$,
$u>0$, obtained by integrating the normal hazard bound
$\phi(x)/\bar\Phi(x)\ge x$ from $u$ to $au$. Because $n$ grows
only polynomially in $u$ whereas $u^2h^2=(\log u)^2$, the sum of
the displayed probabilities over ordered pairs $i\ne j$ is $o(q)$.

Take $s=\lfloor n/\log n\rfloor$ and $\tau=q/(n+s)$. If exactly
one true-null $p$-value belongs to $((j-1)\tau,j\tau]$ with
$1\le j\le s+1$ and all other null $p$-values exceed $q$, set
$j-1$ alternative $p$-values to zero and the others to one. Set all
alternatives to one outside these disjoint events. On each indicated
event BH rejects exactly $j$ hypotheses and the FDP equals $1/j$.
With $H_r=\sum_{j=1}^r1/j$, uniform null marginals therefore give
\[
 \mathrm{FDR}\ge n\tau H_{s+1}
 -\sum_{i\ne l}\Pr(P_i\le q,\ P_l\le q)
 =\pi_0 q H_{s+1}-o(q).
\]
Indeed, for each $i$ the bins are disjoint subsets of $(0,q]$ and
each loss has coefficient $1/j\le1$. Here $\pi_0\to1$ and
$H_{s+1}\sim\log s\sim(d-1)\log u$, proving the asserted matching
lower bound. This construction uses alternatives chosen after
observing the nulls. For a prescribed fixed limiting null fraction
$0<\pi<1$, the same construction with $s\sim n(1-\pi)/\pi$ gives
the corresponding coefficient $\pi(d-1)/2$.

\paragraph{Two-sided rank one.}
Every two-sided null $p$-value is the same uniform random variable
$U$ because $Z_i=\pm G$ for a scalar $G$. If any null is rejected,
all $n$ are rejected. On $U\in((j-1)q/m,jq/m]$, the rejection
count must then be at least $\max(n,j)$; on $U>q$ no null can be
rejected. Thus
\[
 \mathrm{FDR}\le\frac{nq}{m}\{1+H_m-H_n\}
 \le\frac{nq}{m}\{1+\log(m/n)\}\le q.
\]
Taking $s=0$ attains equality. The finite-$n,s$ upper bound for a
tied null block is also given by \citet{LyuWang} for repeated null $p$-values.

\section{Sharp FDR with correlations decaying in the index distance}
\label{S:sec:decay-sharp}\label{S:app:benefit-decay}

Fix $D,\epsilon>0$. Write $\mathcal D_n(D,\epsilon)$ for the correlation matrices
of the ordered true-null coordinates satisfying
\begin{equation}\label{S:eq:logdecay-up}
 |\Sigma_{ij}|\le D[\log(e+|i-j|)]^{-1-\epsilon},\qquad i\ne j.
\end{equation}
The class permits arbitrary signs and perfect local correlation if the
envelope permits it. A principal submatrix of a matrix satisfying
\eqref{S:eq:logdecay-up} satisfies the same envelope when its indices are
listed in their original order.
Let $W^{\rm dec}_{n,s}(q)$ denote the worst FDR in
Theorem~\ref{S:thm:finite}, with the additional condition
$\Sigma\in\mathcal D_n(D,\epsilon)$ on the true nulls. The sides of the individual
tests are fixed in advance. The following limit also holds uniformly
over these assignments.

\begin{theorem}\label{S:thm:logdecay-sharp}
Put $a=nq/(n+s)$ and $L=\min\{s+1,1/a\}$. Then
\[
 \lim_{q\downarrow0}\ \sup_{n\ge1,\ s\ge0}
 \left|\frac{W^{\rm dec}_{n,s}(q)}
 {a H_{\lfloor L\rfloor}}-1\right|=0 .
\]
In particular, under the global null the supremum of the Simes
rejection probability is $q(1+o(1))$ uniformly in dimension. For fixed
$s$, the coefficient of $a$ is $H_{s+1}$ at every dimension sequence.
If $s$ grows without bound, the leading term is $a\log L$.
\end{theorem}

The FDR upper bound also applies to jointly Gaussian alternatives.
Attainment in the theorem uses alternatives that may depend on the
true-null values. The proof combines the Gaussian count estimates
of \citet{ShangGuo} with local control of the small order statistics.

\begin{corollary}\label{S:cor:decay-rate}
With the notation of Theorem~\ref{S:thm:logdecay-sharp}, write
$u=\log(1/a)$. There are constants $C,u_0$ depending only on
$D,\epsilon$ such that for every $0<q<1$, $n\ge1$ and $s\ge0$
with $u\ge u_0$,
\[
 \left|\frac{W^{\rm dec}_{n,s}(q)}{aH_{\lfloor L\rfloor}}-1\right|
 \le C u^{-1/3}.
\]
Let $A(x)$ be the supremum of the global-null size at level $x$ over these null configurations and all dimensions. Under the global null the stronger estimate
$0\le A(x)/x-1\le C[\log(e+v)]^2/v$ holds for
$v=\log(1/x)\ge u_0$.
\end{corollary}

\subsection{Pointwise comparison with the bound for decaying correlations}
\label{S:sec:sg-comparison}

Write $H_m=\sum_{k=1}^m k^{-1}$ and $p=\pi_0$. Both the harmonic bound of \citet{BY} and Theorem~\ref{thm:main} apply to the same ordinary BH procedure. Combining them gives an explicit bound that does not involve correlation-decay parameters.

\begin{proposition}\label{S:prop:sg-dom}
Let $C=e$ for two-sided tests and $C=2e$ for one-sided tests. Throughout the range $0<q\le1/2$ of \citet[Theorem~1(i), equation~(8)]{ShangGuo}, the bound
\[
 U_m(q,p)=\min\{pqH_m,\ B(Cpq)\}
\]
is at most their stated bound, for every admissible covariance class and every choice of its tuning parameters. The inequality is strict if $p>0$, $Cpq<1$, and $B(Cpq)<pqH_m$.
\end{proposition}
\begin{proof}
For $m\ge2$ their second branch is $q$ times an infimum of the expression below, or a special case with fixed $K,E$.
\[
 \Psi_p(K,E,q)=2p+\int_{2p}^{1/q}
 \min\!\left\{1,\frac{6p\sqrt K}{t-2p}+\frac{4pE}{(t-2p)^2}\right\}dt ,
 \qquad K\ge1,\ E\ge0 .
\]
Their first branch is $pqH_m$. For $p>0$, equation~\eqref{eq:fdr} gives
\[B(Cpq)/q=\int_0^{1/q}\min\{1,Cp/t\}\,dt.\] For $t\le2p$ the integrand is at most one. On $t>2p$, $6p/(t-2p)>Cp/t$ since $C\le2e<6$. Thus $B(Cpq)\le q\Psi_p(1,0,q)$ even before optimizing the second branch. For $Cpq<1$ strictness holds on a set of positive measure in $t>Cp$; here $C>2$. If $m=1$ both stated bounds agree with $pq$ and if $p=0$ both vanish.
\end{proof}

For example, with $m=10^5$, $p=0.8$ and $q=0.05$, the two-sided bound is $0.34999$. The harmonic bound is $0.48361$. The second branch of their bound is at least $q\Psi_p(1,0,q)=0.64250$, so its minimum is exactly the harmonic bound. These values compare upper bounds for the same parameters. The correlation-specific FDR can be much smaller.

\subsection{Proof of Theorem~\ref{S:thm:logdecay-sharp}}
\begin{proof}
We work at the boundary means. Decreasing a one-sided true-null mean
after orienting that test to the right can only increase its $p$-value.
The maximal FDP conditional on the null values, given by the finite
oracle in Section~\ref{S:app:window}, is coordinatewise nonincreasing.
This justifies the same upper bound for composite one-sided nulls.

\emph{Large witness ranks.}
For a centered null sample of size $n$, write $P_{(k)}$ for its order
statistics. There are constants $C$ and an integer $d$ depending only
on $D,\epsilon$ such that, uniformly in the dimension, for sufficiently
small $x>0$ and integers $R\ge4$,
\begin{equation}\label{S:eq:decay-large-ranks}
 \Pr\{P_{(k)}\le xk/n\ \text{for some }k\ge R\}
 \le Cx/R+Cx^2[1+\log(1/x)]^d .
\end{equation}
For homogeneous one- or two-sided tests, take $w_j=x2^{-j}$,
$c_j=nw_j$, $0\le j<\lceil\log_2 n\rceil$, and
$N(w)=\#\{i:P_i\le w\}$. Any witness with $k\ge R$ yields a grid
point $xk/n\le w_j<2xk/n$ with $c_j\ge xR$ and
$N(w_j)\ge c_j/(2x)$. If $n<R$ the event is empty.
Use the common low-order Hermite event and the individual remainder
events in \citet[proof of Proposition~2, Steps 1--2]{ShangGuo}, taking
their $\pi_0=1$, $\theta=1/2$ and $t=1/x$. In their notation, for
any subset $J'$ of the dyadic indices, these bounds give
\[
 \Pr\{N(w_j)\ge c_j/(2x)\text{ for some }j\in J'\}
 \le C x^2\left\{
  \Lambda_{\rho,n}(x)+
  \sum_j e_{\rho,n}(z_j)+K_{\rho}\sum_{j\in J'}c_j^{-1}
 \right\}.
\]
Indeed the common event is charged once, and the variance estimate for
the $j$th remainder contributes $C x^2(K_{\rho}/c_j+e_{\rho,n}(z_j))$.
Their Section~4.5, equations (27)--(28) and Step 4, bound the first
two terms in braces by $C[1+\log(1/x)]^d$ and $K_{\rho}$ by $C$,
uniformly over the envelope. One may take an even integer
$d\ge \max\{2,1/\epsilon\}$. For $J'=\{j:c_j\ge xR\}$ a
geometric sum gives $\sum_{j\in J'}c_j^{-1}\le2/(xR)$.
This proves \eqref{S:eq:decay-large-ranks} in each homogeneous family.
A left-sided orientation changes only correlation signs. In a mixed
family, some homogeneous group of size $n_g$ contains at least half
the $k$ witness coordinates. If their number in this group is $l$,
the group has a witness at rank $l\ge\lceil R/2\rceil$ with level
$x_g=2xn_g/n$, because $xk/n\le x_g l/n_g$. Apply the homogeneous
bound to the two nonempty groups and sum, using $\sum_g x_g=2x$.
The same proof for homogeneous tests works already for $R\ge2$.
For small $x$, $u^2[1+\log(1/u)]^d$ is increasing on $0<u\le 2x$,
which controls the group remainders.

\emph{Localization of small values.}
Choose an integer $h$ so that the right side of \eqref{S:eq:logdecay-up}
is at most $1/4$ for lags at least $h$. For each fixed $M\ge1$,
uniformly in $n$ and $\Sigma\in\mathcal D_n(D,\epsilon)$,
\begin{equation}\label{S:eq:decay-pair}
 \Pr\{P_i,P_j\le Mx/n\ \text{for some }|i-j|\ge h\}=O_M(x^{8/5}).
\end{equation}
For two coordinates with absolute correlation $r\le1/4$ and small $t$,
a union over the four sign pairs gives
$\Pr(P_i,P_j\le t)\le C t^{2/(1+r)}$. To see this, set
$u=\bar\Phi^{-1}(t)\ge1$; every selected signed coordinate exceeds
$u$. The bivariate Gaussian density in any such quadrant is bounded
by $C\exp[-(z_1^2+z_2^2)/(2(1+r))]$. Integration bounds each quadrant
by $C u^{-2}\exp[-u^2/(1+r)]$; the lower Mills inequality and
$u^{-2+2/(1+r)}\le1$ give the assertion. For
$h\le |i-j|<\sqrt n$ there are at most $n^{3/2}$ unordered pairs
and $2/(1+r)\ge8/5$. For $|i-j|\ge\max(h,\sqrt n)$, the envelope
gives $r\log n\le C$. Thus $n^{2-2/(1+r)}\le e^{2r\log n}\le C$,
and the $n^2$ pairs contribute at most $C(Mx)^{8/5}$ for $Mx<1$.
This proves \eqref{S:eq:decay-pair}.

\emph{Finitely many alternatives.}
Fix $s$ and let $\tau=q/(n+s)=a/n$. Take integers $B>h$ and
$R\ge\max(4,s)$. If the BH outcome has $k\ge R$ rejected true nulls,
then $P_{(k)}\le(k+s)\tau\le a(1+s/R)k/n$. Consequently
\eqref{S:eq:decay-large-ranks} bounds this event by
$Ca/R+o(a)$, uniformly in $n$ for fixed $R,s$.
If $1\le k<R$, every rejected true null belongs to
$T=\{i:P_i\le(R+s)\tau\}$. By \eqref{S:eq:decay-pair}, except on an
event of probability $o(a)$, the diameter of $T$ is less than $h$
when $T$ is nonempty.

Partition the integer line into successive blocks of $B$ indices with
an independent uniform random shift in $\{0,\ldots,B-1\}$, and intersect
each block with $\{1,\ldots,n\}$. Given the data on this last event,
the probability that $T$ is contained in one block is at least $1-h/B$.
For a nonempty block $b$ let
\[
 G_b=\max_{1\le l\le |b|,\ 0\le j\le s}
       \frac{l}{l+j}\,
       {\bf1}\{P^b_{(l)}\le(l+j)\tau\},
\]
where $P^b_{(l)}$ is its $l$th order statistic.
When $1\le k<R$ and $T$ is inside one block, the global FDP is at
most that block's $G_b$: use $l=k$ and let $j$ be the number of
rejected alternatives. Conditional on the partition,
Theorem~\ref{S:thm:window} for $|b|$ true nulls, $s$ alternatives, and
level $q_b=(|b|+s)\tau\le q$ implies
$\mathbb E G_b\le(1+o(1))|b|\tau H_{s+1}$. The oracle equality in
Section~\ref{S:app:window} justifies applying the theorem to $G_b$.
The error is uniform for $0<q_b\le q$, because there are only finitely
many block sizes and side assignments for fixed $B,s$. Summing over
blocks and averaging over the shift gives
\[
 W^{\rm dec}_{n,s}(q)\le
  \frac{(1+o(1))aH_{s+1}}{1-h/B}+Ca/R+o(a).
\]
Take $q\downarrow0$, then $B,R\to\infty$. The independent true-null
construction in \eqref{S:eq:iidfixed} supplies a matching lower bound
uniform in $n$. This proves the result for each fixed $s$, including
the global-null assertion.

\emph{Growing numbers of alternatives.}
Let $A(t)$ be the supremum of the global-null Simes rejection
probability over dimensions, side assignments and the covariance class;
set $A(t)=1$ for $t\ge1$. The case $s=0$ gives $A(t)/t\to1$.
In particular $A(t)\le C t$ for $t>0$. For any $\eta>0$ choose
$0<\delta<1$ so that $A(t)\le(1+\eta)t$ for $0<t\le\delta$.
Su's pathwise comparison gives $\Pr(\FDP\ge y)\le A(a/y)$ for $y>0$.
Put $y_0=1/(s+1)$. The contribution to $\mathbb E\FDP$ from
integrating its upper tail on $(0,y_0)$ is at most
$y_0 A(a/y_0)\le Ca$, since positive FDP is at least $y_0$.
The part of $[y_0,1]$ with $y<a/\delta$ contributes at most
$a/\delta$; the remainder contributes at most
$(1+\eta)a\log L$, by integration of $1/y$ from
$\max(y_0,a/\delta)$ to one when this interval is nonempty.
Therefore when $L\to\infty$ the upper leading coefficient is at most
one, on letting $\eta\downarrow0$. Equation~\eqref{S:eq:iid} proves the
reverse bound. Finally $a\le q$, so along any sequence $q\downarrow0$
one can extract a subsequence on which $s$ is constant or tends to
infinity. In the latter case $L\to\infty$ and
$H_{\lfloor L\rfloor}\sim\log L$. Applying the two cases to any
putative sequence violating the uniform limit proves the theorem.
\end{proof}

\subsection{A quantitative remainder under correlation decay}
\begin{proof}
First suppose $\log(s+1)\le u^{1/3}$.
In the localization proof choose $B=\lceil u\rceil$ and
$R=\max\{B,s,4\}$. For all large $u$, $B>2h$. The large-rank
probability divided by $a$ is $O(1/u+a(1+u)^d)$.
The proof of \eqref{S:eq:decay-pair} bounds the far-pair probability
by $C[(R+s)a]^{8/5}$ as soon as $(R+s)a$ is smaller than an
absolute constant, with $C$ independent of $R,s$. Dividing by $a$
gives an exponentially decreasing bound in $u$, since
$R+s\le 2(\lceil u\rceil+e^{u^{1/3}})$.

For each block of size $b\le B$ use \eqref{S:eq:window-upper}
with local total dimension $m_b=b+s$, common threshold increment
$\tau=a/n$, and $\delta=1/u$. In the notation of its proof the local
cutoff is $L_b=bm_b^2u$, so
$\log(L_b/\delta)$ and $1+T+M^2$ are bounded by
$C\{\log(e+u)+\log(s+1)\}$. Also
$\log(1/(\tau L_b))\ge u-\log L_b\ge u/2$ for large $u$.
The Mills inequalities give $y_*^2\ge c u$ and verify the
condition for \eqref{S:eq:window-upper}. Its relative error, including
the discarded $3\delta$ term, is at most
\[
 C\,\frac{[\log(e+u)+\log(s+1)]^2}{u}.
\]
Averaging over the partition costs a further factor $1+O(1/u)$.
Here $L=s+1$. The construction in \eqref{S:eq:lower-sum} and
Bernoulli's inequality give
$W^{\rm dec}_{n,s}(q)\ge aH_{s+1}-a^2(s+1)$.
These inequalities prove the claimed $O(u^{-1/3})$ bound when
$\log(s+1)\le u^{1/3}$ and the sharper assertion when $s=0$.
In particular $A(x)\le Cx$ for all $x>0$, using the convention
$A(x)=1$ for $x\ge1$.

For the remaining case $\log(s+1)>u^{1/3}$, let $v=\log L$,
so $v\ge u^{1/3}$. Integrating Su's comparison as in the theorem
now gives
\[
 W^{\rm dec}_{n,s}(q)\le a\log L+Ca+
 Ca\int_{u-v}^{u}\frac{[\log(e+w)]^2}{1+w}\,dw .
\]
Indeed, put $y_0=1/(s+1)$. The part below $y_0$ is at most $Ca$ and
the part of $[y_0,1]$ below $a$ is at most $a$. On the remaining
interval change variables $w=\log(y/a)$; enlargement of $C$ extends
the global-null relative error estimate to $w\ge0$.
When $v\le u/2$ the integral is at most
$Cv[\log(e+u)]^2/u$. When $v>u/2$ it is at most
$C[\log(e+u)]^3$. In both cases the upper bound divided by
$aH_{\lfloor L\rfloor}$ exceeds one by at most $C u^{-1/3}$,
using $H_{\lfloor L\rfloor}=v+O(1)$ and $v\ge u^{1/3}$.
Equation~\eqref{S:eq:iid} gives the lower bound with relative
error $O(1/v)$. This proves the corollary.
\end{proof}

\section{Exact control for two-sided Gaussian Simes tests}\label{S:sec:planar-main}

We now return to the global null and consider the sharper question $\Prob\{\Simes(P)\le q\}\le q$. In a fixed dimension, \citet[Theorem~2.1]{Zhang} proves the inequality uniformly over correlations at sufficiently small levels and characterizes equality by independent groups of signed copies. For growing dimension, Zhang obtains the unit constant on $0<q\le a$ for any fixed $a<1$ when covariance graphs have bounded degree or when independent block sizes meet the growth condition in \citet[Theorems~5.2 and~6.1]{Zhang}. The model below has a dense covariance graph and admits an explicit representation of its rejection probability.

Let $u_j=(\cos(j\pi/m),\sin(j\pi/m))$ for $0\le j<m$.
For $G\sim N_2(0,I_2)$ put $P_j=S(|u_j^\top G|)$, where
$S(z)=2\bar\Phi(z)$, $z\ge0$, and write $\zeta=S^{-1}$.
Thus every test uses the same two Gaussian coordinates. The
correlation of tests $i$ and $j$ before taking absolute values is
$\cos((i-j)\pi/m)$.

\begin{theorem}[Regular planar directions]\label{S:thm:polygon}
There is an integer $M$ such that for every $m\ge M$ and
every $0<\alpha<1$,
\[
 \Prob\{\Simes(P)\le\alpha\}<\alpha .
\]
For each such $m$, the same strict inequality holds at every level in
a relative neighborhood of the polygon correlation matrix among
$m$-dimensional correlation matrices.
\end{theorem}

These correlations stay large over many nearby directions as $m$ grows and each test can be correlated with almost every other test. The designs violate positive regression dependence. The neighborhood statement includes dense positive definite examples. The proof gives existence of $M$ and of the neighborhoods.

Writing $\vartheta$ for the angle of $G$, let $c_k(\vartheta)$ be the $k$th largest of $|\cos(\vartheta-j\pi/m)|$, $0\le j<m$. Rotational symmetry and the Rayleigh tail of $\|G\|$ give the finite-dimensional formula
\[
 \Pr\{\Simes(P)\le\alpha\}=\frac{2m}{\pi}\int_0^{\pi/(2m)}
 \exp\!\left[-\frac12\min_{1\le k\le m}
   \left(\frac{\zeta(\alpha k/m)}{c_k(\vartheta)}\right)^2\right]d\vartheta .
\]

A zero denominator gives an infinite ratio. The integral expresses the size as an average over the angle of a single Gaussian vector. Section~\ref{S:app:polygon} controls it uniformly near both endpoints of $(0,1)$ and on compact subintervals. The neighborhood assertion uses the covariance-uniform small-level theorem of \citet{Zhang}.

\subsection{Uniform exactness in a neighborhood of independence}

Write $F_\Sigma(q)$ for the rejection probability of the global-null Simes test applied to $P_i=2\bar\Phi(|Z_i|)$, $Z\sim N_m(0,\Sigma)$. For a fixed dimension, the next result gives a single neighborhood of independence on which exact control holds at every level. Its radius can depend on $m$.

\begin{theorem}\label{S:thm:local-all-level}
For every $m\ge2$ there is $\delta_m>0$ such that
\[
   0<\|\Sigma-I_m\|_{\rm F}<\delta_m
   \quad\Longrightarrow\quad
   F_\Sigma(q)<q \qquad (0<q<1).
\]
For every $0<a<1$, uniformly for $a\le q<1$ and $\Sigma=I_m+E$ near $I_m$,
\[
   F_\Sigma(q)=q-D_m(q)\sum_{i<j}E_{ij}^2
   +O_{m,a}\!\left((1-q)\|E\|_{\rm F}^3\right),
   \qquad \inf_{a\le q<1}\frac{D_m(q)}{1-q}>0 .
\]
\end{theorem}

\subsection{The upper endpoint in three dimensions}

The acceptance probability close to level one gives a second comparison with the low-level theorem of \citet{Zhang}. Write $F_\Sigma(q)$ as above. For $m\ge2$ and a matrix with no perfectly correlated pair, define
\[
 A_i(\Sigma)=\Prob\!\left\{P^{(-i)}_{(k)}> k/m,\ 1\le k\le m-1
                   \,\middle|\, Z_i=0\right\},
\]
where the remaining $p$-values retain their original marginal standardization.

\begin{proposition}\label{S:prop:trivariate-upper-end}
For each fixed $m$ and $\epsilon>0$, uniformly on the correlation matrices satisfying $\max_{i<j}|\Sigma_{ij}|\le1-\epsilon$,
\[
 \frac{1-F_\Sigma(q)}{1-q}
     =\sum_{i=1}^m A_i(\Sigma)+o(1),\qquad q\uparrow1 .
\]
For $m=3$ the sum is at least one, and it is strictly larger than one for every $\Sigma\ne I_3$ in this class. Consequently there exists $b_\epsilon<1$ such that, uniformly on this class, $F_\Sigma(q)\le q$ for $b_\epsilon<q<1$, with equality only at independence.
\end{proposition}

\begin{corollary}\label{S:cor:rank-two-upper-end}
Fix $m\ge2$ and $\epsilon>0$. Among correlation matrices of rank at most
two with $\max_{i<j}|\Sigma_{ij}|\le1-\epsilon$, there exists a common
$b_{m,\epsilon}<1$ such that $F_\Sigma(q)\le q$ for all
$b_{m,\epsilon}<q<1$, with strict inequality except at independence.
\end{corollary}

\subsection{Open neighborhoods of Gaussian PRDS models}

\begin{theorem}\label{S:thm:prds-interior}
Fix $m\ge2$. Let $\Sigma_0$ be a positive definite correlation matrix
whose folded Gaussian $p$-values are PRDS. There exists an open
neighborhood $\mathcal U$ of $\Sigma_0$ in the space of correlation
matrices such that every $\Sigma\in\mathcal U$ satisfies
$F_\Sigma(q)\le q$ simultaneously for $0<q<1$. The inequality is
strict at every level for every $\Sigma\in\mathcal U\setminus\{I_m\}$.
The neighboring matrices are subject only to membership in $\mathcal U$.
\end{theorem}

\begin{corollary}\label{S:cor:non-prds-interior}
For every $m\ge3$ there are open sets of positive definite, fully
dense correlation matrices for which Simes is strictly conservative
at every level although the folded $p$-values fail PRDS. An off-diagonal
correlation can be chosen arbitrarily close to one.
\end{corollary}

\subsection{Dense spherical designs in arbitrary fixed rank}
Let $d\ge2$ be fixed, $G\sim N_d(0,I_d)$, and let
$u_{m,1},\ldots,u_{m,m}\in S^{d-1}$ be deterministic unit vectors.
Test the global null with $P_i=S(|u_{m,i}^{\top}G|)$, where
$S(z)=2\bar\Phi(z)$, and write $S_m$ for their Simes $p$-value.
For uniform $U\in S^{d-1}$, let $D_d(x)=\Pr(|U_1|\ge x)$.
Define the discrepancy of the double caps by
\[
 \Delta_m=\sup_{v\in S^{d-1},\,0\le x\le1}
 \left|m^{-1}\sum_{i=1}^m\mathbf 1\{|u_{m,i}^{\top}v|\ge x\}-D_d(x)\right|.
\]
Put $h=(d-1)/2$ and $\kappa_d=(h/e)^h/\Gamma(h+1)$.
\begin{theorem}[Uniform Simes validity for dense spherical designs]
\label{S:thm:spherical-design}
There is a constant $C_d<\infty$ such that, if $0<q_m<1$ and
$\Delta_m[1+\log(1/q_m)]^h\to0$, then for all sufficiently large $m$
\[
             \Pr(S_m\le q)<q\qquad(q_m\le q\le 1-C_d\Delta_m).
\]
If additionally $q_m\to0$, then
$\Pr(S_m\le q_m)/q_m\to\kappa_d$. The convergence and the sufficient
dimension threshold are uniform over configurations obeying a common
upper bound on the discrepancy with the stated decay.
\end{theorem}
Here $\kappa_2=\sqrt{2/(\pi e)}$ and $\kappa_3=1/e$.
In particular this gives strict nominal control, at every fixed
conventional range of levels, for arbitrary asymptotically uniform
spherical designs. Their covariance graphs may be complete.

For example, deterministic arrays with
$\Delta_m=O_d(\sqrt{\log(m)/m})$ exist. To see this, sample independent
uniform directions and use nets of angular mesh $1/m$ for the cap
axis and angular radius. There are polynomially many net caps for
fixed $d$. Hoeffding's inequality, a union bound and Borel--Cantelli
give the asserted error on the nets, for a sufficiently large fixed
constant. Enlarging and shrinking each net cap by $O(1/m)$ in angle
brackets every cap; its spherical mass changes by $O_d(1/m)$ because
the angle density on $[0,\pi/2]$ is proportional to $\sin^{d-2}$.
Almost every resulting array has a complete covariance graph.
For these arrays one may allow
$\log(1/q_m)=o((m/\log m)^{1/(d-1)})$.

\subsection{Proof of the planar Simes inequality}\label{S:app:polygon}

\begin{proof}
Let $R=\|G\|$. Its direction is independent and uniform, and
$\Prob(R>r)=e^{-r^2/2}$. For a specified direction let
$d_1\le\cdots\le d_m$ be its angles to the $m$ lines, measured in
$[0,\pi/2]$. The $k$th largest absolute projection has size $R\cos d_k$.
Counting grid points on a projective arc of length $2d$ gives
\begin{equation}
 d_1\le\frac{\pi}{2m},\qquad
 d_k\ge t_k:=\frac{(k-1)\pi}{2m}\quad(1\le k\le m).
 \label{S:eq:polygon-spacing}
\end{equation}
Write $r_\alpha=\sqrt{-2\log\alpha}$. If no direction can reject at
radius $r_\alpha$, then the rejection probability is strictly below
$\alpha$. Indeed the finitely many continuous ordered projections are
strictly below their respective thresholds on the compact circle. The
same holds on a slightly larger circle, and all smaller radii satisfy
the same inequalities.

We show that every sequence $m\to\infty$ admits no sequence of failing
levels. Upon taking a subsequence suppose $\alpha=\alpha_m$ converges.

\emph{A limiting level in $(0,1)$.}
For $r,z\ge0$ let $D(r,z)=\Prob_\Theta(r|\cos\Theta|\ge z)$ for a
uniform angle $\Theta$. For fixed $r>0$, $D(r,z)=2\arccos(z/r)/\pi$
on $0\le z\le r$, and $D(r,z)=0$ for $z>r$. At a fixed positive
$r$ and every $0\le z\le r$,
\[
 S(z)=\E[D(R,z)]>\Prob(R\ge r)D(r,z).
\]
At $z=0$ the inequality reads $1>e^{-r^2/2}$. For $z>0$, $D(R,z)$
strictly increases for $R>r$. By continuity the difference has a
positive minimum on $0\le z\le r$. This gap is uniform when $r$ lies
in a sufficiently small compact neighborhood of
$\lim r_{\alpha_m}\in(0,\infty)$, by parameterizing $z=rs$,
$0\le s\le1$. For any direction, the fraction of grid lines with
projection at least $z$ differs from $D(r,z)$ by at most $1/m$,
by counting points in a single arc of the projective circle.
At $r=r_\alpha$ and all sufficiently large $m$ this fraction is
strictly smaller than $S(z)/\alpha$ for all $z\in[0,r]$.
There is consequently no rejection on this circle.

\emph{Levels tending to one.}
Set $c=\sqrt{2/\pi}$, choose $0<r_0<1/(\pi c)$, and put
$p_0=1-e^{-r_0^2/2}$. Since $1-S(z)\le cz$, we have
$\zeta(1-v)\ge v/c$. On $R\le r_0$, for $k=m-j$, $1\le j<m$,
\[
 R\cos d_k\le r_0\sin\frac{(j+1)\pi}{2m}
       <\frac{j}{cm}\le \zeta(\alpha k/m).
\]
Only rank $m$ can reject there. At radius $r\le r_0$, the event
$\min_j|u_j^\top G|<\zeta(\alpha)$ consists, up to endpoints, of
equally spaced arcs about the perpendiculars. Its fraction of the
circle is at least $\min\{1,2m \zeta(\alpha)/(\pi r_0)\}$, since each
half-width is $\arcsin(\min\{1,\zeta(\alpha)/r\})$.
With $x=1-\alpha$, the probability of acceptance is therefore at least
\[
 p_0\min\{1,2mx/(\pi r_0c)\}.
\]
It exceeds $x$ whenever $0<x<p_0$ and $2mp_0>\pi r_0c$.

\emph{Levels tending to zero with $r_\alpha/m\le 2/3$.}
Put $r=r_\alpha\to\infty$. We prove
$S(r\cos t_k)>k\alpha/m$ for every rank eventually.
For $k=1$, the normal tail asymptotic
$S(r)\sim c e^{-r^2/2}/r$ gives
$\limsup\alpha/(m S(r))\le 2/(3c)<1$.
For any putative sequence of failing ranks $k\ge2$, let
$v=r(k-1)/m$. The elementary lower bound
$S(y)\ge C e^{-y^2/2}/(1+y)$ for $y\ge0$ follows by integrating
the half-normal density from $y$ to $y+1/(1+y)$.
Since $\sin t_k\ge(k-1)/m$ and $k/(k-1)\le2$, if $v\to\infty$ then
\[
 \frac{k\alpha}{mS(r\cos t_k)}
 \le C'v(1+1/r)e^{-v^2/2}\longrightarrow0 .
\]
Otherwise pass to a subsequence on which $v$ is bounded. Then
$t_k=\pi v/(2r)$, and the normal tail asymptotic gives
\[
 \frac{k\alpha}{mS(r\cos t_k)}
 =(1+o(1))\frac{1}{c}\frac{k}{k-1}\,
       v\exp(-\pi^2 v^2/8).
\]
Its limit superior is at most $4/(\pi c\sqrt e)=\sqrt{8/(\pi e)}<1$,
by maximizing over $v\ge0$. These bounds exclude every failing rank.
By \eqref{S:eq:polygon-spacing} they exclude every rejecting direction
at radius $r$.

\emph{Levels tending to zero with $r_\alpha/m\ge 2/3$.}
Put $A=\zeta(\alpha/m)$ and $B=\zeta(2\alpha/m)$. We first show that, for
all sufficiently large $m$ in this regime, a rejection at any rank
$k\ge3$ entails one at rank one. For $r=r_\alpha\to\infty$ and
every fixed $0<a<1$, the normal tail asymptotic gives
$S(ar)/\alpha\to\infty$. Thus $A>ar$ eventually and
$\liminf A/m\ge2/3$. Write $\eta=\pi/(2m)$ and
$b_k=\cos(t_k)/\cos\eta\in(0,1)$ for $k\ge3$.
The normal hazard is at least its argument on $[0,\infty)$, so
\[
 \frac{S(A b_k)}{S(A)}
 \ge\exp\{A^2(1-b_k^2)/2\}.
\]
The exponent exceeds $\log k$ for all $3\le k\le m$ eventually.
Indeed for any fixed $k\ge3$ its limit inferior is at least
$\pi^2 k(k-2)/18>\log k$, by Taylor expansion of sine at zero.
For $k\to\infty$ use
\[
 1-b_k^2=\frac{\sin^2 t_k-\sin^2\eta}{\cos^2\eta}
       \ge\frac{(k-1)^2-\pi^2/4}{m^2}.
\]
The resulting lower bound grows faster than $\log k$. The two cases
exclude any sequence of exceptional ranks. Hence
$\zeta(k\alpha/m)>A b_k$. If $R\cos d_k\ge \zeta(k\alpha/m)$, the spacing
bounds give $R\cos d_1>A$ (if $\cos d_k=0$ the premise is impossible).

It remains to compare the first two ranks. Use all $2m$ oriented
directions $u_j$, now with angles $j\pi/m$, $0\le j<2m$, and write
$X_j=u_j^\top G$ with circular indices. Put
$N_A=\#\{0\le j<m:|u_j^\top G|\ge A\}$. Exactly $N_A$ of the
$2m$ oriented projections satisfy $X_j\ge A$, and they are consecutive:
the superlevel set in angle is an arc of length strictly below $\pi$.
The number of grid edges whose endpoints both exceed $A$ is therefore
$(N_A-1)_+$. A rank-two rejection when $N_A=0$ has an adjacent pair
with both oriented projections in $[B,A)$: the two largest absolute
projections correspond to the two oriented directions nearest to $G$.
Angular ties have probability zero.

Each adjacent pair $(X,Y)$ is exchangeable bivariate normal with
correlation $\cos(\pi/m)>0$. Since $\Prob(X\ge B)=2\Prob(X\ge A)$
and $x\mapsto\Prob(Y\ge B\mid X=x)$ strictly increases,
\[
 \Prob(X\ge B,Y\ge B)<2\Prob(X\ge A,Y\ge B).
\]
Expanding the probability of $[B,A)^2$ and using exchangeability
gives $\Prob((X,Y)\in[B,A)^2)<\Prob(X\ge A,Y\ge A)$.
Summation over all $2m$ edges yields a rejection probability strictly
less than $\Prob(N_A>0)+\E[(N_A-1)_+]=\E N_A=\alpha$.
Together the four regimes show that no sequence $m\to\infty$ with a
failing level exists, proving the first assertion.

For the neighborhood assertion fix an $m$ for which the first assertion
holds and increase $M$ if necessary. The covariance-uniform small-level
theorem of \citet[Theorem~2.1]{Zhang} gives strict control for all
sufficiently small levels in a neighborhood of the polygon matrix:
its finitely many equality matrices (independent groups of signed
copies) exclude the polygon.
At the other endpoint let $\Sigma$ tend to the polygon matrix and
$\alpha\uparrow1$. Pairwise intersections of the central strips
$|Z_i|<\zeta(\alpha)$ have probability $O((1-\alpha)^2)$ uniformly
near that matrix, since all two-coordinate Gaussian densities are
uniformly bounded there for this fixed $m$. Conditional on an individual
strip, the remaining coordinates converge jointly to their Gaussian
law at the polygon given $Z_i=0$, by the Gaussian conditional formula
and $\zeta(\alpha)\to0$. In that limiting law every other coordinate has
positive variance, so the strictly positive threshold boundaries at
ranks $1,\ldots,m-1$ have probability zero. The residual latent vector
has the form $V u_i^\perp$ with $V\sim N(0,1)$. For $|V|<r_0$ all
those ranks strictly miss at level one by the bound for $k=m-j$ above.
Consequently the joint limit inferior of the acceptance probability
divided by $1-\alpha$ is at least $m\Prob(|V|<r_0)>1$ for all large
enough $m$, by adding the individual strip contributions and subtracting
the pairwise overlaps. On the remaining compact level interval,
strictness persists by joint continuity in $(\Sigma,\alpha)$:
a common Gaussian square-root representation gives convergence almost
surely and the union of threshold boundaries has probability zero.
This proves the second assertion.
\end{proof}

These designs lie outside positive regression dependence, including
the version defined by conditioning on a null $p$-value falling below
a cutoff. Choose two lines at angles $\pm\delta$ to $u_0$ with
$\pi/6<\delta<\pi/2$, and put $a=\cos\delta$, $b=\sin\delta$.
Given $V=|u_0^\top G|=v$, their absolute projections have the joint
law of $|av+bW|,|av-bW|$ for $W\sim N(0,1)$. Let $D$ be the event
that both exceed $t>0$. For $0\le v<t/a$,
$h(v)=\Prob(D\mid V=v)=S((t+av)/b)$. In particular
$h(0)=1-ct/b+O(t^3)$ as $t\downarrow0$, where $c=\sqrt{2/\pi}$.
The two projections have a bounded bivariate density, so
$\Prob(D)=1-2ct+O(t^2)$. Because $b>1/2$, for some fixed small $t$
we have $h(0)>\Prob(D)$. The right derivative at zero of
$v\mapsto\Prob(D\mid V\ge v)$ is therefore
$c\{\Prob(D)-h(0)\}<0$. This contradicts positive regression
dependence for the upper-tail $p$-values, by taking the increasing
$p$-value event $D^c$. The strict comparison at two fixed cutoffs
persists under small covariance perturbations by continuity of the
Gaussian probabilities. In particular, the neighborhoods in
Theorem~\ref{S:thm:polygon} include positive definite examples with every
off-diagonal correlation nonzero and strict Simes control at every level.

\subsection{Local and spherical results}\label{S:app:benefit-simes}
\subsection{Proof of Theorem~\ref{S:thm:local-all-level}}
\begin{proof}
Put $\tau=q/m$, $z_u=\bar\Phi^{-1}(u/2)$,
$L(u)=2z_u\phi(z_u)$ and $h(u)=L(u)/u$. For $0<u<1$, $h(u)=z_u\phi(z_u)/\bar\Phi(z_u)$ strictly decreases: $z_u$ decreases and the Gaussian hazard increases on $[0,\infty)$. Also $h(1)=0$ by continuity and $L(q)/(1-q)\to1$ as $q\uparrow1$.

We compute the second-order coefficient by giving just the first pair correlation $\rho$. Let $C_\rho(s,t)$ denote its folded-normal $p$-value copula. Given the remaining $m-2$ independent uniform $p$-values, let $K$ and $J$ be the ordinary BH rejection counts when $(P_1,P_2)$ is replaced by $(0,0)$ and $(0,1)$, respectively. For $q<1$,
$1\le J<K\le m$: taking the last feasible rank with $(0,1)$ and replacing 1 by 0 makes the next rank feasible. With only $P_1$ replaced by 0, the rejection count is $K$ when $P_2\le K\tau$ and is $J$ otherwise. This follows from the step-up rule: a rejected $p$-value can be set to zero, or an unrejected one to one, without changing the rejection set. Thus the contribution of hypothesis 1 to the FDR, conditional on the remaining $m-2$ values and using the usual leave-one-out identity, is
\[
 \tau-C_\rho(J\tau,K\tau)/J+C_\rho(K\tau,K\tau)/K.
\]
By exchangeability, the second hypothesis contributes the same. Each independent coordinate contributes $\tau$. Writing $F_{\rho}(q)$ for this one-pair model gives
\[
 q-F_{\rho}(q)
 =2\,\mathbb E\left[\frac{C_\rho(J\tau,K\tau)}{J}
                 -\frac{C_\rho(K\tau,K\tau)}{K}\right].
\]
The bivariate normal density at $\rho=0$ has first derivative $xy\phi(x)\phi(y)$ and second derivative $(x^2-1)(y^2-1)\phi(x)\phi(y)$. On the folded tail the first term integrates to zero and
$\int_{|x|\ge z_u}(x^2-1)\phi(x)\,dx=L(u)$. Consequently the coefficient of $\rho^2$ in $q-F_{\rho}(q)$ equals
\[
 D_m(q)=\tau\,\mathbb E\!\left[L(K\tau)
                 \{h(J\tau)-h(K\tau)\}\right] .
\]
All summands are nonnegative. On the event that every remaining $p$-value is at most $q$, which has probability $q^{m-2}$, we have $K=m$, $J\le m-1$. Hence
\[
 D_m(q)\ge q^{m-2}\tau L(q)\{h((m-1)q/m)-h(q)\}.
\]
Dividing this lower bound by $1-q$ gives a continuous positive function on $[a,1]$ when the value at 1 is assigned by continuity.

For a general zero-diagonal symmetric perturbation $E$, the event defining $F_\Sigma$ is invariant under independent coordinate sign changes. Thus $F_{I+E}=F_{I+SES}$ for every diagonal sign matrix $S$. Its Taylor polynomial in the off-diagonal entries therefore has zero linear coefficients and zero mixed quadratic coefficients: two distinct edges have a vertex occurring just once. Permutation symmetry identifies each pure quadratic coefficient with $-D_m(q)$, and $F_I(q)=q$ by the independent Simes identity.

For the uniform remainder, integrate the density on the acceptance event, a subset of $\bigcup_i\{|Z_i|<z_q\}$ up to a null set. For $\|E\|_{\rm F}<1/2$ the third derivative in $t$ of the $N(0,I_m+tE)$ density at $x$ is bounded on $0\le t\le1$ by
\[
 C_m\|E\|_{\rm F}^3(1+\|x\|^6)\exp(-c_m\|x\|^2).
\]
Indeed $\|(I+tE)^{-1}\|_{\rm op}\le2$ and each differentiation of the log density contributes determinant traces and quadratic forms, with total degree in $E$ equal to the number of differentiations. The integral of the displayed envelope over any one strip is $O_{m,a}(z_q\|E\|_{\rm F}^3)$, and $z_q=O_a(1-q)$ on $a\le q<1$. This proves the expansion by Taylor's formula. Since $\sum_{i<j}E_{ij}^2=\|E\|_{\rm F}^2/2$, it also proves strictness on $[a,1)$ for all sufficiently small nonzero perturbations.

Finally, take $a$ inside the positive fixed-dimensional interval in \citet[Theorem~2.1]{Zhang}. That theorem controls all correlation matrices on $(0,a]$, with equality precisely for independent groups of signed duplicates. If $\|E\|_{\rm F}<1/2$ the only such matrix is $I_m$. Combining the two intervals proves the first assertion.
\end{proof}

\subsection{Proof of Proposition~\ref{S:prop:trivariate-upper-end}}
\begin{proof}
Acceptance requires some $P_i>q$. On that central strip all other acceptance conditions are $P^{(-i)}_{(k)}>kq/m$ for $1\le k<m$. Two strips intersect with probability $O_\epsilon((1-q)^2)$, since the bivariate Gaussian densities are uniformly bounded and $z_q=O(1-q)$. Conditional on $|Z_i|<z_q$, $Z_i$ tends to zero and the conditional law of $Z_{-i}$ given $Z_i=x$ is Gaussian with mean $\Sigma_{-i,i}x$ and covariance $\Sigma_{-i,-i}-\Sigma_{-i,i}\Sigma_{i,-i}$. Its marginal variances are uniformly positive. If $\Sigma$ converges in the indicated compact class while $q\uparrow1$, Gaussian square-root continuity and the zero probability of every threshold boundary give convergence of the conditional acceptance probability to $A_i$ at the limiting matrix. The same argument makes $A_i$ continuous there. Inclusion--exclusion over strips gives the uniform expansion.

For $m=3$, given $Z_i=0$ write the two remaining coordinates as $\sigma_j Y_j$, $\sigma_k Y_k$, where $Y_j,Y_k$ are standardized Gaussians and $\sigma_l^2=1-\Sigma_{il}^2\le1$. Their original $p$-values dominate $U_l=2\bar\Phi(|Y_l|)$. The folded bivariate Simes inequality at level $2/3$ gives
$\Prob\{\min(U_j,U_k)>1/3,\ \max(U_j,U_k)>2/3\}\ge1/3$. For completeness, for bivariate correlation $|\rho|<1$ the conditional probability $\Prob(U_k\le2/3\mid U_j=u)$ is nonincreasing in $u$: the probability of a Gaussian shift outside a centered interval increases with the absolute shift. Then $2C_\rho(1/3,2/3)\ge C_\rho(2/3,2/3)$, which is exactly the stated Simes inequality. At $|\rho|=1$ the two uniforms coincide, with equality.

If $\sigma_j<1$ and the standardized pair is nonsingular, take $U_j$ in a sufficiently short interval just below $1/3$ and $U_k>2/3$. This open set has positive probability. Shrinkage strictly increases the $j$th $p$-value beyond $1/3$, making the full acceptance event true when the standardized event fails. If the standardized pair has rank one then $U_j=U_k$; take their common value just below $2/3$. Both are above $1/3$, and the shrunken $j$th coordinate now gives a $p$-value above $2/3$. Thus $A_i>1/3$ whenever $i$ is incident to a nonzero correlation. For any $\Sigma\ne I_3$ at least one such $i$ exists. The uniform expansion and compactness give strict control near level one off a neighborhood of $I_3$. Theorem~\ref{S:thm:local-all-level} handles that neighborhood.
\end{proof}

\subsection{Proof of Corollary~\ref{S:cor:rank-two-upper-end}}
\begin{proof}
Conditionally on $Z_i=0$, the remaining absolute statistics have the
form $\sigma_{ij}|Y|$ for a single standard normal $Y$, where
$\sigma_{ij}=\sqrt{1-\Sigma_{ij}^2}>0$. Write
$\sigma_{i,[1]}\ge\cdots\ge\sigma_{i,[m-1]}$ for these numbers in
decreasing order. Then
\[
 A_i(\Sigma)=1-2\bar\Phi\!\left(
       \min_{1\le k<m}\frac{z_{k/m}}{\sigma_{i,[k]}}\right)
       \ge\frac1m.
\]
For $k<m-1$, $z_{k/m}>z_{(m-1)/m}$. At $k=m-1$ the inequality
$z_{k/m}/\sigma_{i,[k]}>z_{(m-1)/m}$ holds whenever $i$ is incident
to a nonzero correlation. Thus $\sum_i A_i(\Sigma)>1$ for
$\Sigma\ne I_m$. When $m>2$ the compact rank-constrained class excludes
$I_m$, so Proposition~\ref{S:prop:trivariate-upper-end} gives a single
endpoint interval for that class. For $m=2$ the bivariate Simes
inequality and Theorem~\ref{S:thm:local-all-level} give the conclusion,
using compactness away from $I_2$ as above.
\end{proof}

\subsection{Proof of Theorem~\ref{S:thm:prds-interior}}
\begin{proof}
Theorem~\ref{S:thm:local-all-level} gives the assertion for $\Sigma_0=I_m$.
Fix a dependent $\Sigma_0$ and work under its distribution. In the
leave-one-out identity write $t_k=kq/m$, let $R_i$ be the rejection count
after forcing $P_i=0$, and set $D_{ik}=\{R_i\ge k\}$. Direct telescoping
gives
\[
 q-F_{\Sigma_0}(q)=\frac qm\sum_i\sum_{k=2}^m
 \big\{\Pr(D_{ik}\mid P_i\le t_{k-1})
           -\Pr(D_{ik}\mid P_i\le t_k)\big\}.
\]
Each brace is nonnegative by PRDS since $D_{ik}$ is a decreasing
event in the remaining $p$-values. At $k=m$, $D_{im}$ is just
$D_i(q)=\{P_j\le q \ \hbox{for all }j\ne i\}$.

Put $z_s=\bar\Phi^{-1}(s/2)$ and
$g_i(x)=\Pr(D_i(q)\mid Z_i=x)$. Conditional on $Z_i=x$, the
other statistics have nonsingular Gaussian covariance independent
of $x$ and mean $(\Sigma_{0,ji}x)_{j\ne i}$. Thus $g_i$ is even
and real analytic in $x$. For example, writing the shifted Gaussian
density as the density at zero mean times an exponential linear in the
integration variable and quadratic in $x$ justifies its analytic
extension to complex $x$ on every bounded disk, by Gaussian domination.
If $i$ has a nonzero covariance with another index, $g_i$ is
nonconstant. Indeed, as $x\to\infty$, every coordinate with nonzero
shift coefficient exceeds $z_q$ in absolute value with probability
tending to one, and the joint probability for all unshifted coordinates
is constant. At $x=0$ this limiting probability is strictly larger
than $g_i(0)$, since the conditional covariance is nonsingular and
the open event with all unshifted coordinates outside $[-z_q,z_q]$
and a chosen shifted coordinate inside $(-z_q,z_q)$ has positive
probability. This includes the case of no unshifted coordinates.
The function
\[
 H_i(s)=\Pr(D_i(q)\mid P_i\le s)
       =\frac1s\int_{z_s}^\infty 2\phi(x)g_i(x)\,dx
\]
is therefore real analytic and nonconstant on $(0,1)$: if it were
constant, differentiation would force $g_i$ to be constant as well.
By PRDS, $H_i$ is nonincreasing, so it is strictly decreasing on
$(0,1)$. Hence the $k=m$ brace is strictly positive for any incident
$i$, giving $F_{\Sigma_0}(q)<q$ at every $0<q<1$.

We also need a strict margin at the upper endpoint. Set
$s_0=(m-1)/m$ and $\sigma_{ij}=(1-\Sigma_{0,ij}^2)^{1/2}$.
For any sequence $s(q)\to s_*\in(0,1]$, conditioning on $Z_j$
in the strip $|Z_j|<z_q$ gives
\[
 \frac{\Pr(P_j>q\mid P_i\le s(q))}{1-q}
       \ \longrightarrow\
       \frac{2\bar\Phi(z_{s_*}/\sigma_{ij})}{s_*}.
\]
Pair intersections among the events $\{P_j>q\}$, conditionally on
$P_i\le s(q)$, have probability $O((1-q)^2)$ by the bounded bivariate
densities and $s_*>0$. Consequently the last-rank terms yield
\[
 \liminf_{q\uparrow1}\frac{q-F_{\Sigma_0}(q)}{1-q}
 \ge \frac1m\sum_{i}\sum_{j\ne i}
       \left\{1-\frac{2\bar\Phi(z_{s_0}/\sigma_{ij})}{s_0}\right\}>0 .
\]
Indeed each summand is nonnegative, and strictness holds whenever
$\Sigma_{0,ij}\ne0$. The uniform expansion and continuity in
Proposition~\ref{S:prop:trivariate-upper-end} transfer this positive
upper-endpoint margin to a single interval and one neighborhood of
$\Sigma_0$. The low-level theorem of \citet[Theorem~2.1]{Zhang}
gives a fixed lower-endpoint interval with strict inequality for
positive definite dependent matrices. On the remaining compact
interval the strict gap just proved has a positive minimum, by joint
continuity of $F_\Sigma(q)$. Intersecting these neighborhoods yields
the assertion.
\end{proof}

\subsection{Proof of Corollary~\ref{S:cor:non-prds-interior}}
\begin{proof}
The three-coordinate family below extends the obstruction at
$a=b=1/4$ in \citet[Lemma~F.11]{Zhang}. Choose $0<b<1$, take
independent standard normals $X,Y,E$, and set
$Z_1=X$, $Z_3=Y$, and $Z_2=aX+bY+sE$ with
$s=(1-a^2-b^2)^{1/2}>0$. Append independent standard normals if
$m>3$. At $a=0$ the two possibly dependent coordinates form a folded
bivariate Gaussian pair, which satisfies PRDS by the monotonicity
of the two-sided Gaussian-shift tail. Independent coordinates preserve
PRDS. By Theorem~\ref{S:thm:prds-interior}, all sufficiently small
$a>0$ and neighboring covariances satisfy strict Simes control at
every level.

For any such $a>0$, consider the decreasing $p$-value event
$D_T=\{|Z_2|>1,\ |Z_3|>T\}$. Then
\[
 \frac{\Pr(D_T)}{\Pr(|Y|>T)}\longrightarrow1,\qquad
 \frac{\Pr(D_T\mid |X|>bT/a)}{\Pr(|Y|>T)}
       \longrightarrow\frac{1+2\bar\Phi(1/s)}2<1.
\]
For the first limit, condition on $|Y|>T$ and use independence of
$aX+sE$. For the second, also condition on $|X|>bT/a$. The Gaussian
overshoot above a threshold tending to infinity converges in
probability to zero, from
$\bar\Phi(t+h)/\bar\Phi(t)\le e^{-th}$ for $h>0$. Thus
$a|X|-b|Y|\to0$ conditionally in probability. The two signs agree
or disagree with equal probability. For agreeing signs
$|aX+bY|\to\infty$, and for opposite signs it tends to zero. The
residual $sE$ remains independent. These are the two limits displayed
above.
At a sufficiently large finite $T$, therefore,
$\Pr(D_T\mid P_1\le2\bar\Phi(bT/a))<\Pr(D_T)$, in contradiction
to PRDS. This strict inequality persists on an open neighborhood
by continuity of Gaussian probabilities. Intersect this neighborhood
with the one from Theorem~\ref{S:thm:prds-interior} and choose a ball
around a matrix with all off-diagonal entries nonzero. The construction
applies to every $b<1$.
\end{proof}

\subsection{Proof of Theorem~\ref{S:thm:spherical-design}}
\begin{proof}
Write $Q_d(r)=\Pr(\|G\|\ge r)$. For $r>0$ and $v\in S^{d-1}$ set
\[
 A_m(r,v)=\sup_{0\le x\le1}
   \frac{m^{-1}\sum_i\mathbf 1\{|u_{m,i}^{\top}v|\ge x\}}{S(rx)},
 \quad
 A(r)=\sup_{0\le x\le1}\frac{D_d(x)}{S(rx)} .
\]
The first supremum is the reciprocal of the Simes $p$-value at $G=rv$:
between successive observed absolute projections the numerator is
constant and $1/S(rx)$ increases. The crossing event is increasing
in $r$ on each ray. Uniformly in $v$,
\[
 |Q_d(r)A_m(r,v)-Q_d(r)A(r)|\le
            \Delta_m Q_d(r)/S(r)\le C_d\Delta_m(1+r)^{d-1}.
 \tag{S}
\]
The last inequality follows from the normal and $\chi_d$ tail
asymptotics and continuity on bounded intervals. For $r>0$,
$Q_d(r)A(r)<1$. In fact $S(rx)=\E[D_d(rx/\|G\|)]$, with $D_d=0$
above one. For $0<x<1$ the integrand strictly exceeds $D_d(x)$
on $\|G\|>r$, and the cases $x=0,1$ are immediate. The supremum is
attained by continuity on $[0,1]$. The map $r\mapsto Q_d(r)A(r)$
is continuous.

We next prove $Q_d(r)A(r)\to\kappa_d$ as $r\to\infty$. Write
$w=r^2(1-x^2)$. Since $1-U_1^2$ has the $\operatorname{Beta}(h,1/2)$
distribution,
\[
 D_d(x)\sim\frac{(1-x^2)^h}{hB(h,1/2)}\quad(x\uparrow1),\qquad
 D_d(x)\le C_d(1-x^2)^h\quad(0\le x\le1).
\]
Using $Q_d(r)\sim 2^{1-d/2}r^{d-2}e^{-r^2/2}/\Gamma(d/2)$ and
$S(z)\sim\sqrt{2/\pi}\,e^{-z^2/2}/z$, we obtain, uniformly on each
bounded $w$-interval,
\[
 \frac{Q_d(r)D_d(x)}{S(rx)}
       =\frac{w^h e^{-w/2}}{2^h\Gamma(h+1)}+o(1).
\]
Uniformity at zero follows by expressing $D_d(x)/(1-x^2)^h$ by its
continuous extension. For all $w$ and $r\ge1$ the left side is at most
$C_dw^h e^{-w/2}$ by the above cap bound and
$S(z)\ge c e^{-z^2/2}/(1+z)$. The limiting maximum occurs at
$w=2h$. Moreover $\kappa_d<1$, since
$\Gamma(h+1)>\int_h^\infty h^h e^{-t}dt=(h/e)^h$.

Near radius zero, $A(r)=1$ throughout some $0\le r\le r_0$.
Indeed, $\inf_{0<x\le1}(1-D_d(x))/x>0$, since $|U_1|$ has positive
continuous density at zero, whereas $1-S(rx)\le\sqrt{2/\pi}\,rx$.
For $r\ge r_0$ the preceding strict inequality, continuity and the
limit at infinity give $\sup_{r\ge r_0}Q_d(r)A(r)<1$.
Take $r=Q_d^{-1}(q)$ and $q\ge q_m$; then
$r^2=O_d(1+\log(1/q_m))$. If $r\ge r_0$, (S) gives
$qA_m(r,v)<1$ uniformly for large $m$. If $r<r_0$, then
$A_m(r,v)\le1+\Delta_m/S(r_0)$ and $qA_m(r,v)<1$ whenever
$q\le1-C_d\Delta_m$, for example with $C_d=2/S(r_0)$.
No direction rejects at this radius; continuity of the finitely many
ordered projections on the compact sphere gives the same on a slightly
larger sphere. This proves strict validity.

If $q_m\to0$, fix $t>0$ and put $r_t=Q_d^{-1}(tq_m)$ for $m$ large.
Uniformly in $v$, (S) gives $q_m A_m(r_t,v)\to\kappa_d/t$. For
$t>\kappa_d$ no direction rejects at $r_t$, so size is at most
$tq_m$. For $0<t<\kappa_d$ all directions reject there, so size is at
least $tq_m$. Let $t$ decrease or increase to $\kappa_d$.
\end{proof}

\section{Finite Gaussian examples}\label{S:sec:numerics}

\subsection{An exact comparison with one true null}\label{S:sec:one-null}

For one true null the dependence on the alternative distribution can be calculated explicitly. Throughout this subsection all statistics have unit marginal variance. Let $H_m=\sum_{j=1}^m1/j$.

\begin{proposition}\label{S:prop:one-null}
For one true null and $m-1$ nonzero Gaussian mean shifts, the supremum of the one-sided FDR over joint Gaussian laws is $(q/m)H_m$ for $0<q\le1/2$. For two-sided tests with one true null and one Gaussian mean shift, let $D(q)$ denote the corresponding supremum. Then for every $0<q<1$,
\[
 \frac{5q}{8}-\frac12\bar\Phi\bigl(t_{q/2}+2t_q\bigr)
 \ \le\ D(q)\ \le\ \frac{5q}{8},
 \qquad t_x=\bar\Phi^{-1}(x/2).
\]
In particular $D(q)\sim5q/8$ as $q\downarrow0$.
\end{proposition}

For two tests and one unrestricted alternative, the worst FDR is exactly $3q/4$. For one-sided Gaussian tests, \citet[Section~4.1.2]{RV} obtained this value for two perfectly negatively correlated coordinates, one with a nonzero mean. The construction below gives the corresponding harmonic value for arbitrary $m$. Proposition~\ref{S:prop:one-null} shows how this changes under a Gaussian mean shift. Conditional on the sign of the null statistic, the signal can favor rejection on one tail while the opposite tail restricts its ability to follow the BH thresholds. In the one-sided case a sequence of Gaussian mean shifts realizes the harmonic value exactly. Section~\ref{S:app:one-null} gives both constructions.

\subsection{Size and FDR at fixed nominal levels}\label{S:sec:examples}

Take $m=200$ and independent standard normals $Z,\varepsilon_1,\ldots,\varepsilon_{199}$. At every nominal level use the Gaussian statistics
\begin{equation}\label{S:eq:example-model}
 T_i=0.3Z+\sqrt{0.91}\,\varepsilon_i\quad(1\le i\le199),
 \qquad T_{200}=5.5-Z .
\end{equation}
The map from the 200 independent Gaussian inputs to $(T_1,\ldots,T_{200})$ is invertible. Its covariance is therefore positive definite. Each pair of null statistics has correlation $0.09$, while each null has correlation $-0.3$ with the signal. This parameter choice illustrates the common-factor mechanism studied by \citet{Dobriban} and \citet{Lei}.

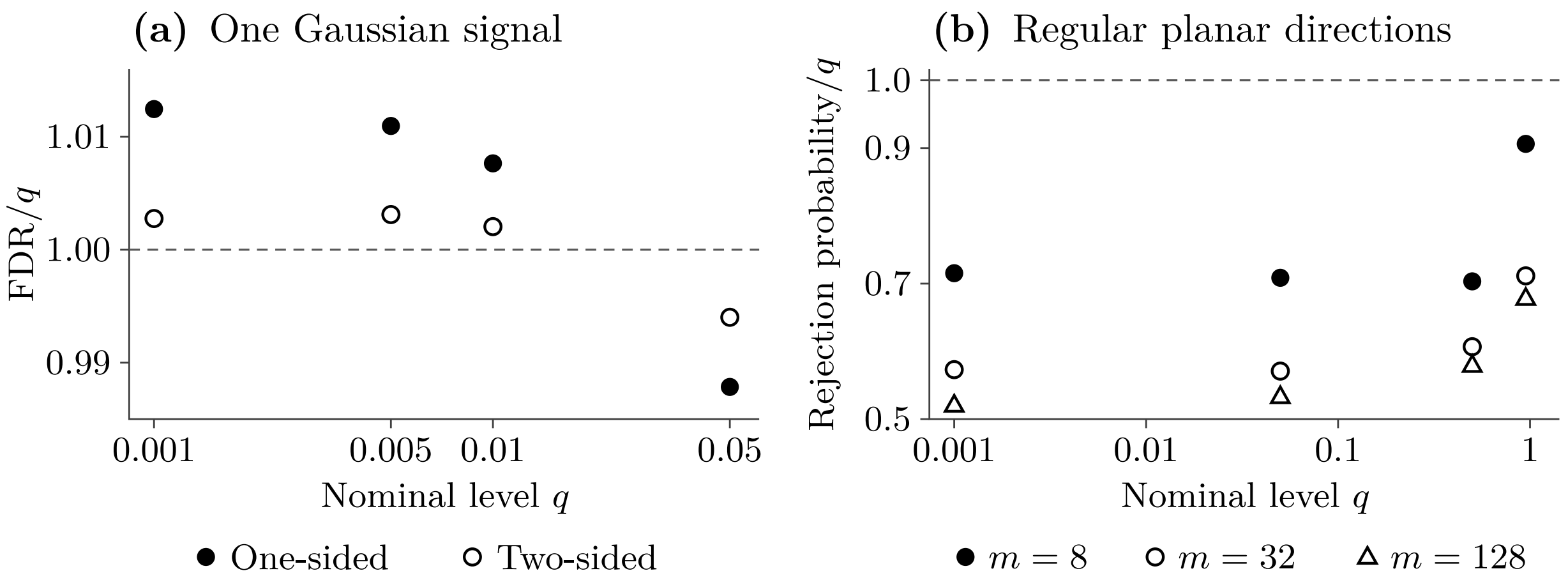
\begin{figure}[tbp]
\centering
\begingroup
\tikzset{suaxis/.style={draw=black!75,line width=.45pt},
  sureference/.style={draw=black!65,dash pattern=on 2.7pt off 1.9pt,line width=.5pt},
  sudot/.style={draw=black,fill=black,line width=.45pt},
  suring/.style={draw=black,fill=white,line width=.7pt}}
\newcommand{\suTriangle}[2]{\begin{scope}[shift={(#1,#2)}]
  \draw[suring] (0,2.4pt)--(-2.25pt,-1.6pt)--(2.25pt,-1.6pt)--cycle;
  \end{scope}}
\begin{minipage}[t]{.49\linewidth}\centering
\resizebox{\linewidth}{!}{\begin{tikzpicture}[x=1cm,y=1cm,font=\footnotesize,inner sep=1pt]
  \path[use as bounding box] (-1.16,-1.50) rectangle (5.65,3.67);
  \node[anchor=west,font=\small] at (-.02,3.43) {\textbf{(a)}\enspace One Gaussian signal};
  \draw[suaxis] (0,3.10)--(0,0)--(5.58,0);
  \draw[sureference] (0,1.5)--(5.58,1.5);
  \foreach \y/\labeltext in {.5/0.99,1.5/1.00,2.5/1.01}{
    \draw[suaxis] (0,\y)--(-.08,\y);
    \node[anchor=east] at (-.12,\y) {$\labeltext$};}
  \node[rotate=90] at (-.93,1.55) {$\mathrm{FDR}/q$};
  \foreach \t/\labeltext in {1/0.001,5/0.005,10/0.01,50/0.05}{
    \pgfmathsetmacro{\xx}{.22+5.1*ln(\t)/ln(50)}
    \draw[suaxis] (\xx,0)--(\xx,-.08);
    \node[anchor=north] at (\xx,-.13) {$\labeltext$};}
  \node at (2.8,-.70) {Nominal level $q$};
  \foreach \t/\ratio in {1/1.012445575,5/1.010955508,10/1.007642529,50/0.987848916}{
    \filldraw[sudot] ({.22+5.1*ln(\t)/ln(50)},{100*(\ratio-.985)}) circle[radius=2.05pt];}
  \foreach \t/\ratio in {1/1.002760907,5/1.003111584,10/1.002048066,50/0.994025681}{
    \draw[suring] ({.22+5.1*ln(\t)/ln(50)},{100*(\ratio-.985)}) circle[radius=2.05pt];}
  \filldraw[sudot] (.68,-1.22) circle[radius=2.05pt];
  \node[anchor=west] at (.85,-1.22) {One-sided};
  \draw[suring] (3.04,-1.22) circle[radius=2.05pt];
  \node[anchor=west] at (3.21,-1.22) {Two-sided};
\end{tikzpicture}}
\end{minipage}\hfill
\begin{minipage}[t]{.49\linewidth}\centering
\resizebox{\linewidth}{!}{\begin{tikzpicture}[x=1cm,y=1cm,font=\footnotesize,inner sep=1pt]
  \path[use as bounding box] (-1.16,-1.50) rectangle (5.65,3.67);
  \node[anchor=west,font=\small] at (-.02,3.43) {\textbf{(b)}\enspace Regular planar directions};
  \draw[suaxis] (0,3.1)--(0,0)--(5.58,0);
  \draw[sureference] (0,3.0)--(5.58,3.0);
  \foreach \y/\labeltext in {0/0.5,1.2/0.7,2.4/0.9,3/1.0}{
    \draw[suaxis] (0,\y)--(-.08,\y);
    \node[anchor=east] at (-.12,\y) {$\labeltext$};}
  \node[rotate=90] at (-.93,1.55) {Rejection probability/$q$};
  \foreach \t/\labeltext in {1/0.001,10/0.01,100/0.1,1000/1}{
    \pgfmathsetmacro{\xx}{.22+5.1*ln(\t)/ln(1000)}
    \draw[suaxis] (\xx,0)--(\xx,-.08);
    \node[anchor=north] at (\xx,-.13) {$\labeltext$};}
  \node at (2.8,-.70) {Nominal level $q$};
  \foreach \t/\ratio in {1/.7153257,50/.7085657,500/.7032946,950/.9061615}{
    \filldraw[sudot] ({.22+5.1*ln(\t)/ln(1000)},{6*(\ratio-.5)}) circle[radius=2.05pt];}
  \foreach \t/\ratio in {1/.5730574,50/.5708959,500/.6068712,950/.7114331}{
    \draw[suring] ({.22+5.1*ln(\t)/ln(1000)},{6*(\ratio-.5)}) circle[radius=2.05pt];}
  \foreach \t/\ratio in {1/.5209770,50/.5335245,500/.5797665,950/.6786515}{
    \pgfmathsetmacro{\xx}{.22+5.1*ln(\t)/ln(1000)}
    \pgfmathsetmacro{\yy}{6*(\ratio-.5)}
    \suTriangle{\xx}{\yy}}
  \filldraw[sudot] (.32,-1.22) circle[radius=2.05pt];
  \node[anchor=west] at (.48,-1.22) {$m=8$};
  \draw[suring] (2.00,-1.22) circle[radius=2.05pt];
  \node[anchor=west] at (2.16,-1.22) {$m=32$};
  \suTriangle{3.88}{-1.22}
  \node[anchor=west] at (4.04,-1.22) {$m=128$};
\end{tikzpicture}}
\end{minipage}
\endgroup
\caption{False discovery rate and Simes rejection probability relative to the nominal level. Panel (a) gives the FDR for 199 true nulls and one signal in \eqref{S:eq:example-model}, divided by the nominal level $q$. Panel (b) gives the corresponding rejection probability under the global null for $m$ regularly spaced planar directions and two-sided tests. Each point is a numerical integral. The horizontal axes are logarithmic, the vertical scales differ, and the dashed line in each panel marks equality with $q$.}
\label{S:fig:calibration}
\end{figure}

Compare this mixed configuration with the global-null planar models of Section~\ref{S:sec:planar-main}. Figure~\ref{S:fig:calibration} plots each error rate divided by its nominal level. Both calculations use one-dimensional integration. For the model in \eqref{S:eq:example-model}, conditioning on $Z$ makes the null $p$-values independent and gives a recursion for the distribution of the BH rejection count. For the planar models, conditioning on the angle gives the integral in Section~\ref{S:sec:planar-main}.

At $q=0.005$, the one- and two-sided FDRs in \eqref{S:eq:example-model} are $0.00505478$ and $0.00501556$. Both exceed $q$. At $q=0.05$, both fall below $q$ in the same model. The magnitude of the excess also illustrates the distinction between an observed FDR and a worst-case upper bound. With a single signal, Theorem~\ref{S:thm:finite} bounds the two-sided FDR divided by $q$ by $(199/200)e(1+\log 2)\approx4.58$. For this specific covariance, the null Simes inequality holds with constant one \citep[Remark~4.9]{Lei}, yielding the smaller bound $1.685$.

For the global-null planar models, every computed size is below the nominal level. With $m=128$ directions, the size divided by level ranges from $0.521$ at $q=0.001$ to $0.679$ at $q=0.95$. The all-level result for sufficiently large $m$ follows from Theorem~\ref{S:thm:polygon}. Section~\ref{S:app:num} reports the numerical values and the integration checks.

\subsection{Details of the finite examples}\label{S:app:num}

\subsection{Proof of Proposition~\ref{S:prop:one-null}}\label{S:app:one-null}

\begin{proof}
For one-sided tests the single-null upper bound follows from the interval argument in Section~\ref{S:app:finite-bounds}. For a composite one-sided null the stepwise FDP envelope is nonincreasing in its $p$-value, so the boundary-null law attains the largest expectation. The value is known for monotone functions of Gaussian variables \citep[Theorem~5.2]{LyuWang}. To attain it with Gaussian mean shifts, put $h=q/m$, $z_j=\bar\Phi^{-1}(jh)$ and take $Z\sim N(0,1)$. Use $Z$ for the true-null statistic and $X_j=z_j+z_{j+1}-Z$ for the alternatives, $1\le j<m$. Their means are positive. If $U=\bar\Phi(Z)\in((r-1)h,rh)$, then the alternatives with $j<r$ all have $p$-values at most $(j+1)h$, so rank $r$ is feasible. At a larger rank $k$, every alternative with $j\ge k-1$ has $p$-value strictly above $kh$, since $Z>z_{k-1}\ge z_j$. At most $k-1$ hypotheses then pass the rank-$k$ threshold. Hence BH rejects exactly $r$, including the true null. For $U>q$ it cannot reject the true null. Integrating over the $m$ intervals gives $hH_m$.

For the two-sided result let $h=q/2$, $U=2\bar\Phi(|Z|)$, and let $E$ be the event that the signal $p$-value is at most $q$. Directly from the two BH critical values,
\[
 \FDR=h+\tfrac12\{\Prob(h<U\le q,E)-\Prob(U\le h,E)\}.
\]
Represent the signal statistic as $\mu+\rho Z+\sqrt{1-\rho^2}\,\varepsilon$ with $\varepsilon$ standard normal independent of $Z$. On the sign of $Z$ for which $\rho Z$ is aligned with $\mu$, the conditional absolute mean increases with $|Z|$. A shifted normal's probability of falling outside $[-t_q,t_q]$ increases with the magnitude of the shift, by differentiation (also for zero residual variance by direct inspection). Thus the conditional probability of $E$ decreases with $U$ on this sign. The two adjacent $U$-intervals each have conditional mass $h$ because $U$ is uniform and independent of the sign. Their probability difference is nonpositive on the aligned sign and at most $h$ on the other. The two signs are equiprobable, giving $\FDR\le h+h/4$.

For the lower bound use signal statistic $t_h+t_q-Z$. On the negative sign $E$ occurs throughout both intervals. On the positive sign it occurs throughout $h<U\le q$ and, for $U\le h$, exactly when $Z\ge t_h+2t_q$, up to endpoints. Substitution gives the stated lower bound. Since $\bar\Phi(t_h+2t_q)=o(q)$, the first-order coefficient follows.
\end{proof}

\subsection{Rank probabilities and numerical integration}\label{S:app:quadrature}

Tables~\ref{S:tab:finite-fdr} and~\ref{S:tab:polygon} give the values used in Figure~\ref{S:fig:calibration}. The first calculation integrates the distribution of the final BH rank; the second integrates the exact radial formula for the planar model.

\begin{table}[!htbp]
\centering
\caption{BH applied to the fixed model \eqref{S:eq:example-model}. Each row changes only the nominal level. The independent-coordinate reference is $(199/200)q$.}
\label{S:tab:finite-fdr}
\begin{tabular}{r rr rr}
\toprule
& \multicolumn{2}{c}{One-sided} & \multicolumn{2}{c}{Two-sided}\\
$q$ & FDR & FDR/$q$ & FDR & FDR/$q$\\
\midrule
0.050 & 0.049392446 & 0.987849 & 0.049701284 & 0.994026\\
0.010 & 0.010076425 & 1.007643 & 0.010020481 & 1.002048\\
0.005 & 0.005054778 & 1.010956 & 0.005015558 & 1.003112\\
0.001 & 0.001012446 & 1.012446 & 0.001002761 & 1.002761\\
\bottomrule
\end{tabular}
\end{table}

\begin{table}[!htbp]
\centering
\caption{Two-sided global-null Simes rejection probability divided by nominal level for regular planar directions.}\label{S:tab:polygon}
\begin{tabular}{rrrrr}\toprule
$m$ & $\alpha=0.001$ & $\alpha=0.05$ & $\alpha=0.50$ & $\alpha=0.95$\\ \midrule
8   &0.71533&0.70857&0.70329&0.90616\\
32  &0.57306&0.57090&0.60687&0.71143\\
128 &0.52098&0.53352&0.57977&0.67865\\ \bottomrule
\end{tabular}
\end{table}

\paragraph{Exact finite-sample integral.}
Write $n=199$, $t_j=jq/m$, and $b=1$ or $2$ for the number of sides. Given $Z=z$, the null $p$-values are independent with common CDF
\[
 F_j=\bar\Phi\!\left(\frac{u_j-0.3z}{\sqrt{0.91}}\right)
 +\mathbf1_{\{b=2\}}\bar\Phi\!\left(\frac{u_j+0.3z}{\sqrt{0.91}}\right),
 \qquad u_j=\bar\Phi^{-1}(t_j/b),
\]
at $t_j$. Let $A_j$ indicate whether the signal $p$-value is at most $t_j$, let $r_j=j-A_j$, and write $d(h,k,p)=\Pr\{\operatorname{Bin}(h,p)=k\}$, taken as zero for infeasible indices. For $p_j=\Pr(R=j\mid Z=z)$, the downward recursion from $m$ to $1$ gives
\begin{equation}\label{S:eq:example-ranks}
 p_j=d(n,r_j,F_j)-\sum_{k=j+1}^m p_k\,d(r_k,r_j,F_j/F_k),
 \qquad
 \FDR=\int_{\mathbb R}\phi(z)\sum_{j=1}^m\frac{r_j}{j}p_j\,dz .
\end{equation}
If exactly $r_j$ nulls lie below $t_j$, rank $j$ is feasible. Partition this exact-count event by the final rank $k\ge j$. Given $R=k$, there are exactly $r_k$ nulls below $t_k$. Their locations have the null law truncated at $t_k$, since all conditions ruling out larger ranks involve only counts at or above $t_k$. The number of these nulls below $t_j$ is thus binomial with success probability $F_j/F_k$. This proves the recursion, including ranks at which the signal changes status.

\paragraph{Numerical accuracy.}
The accompanying base-R script splits $[-9,9]$ at all crossings of the signal over the BH thresholds and uses eight-point Gauss--Legendre quadrature on subintervals of length at most $0.2$. The omitted Gaussian mass is $2\bar\Phi(9)<2.3\times10^{-19}$. The script may omit ranks that would require the conditional binomial null count at $q$ to exceed a cutoff having tail probability at most $10^{-20}$.

To bound the quadrature error, let $g(z)=\mathbb E[\mathrm{FDP}\mid Z=z]$. Within each subinterval, the comparisons involving the signal are fixed. The remaining expectation is a fixed bounded function integrated against $n$ independent normals whose common mean is $0.3z$. Differentiating this translated density and using the $L^2$ norm of the normal Hermite polynomial gives
$|g^{(k)}(z)|\le \sqrt{k!}\,L^k$ with $L=0.3\sqrt{199/0.91}$.
Fourier inversion gives $\|\phi^{(k)}\|_\infty\le \mathbb E|Z|^k/\sqrt{2\pi}\le\sqrt{k!}/\sqrt{2\pi}$. The product rule therefore bounds the $k$th derivative of $\phi g$ by $\sqrt{k!}(1+L)^k/\sqrt{2\pi}$. The Gauss--Legendre remainder over all panels is at most
\[
 18\,\frac{(8!)^4}{17(16!)^3}\frac{\sqrt{16!}}{\sqrt{2\pi}}
 [0.2(1+L)]^{16}<2.2\times10^{-15}.
\]
This estimate concerns exact-arithmetic quadrature. As a separate floating-point check, the script replaces one true-null $p$-value by zero and recomputes the rank probabilities with $n-1$ random nulls and $r_j=j-1-A_j$. Conditional independence gives the alternative FDR integrand $n\sum_j F_j p_j/j$, using these new rank probabilities. Ten-point quadrature on a different mesh agrees with the direct calculation within $3\times10^{-17}$ in every table entry. Both computations recover the independent-coordinate FDR when the null loading is set to zero.

\section{Further properties of the maximal inequality}\label{S:app:more}

\subsection{The loss in smoothing empirical counts}\label{S:sec:smooth-sharpness}

Fix $m\ge1$ and $0<p_{(1)}\le\cdots\le p_{(m)}\le1$. Write
\[
 T=\max_{1\le k\le m}\frac{k}{m p_{(k)}},\qquad
 g_k=\left(\prod_{j=1}^k p_{(j)}\right)^{1/k}.
\]
Let
\[
 L=\max\left\{1,\sup_{0<s\le1}
     \frac1{ms}\sum_{i=1}^m\left(\log\frac{s}{p_{(i)}}\right)_+
     \right\}.
\]
Thus $T$ is the reciprocal of the Simes statistic. Including the constant
one in $L$ covers samples for which every smoothed count is below one.

\begin{lemma}\label{S:lem:smooth-sharpness}
For every sample,
\[
 L=\max\left\{1,\frac1e\max_{1\le k\le m}\frac{k}{m g_k}\right\},
 \qquad \frac{T}{e}\le L\le T.
\]
The factor $e$ is attained for every $m$. Moreover, $L$ equals the
supremum of
\[
 \frac{m^{-1}\sum_{i=1}^m f(-\log p_{(i)})}
      {\mathbb E f(E)},\qquad E\sim\mathrm{Exp}(1),
\]
over all finite nonnegative, nondecreasing convex functions
$f$ on $[0,\infty)$ for which $0<\mathbb E f(E)<\infty$.
\end{lemma}

\begin{proof}
Denote the smoothed count at $s$ by $M(s)$ and extend its formula to
$s>1$. On this range $M(s)=\{\log s+M(1)\}/s
\le\max\{1,M(1)\}$ because $\log s\le s-1$. Thus extending the
supremum to all $s>0$ preserves $L$. For every $s>0$,
\[
 ms M(s)=\max\left\{0,\max_{1\le k\le m} k\log(s/g_k)\right\}.
\]
Indeed the sum of the positive parts is the maximum of the prefix
sums of the decreasing sequence $\log(s/p_{(i)})$, including the empty
prefix. For each $g>0$, the maximum of $s^{-1}\log(s/g)$ over $s>0$
is $1/(eg)$, attained at $s=eg$. Interchanging the two suprema proves
the formula for $L$. Since $g_k\le p_{(k)}$, it gives $L\ge T/e$.
The empirical distribution function $F(t)=m^{-1}\#\{i:p_{(i)}\le t\}$
satisfies $F(t)\le Tt$, hence
$M(s)=s^{-1}\int_0^s F(t)\,dt/t\le T$ for $s\le1$. Also $T\ge1$.

To prove the assertion about convex scores, let $d=f'_+(0)\ge0$ and let
$\nu$ be the measure of increments of the nondecreasing right derivative
of $f$ on $(0,\infty)$. Integration of the derivative gives
\[
 f(y)=f(0)+dy+\int_{(0,\infty)}(y-b)_+\,\nu(db).
\]
Since $\mathbb E(E-b)_+=e^{-b}$ for $b\ge0$, the normalized empirical
score is an average of $1$, $M(1)$, and $M(e^{-b})$ with nonnegative
weights proportional, respectively, to $f(0)$, $d$, and
$e^{-b}\nu(db)$. Tonelli's theorem justifies the integration. This
bounds the score by $L$. The constant and each hinge $y\mapsto(y-b)_+$
belong to the stated class, proving the reverse inequality.

The formula also characterizes equality in $T\le eL$. It holds exactly
when $T\ge e$ and $g_k\ge k/(mT)$ for every $k$. In particular, if
$p_{(i)}=u\le e^{-1}$ for every $i$, then $T=1/u$ and $L=1/(eu)$.
\end{proof}

The same sample geometry arises under a Gaussian null by taking all
coordinates equal. Their common two-sided $p$-value $U$ is uniform,
so $T=1/U$ and $L=\max\{1,1/(eU)\}$. In particular, for $h>1$,
\[
 \Pr(L\ge h)=\frac1{eh}.
\]
For $0<q<1/e$, the Simes event agrees exactly with
$\{L\ge 1/(eq)\}$ and has probability $q$. In this model the
pathwise comparison attains the factor $e$, while the maximal
probability bound $1/h$ exceeds the actual tail of $L$ by that same
factor. Thus this sharp pathwise constant does not establish sharpness
of the Gaussian Simes probability bound.

The bound $1/(eh)$ on the tail of $L$ already fails for two independent
uniform $p$-values $U,V$. For $h>1$ put $q=1/(eh)$. The displayed
formula for $L$ gives
\[
 \{L\ge h\}=\{\min(U,V)\le q/2\}\ \cup\ \{UV\le q^2\}.
\]
The first event has probability $q-q^2/4$. Outside it, the second event
has area $\int_{q/2}^{2q}(q^2/x-q/2)\,dx=q^2(\log 4-3/4)$, since
$2q<1$. Hence $\Pr(L\ge h)=q+(\log 4-1)q^2>q$, even though the
ordinary Simes test has rejection probability exactly $q$ in this model.

\subsection{Independent replication}\label{S:app:replication}

There is a simple reason to seek a linear global-null bound when the question asks only for convergence to zero. Let $A(t)$ be the supremum of the global-null Simes rejection probability over all dimensions and correlations for a fixed Gaussian testing convention. If $J$ independent copies of one $m$-coordinate model are concatenated, any block rejecting at $t$ forces a rejection in the combined vector at $Jt$ whenever $Jt<1$. Indeed, a witness of rank $k$ in the block is feasible at rank $k$ in the combined vector, since the threshold remains $kt/m$. Independence of the blocks and approximation of the supremum give
\[
 A(Jt)\ge1-\{1-A(t)\}^J.
\]
If $A(t_0)<1$ at any $0<t_0<1$, take $J=\lfloor t_0/t\rfloor$ for $t<t_0/2$. It follows that $A(t)\le 2t[-\log(1-A(t_0))]/t_0$. In particular, uniform convergence of the global-null rejection probability to zero already implies a linear upper bound at small levels. By Su's linking argument this global-null condition is sufficient for mixed-null FDR consistency as well.

\section{Weights and Gaussian subspace tests}\label{S:app:ext}

\subsection{Prior-weighted hypotheses}\label{S:app:weights}

Prior information can enter the BH thresholds through deterministic weights $w_i\ge0$ with $\sum_{i=1}^m w_i=m$. Following \citet{GRW}, order $P_i/w_i$ and use the BH critical values $kq/m$, taking $P_i/w_i=\infty$ when $w_i=0$. Count each discovery once in the FDP. This is the prior-weight version of a self-consistent testing rule \citep{BR}.

\begin{theorem}[Prior-weighted BH]\label{S:thm:weights}
Let $I_0$ denote the true nulls. For $i\in I_0$ set $c_i=1$ for a two-sided test and $c_i=2$ for a prespecified one-sided test. Then the FDP and FDR bounds \eqref{eq:fdx}--\eqref{eq:fdr} hold for the weighted procedure with
\[
 u=\frac{e q}{m}\sum_{i\in I_0}c_iw_i .
\]
\end{theorem}

\begin{proof}
Write $Q_i=2\bar\Phi(|Z_i|)$ for the folded centered null values, so $Q_i\le c_iP_i$. Discard true nulls of weight zero from the sums below, and put $d_i=c_iw_i$ and $D=\sum_{i\in I_0}d_i$. When $D=0$ no true null can be rejected. For $D>0$, use $Y_i=-\log Q_i$ and define, for $t>0$,
\[
 b_i(t)=-\log(e t d_i),\qquad
 J_t=\sum_{i\in I_0}\frac{d_i}{D}
             \frac{(Y_i-b_i(t))_+}{A(b_i(t))},\qquad
 A(b)=\E(E-b)_+,\quad E\sim{\rm Exp}(1).
\]
For $b\ge0$, $A(b)=e^{-b}$; for $b<0$, $A(b)=1-b\le e^{-b}$. This second case permits weights large enough that $etd_i>1$. Every $J_t$ has mean one and centrally symmetric convex sublevel sets. As $t$ decreases, each $b_i(t)$ increases. For any $b'\ge b$, tilting the exponential density by $(y-b')_+/A(b')$ gives a stochastically larger law than tilting by $(y-b)_+/A(b)$. Indeed, on the support of the earlier tilt the ratio $(y-b')_+/(y-b)_+$ is nondecreasing; normalizing the two densities gives the stochastic order by a single crossing. Consequently the conditional section probability $\ell_i(y)$ from Theorem~\ref{lem:max}, which decreases with $y$, has smaller expectation under the later tilt. Summing with coefficients $d_i/D$ and applying the same first-crossing sum as in Theorem~\ref{lem:max} gives $\Prob(\max_{t\in\mathcal T}J_t\ge h)\le1/h$ on every finite grid $\mathcal T$.

When BH makes $R>0$ rejections, set $t=qR/m$. Each of the $V$ rejected true nulls satisfies $Q_i\le d_i t$, hence $Y_i-b_i(t)\ge1$ and its contribution to $J_t$ is at least $1/(e t D)$. Thus $J_t\ge m\,\FDP/(e qD)$. Take $\mathcal T=\{qk/m:1\le k\le m\}$ in the preceding inequality and integrate as in Theorem~\ref{thm:main}.
\end{proof}

With all weights equal to one, this gives the more precise mixed-sidedness factor $u=eq(n_2+2n_1)/m$, where $n_2$ and $n_1$ count the two- and one-sided true nulls. The threshold weights in the theorem are specified before inspecting the test statistics.

\subsection{Overlapping Gaussian chi-square tests}\label{S:sec:groups}

Let $G\sim N_d(0,I_d)$. For each tested subspace $E_i\subseteq\R^d$,
let $\Pi_i$ be the orthogonal projection onto $E_i$, let
$r_i=\dim E_i\ge1$, and write $S_r(v)=\Prob(\chi_r\ge v)$ for the
survival function of the norm of an $r$-dimensional standard Gaussian.
The upper-tail $p$-value for $\|\Pi_iG\|^2$ is
$P_i=S_{r_i}(\|\Pi_iG\|)$. Ranks may differ and subspaces may overlap,
be nested, or coincide. Rank one gives a two-sided Gaussian test.

\begin{theorem}[Gaussian group tests]\label{S:thm:groups}
For any finite collection of these subspaces,
\[
 \Prob\{\Simes(P)\le t\}\le \min(1,et),\qquad 0<t<1 .
\]
If the true-null group statistics have this joint representation and the
other $p$-values are arbitrary, ordinary BH satisfies
\eqref{eq:fdx} and \eqref{eq:fdr} with $u=e\pi_0 q$.
\end{theorem}

The required extension of Theorem~\ref{lem:max} concerns the length of a Gaussian projection.

\begin{lemma}[Radial sections]\label{S:lem:radial-sections}
Let $E$ be a positive-dimensional subspace, $\Pi$ its orthogonal
projection, and $K$ a centrally symmetric convex Borel subset of
$\R^d$. There is a version of $\Prob(G\in K\mid\|\Pi G\|=v)$
which is nonincreasing for $v\ge0$.
\end{lemma}

\begin{proof}
Write $G=Z+W$, with independent standard Gaussians on $E$ and
$E^\perp$, respectively. For $z\in E$, set
$h(z)=\Prob(z+W\in K)$. This function is even. It is log-concave
by Pr\'ekopa's marginal theorem applied to
$\one_K(z+w)\exp(-\|w\|^2/2)$ on $E\times E^\perp$
\citep{Prekopa}. If $E^\perp=\{0\}$, then $h=\one_K$.
Along each ray $v\mapsto v\omega$, the function $h(v\omega)$ is
nonincreasing for $v\ge0$ by evenness and log-concavity. The direction
$\Omega$ of $Z$ is uniform on the unit sphere in $E$ and independent
of $\|Z\|$. Hence $v\mapsto\E[h(v\Omega)]$ is the claimed version.
\end{proof}

\begin{proof}[Proof of Theorem~\ref{S:thm:groups}]
For $v>0$ the $\chi_r$ density is proportional to
\[
 v^{r-1}\exp(-v^2/2).
\]
This function is log-concave for every integer $r\ge1$. Integration over $[v,\infty)$ preserves log-concavity by
Pr\'ekopa's theorem. It follows that $\psi_r=-\log S_r$ is increasing and convex. Consequently
\[
   Y_i(g)=\psi_{r_i}(\|\Pi_i g\|)
\]
is even and convex on $\R^d$.
Each $Y_i(G)$ has the standard exponential distribution.

For a centrally symmetric convex set $K$, Lemma~\ref{S:lem:radial-sections}
gives a nonincreasing version $\ell_i(y)=\Prob(G\in K\mid Y_i(G)=y)$,
since $\psi_{r_i}$ is strictly increasing. Tilting Gaussian measure by
$e^b(Y_i-b)_+$ for $b\ge0$ gives $K$ probability
$\E[\ell_i(b+W)]$, where $W$ has density $w e^{-w}$ on $w>0$.
This decreases with $b$. Consequently
$H_b=m^{-1}e^b\sum_i(Y_i-b)_+$ has convex symmetric sublevel sets
and satisfies \eqref{eq:tilt}. The stopped-sum argument in
Theorem~\ref{lem:max} applies, as does the deterministic rank comparison
\eqref{eq:rank}. This proves the global-null bound. Applying it to the
true-null groups in the proof of Theorem~\ref{thm:main} gives the stated
FDP and FDR inequalities.
\end{proof}

\paragraph{Group contrasts with known covariance.}
Suppose $X\sim N_p(\theta,\Sigma)$, and consider prespecified tests of
$L_i\theta=b_i$ with known within-group covariance
$C_i=L_i\Sigma L_i^\top$ of rank $r_i\ge1$. Under its null the Wald
statistic is
\[
 (L_iX-b_i)^\top C_i^+(L_iX-b_i)
    =G^\top B_i^\top(B_iB_i^\top)^+B_i G,
 \qquad B_i=L_i A,\quad AA^\top=\Sigma ,
\]
where $X=\theta+AG$ and $C_i^+$ denotes the Moore--Penrose inverse.
By singular value decomposition, $B_i^\top(B_iB_i^\top)^+B_i$ projects onto the row space of $B_i$. Theorem~\ref{S:thm:groups} therefore applies jointly to the true groups, including overlapping coordinate sets. Only $C_i$ is needed to compute the marginal $p$-value.

\paragraph{Commuting projections.}
Suppose all groups observe $X=\theta+G$ and test $\Pi_j\theta=0$ using $S_{r_j}(\|\Pi_jX\|)$. If each true-null $\Pi_i$ commutes with every $\Pi_j$, then $\FDR\le\pi_0q$. Overlapping coordinate groups in a common orthonormal basis satisfy this condition. Given $\|\Pi_i X\|=v$ at a true null, write $X=v\Omega+X_\perp$. The joint law of the uniform direction $\Omega\in E_i$ and $X_\perp=(I_d-\Pi_i)X$ is independent of $v$. As each $\Pi_j$ preserves $E_i$ and $E_i^\perp$,
\[
 \|\Pi_j X\|^2=v^2\|\Pi_j\Omega\|^2+\|\Pi_j X_\perp\|^2.
\]
Under this coupling all group statistics are nondecreasing in $v$. The resulting $p$-values satisfy positive regression dependence on the true-null subset, giving the assertion by \citet[Theorem~1.2]{BY}.

Related inequalities use total positivity for ordered tests \citep{SarkarMTP} or for distribution functions of certain Wishart-diagonal chi-square laws \citep{RoyenConvolutions}. Theorem~\ref{S:thm:groups} controls rank crossings for arbitrary configurations of subspaces.

\phantomsection


\addcontentsline{toc}{section}{References}
\begin{thebibliography}{99}

\bibitem[Benjamini and Hochberg(1995)]{BH}
Y.~Benjamini and Y.~Hochberg (1995).
Controlling the false discovery rate: a practical and powerful approach to multiple testing.
\emph{Journal of the Royal Statistical Society, Series B} \textbf{57}, 289--300. \msdoi{10.1111/j.2517-6161.1995.tb02031.x}{DOI}.

\bibitem[Benjamini and Yekutieli(2001)]{BY}
Y.~Benjamini and D.~Yekutieli (2001).
The control of the false discovery rate in multiple testing under dependency.
\emph{Annals of Statistics} \textbf{29}, 1165--1188. \msdoi{10.1214/aos/1013699998}{DOI}.

\bibitem[Blanchard and Roquain(2008)]{BR}
G.~Blanchard and E.~Roquain (2008). Two simple sufficient conditions for FDR control.
\emph{Electronic Journal of Statistics} \textbf{2}, 963--992. \msdoi{10.1214/08-EJS180}{DOI}.

\bibitem[Dobriban(2026)]{Dobriban}
E.~Dobriban (2026).
The Benjamini--Hochberg procedure can fail to control the FDR for correlated two-sided Gaussian tests.
\emph{Proceedings of the National Academy of Sciences} \textbf{123}, e2626267123. \msdoi{10.1073/pnas.2626267123}{DOI}.

\bibitem[Genovese et~al.(2006)Genovese, Roeder and Wasserman]{GRW}
C.~R. Genovese, K.~Roeder and L.~Wasserman (2006).
False discovery control with $p$-value weighting. \emph{Biometrika} \textbf{93}, 509--524. \msdoi{10.1093/biomet/93.3.509}{DOI}.

\bibitem[Lei(2026)]{Lei}
L.~Lei (2026).
How much can Gaussian dependence inflate the Benjamini--Hochberg procedure's FDR?
\msarxiv{2607.14812v4}{2607.14812v4}.

\bibitem[Lyu and Wang(2026)]{LyuWang}
C.~Lyu and R.~Wang (2026). On the false discovery rate of one-sided $p$-values under Gaussian dependence.
\msarxiv{2611.19589v1}{2611.19589v1}.

\bibitem[Pr\'ekopa(1973)]{Prekopa}
A.~Pr\'ekopa (1973).
On logarithmic concave measures and functions.
\emph{Acta Scientiarum Mathematicarum (Szeged)} \textbf{34}, 335--343. \href{https://rutcor.rutgers.edu/Prekopa/pdf/SCIENT2.pdf}{Article}.

\bibitem[Roquain and Villers(2011)]{RV}
E.~Roquain and F.~Villers (2011).
Exact calculations for false discovery proportion with application to least favorable configurations.
\emph{Annals of Statistics} \textbf{39}, 584--612. \msdoi{10.1214/10-AOS847}{DOI}.

\bibitem[Royen(2026)]{RoyenConvolutions}
T.~Royen (2026).
A probability inequality for convolutions of $\mathrm{MTP}_2$-distribution functions.
\emph{Statistics \& Probability Letters} \textbf{232}, 110641.
\msdoi{10.1016/j.spl.2026.110641}{DOI}.

\bibitem[Sarkar(1998)]{SarkarMTP}
S.~K. Sarkar (1998).
Some probability inequalities for ordered $\mathrm{MTP}_2$ random variables:
a proof of the Simes conjecture.
\emph{Annals of Statistics} \textbf{26}, 494--504.
\msdoi{10.1214/aos/1028144846}{DOI}.

\bibitem[Shang and Guo(2026)]{ShangGuo}
Z.~Shang and W.~Guo (2026).
A dimension-free bound on the FDR of the Benjamini--Hochberg procedure under slowly decaying Gaussian correlations.
\msarxiv{2611.32898v1}{2611.32898v1}.

\bibitem[Simes(1986)]{Simes}
R.~J. Simes (1986).
An improved Bonferroni procedure for multiple tests of significance.
\emph{Biometrika} \textbf{73}, 751--754. \msdoi{10.1093/biomet/73.3.751}{DOI}.

\bibitem[Su(2018)]{Su}
W.~J. Su (2018).
The FDR-linking theorem.
\msarxiv{1812.08965}{1812.08965}.

\bibitem[Zhang(2026)]{Zhang}
X.~Zhang (2026).
Exact Simes inequalities for two-sided Gaussian tests.
\msarxiv{2611.02513v1}{2611.02513v1}.

\bibitem[Zhang et~al.(2020)Zhang, Tong, Landers and Wu]{TFisher}
H.~Zhang, T.~Tong, J.~Landers and Z.~Wu (2020).
TFisher: a powerful truncation and weighting procedure for combining $p$-values.
\emph{Annals of Applied Statistics} \textbf{14}, 178--201. \msdoi{10.1214/19-AOAS1302}{DOI}.
\end{thebibliography}
\end{document}